\documentclass{amsart}
\usepackage{amsfonts}

\usepackage{amsthm}
\usepackage{eqnarray,amsmath}

\providecommand{\MR}{\relax\ifhmode\unskip\space\fi MR }

\usepackage{tikz-cd}
\usepackage{amssymb}
\usepackage{float}
\usepackage{microtype}
\usepackage[shortlabels]{enumitem}
\usepackage{tikz}
\usetikzlibrary{arrows.meta}
\usepackage{comment}
\usepackage{mathrsfs}
\usepackage{mathtools}
\usepackage{hyperref}
\usepackage{bbm}
\usepackage{float}
\usepackage[all]{xy}

\hypersetup{colorlinks,linkcolor={blue},citecolor={olive},urlcolor={red}}

\numberwithin{equation}{section}

\usepackage[top=1.3in, left=1.3in, right=1.3in, bottom=1.3in]{geometry}

\usepackage[normalem]{ulem}

\newcommand{\E}{\mathbb{E}}
\newcommand{\Var}{\mathop{\mathrm{Var}}}

\renewcommand{\P}{\mathbb{P}}

\newcommand{\one}{\mathbf{1}}

\newcommand{\of}[1]{\left(#1\right)}

\newcommand{\Z}{\mathbb{Z}}

\newcommand{\R}{\mathbb{R}}

\newcommand{\bZ}{{\mathbb Z}}

\DeclarePairedDelimiter\floor{\lfloor}{\rfloor}

\newtheorem{innercustomlemma}{Lemma}
\newenvironment{customlemma}[1]
{\renewcommand\theinnercustomlemma{#1}\innercustomlemma}
{\endinnercustomlemma}

\usepackage{theoremref}

\def\R{\mathbf{r}}

\newtheorem{theorem}{Theorem}
\newtheorem{lemma}{Lemma}[section]
\newtheorem{prop}[lemma]{Proposition}
\newtheorem{corollary}[lemma]{Corollary}

\theoremstyle{definition}
\newtheorem{remark}[lemma]{Remark}
\newtheorem{define}[lemma]{Definition}

\usepackage[font=small]{caption}

\makeatletter
\def\thmref@flush{%
	\ifx\thmref@last\empty\else
	\ifthmref@comma, \thmref@finaltrue\fi \thmref@commatrue
	\thmref@last \ifx\thmref@stack\empty\else s\fi \thmref@num 0
	\let\do\thmref@one \thmref@stack
	\ifcase\thmref@num\or\space and\else\thmref@finaltrue, and\fi
	~\ref{\thmref@head}\let\thmref@stack\empty\fi}
\def\thmref@one#1{\ifnum\thmref@num>0,\fi
	\space\ref{#1}\advance\thmref@num 1\relax}
\makeatother

\definecolor{hancolor}{rgb}{0.0 0.0, 1.0}
\newcommand{\commHL}[1]{{\textcolor{hancolor}{#1$_{-[HL]}$}}} 

\usepackage{graphicx}
\graphicspath{ {./images/} }
\usepackage{graphicx}
\usepackage{caption}
\usepackage{subcaption}
\usepackage{mwe}

\title[One-dimensional excitable media with variable interaction range]{Phase transition in one-dimensional excitable media \\ with variable interaction range}

\author[A.~Aguirre]{Ander Aguirre}
\address{Ander Aguirre,  Department of Mathematics, The Ohio State University, Columbus, OH 43210}
\email{\texttt{aguirre.93@osu.edu}}

\author[H.~Lyu]{Hanbaek Lyu}
\address{Hanbaek Lyu, Department of Mathematics, University of Wisconsin -- Madison, Madison, WI 53706}
\email{\texttt{hlyu@math.wisc.edu}}

\author[D.~Sivakoff]{David Sivakoff}
\address{David Sivakoff, Department of Statistics and Department of Mathematics, The Ohio State University, Columbus, OH 43210}
\email{\texttt{dsivakoff@stat.osu.edu}}

\thanks{}

\begin{document}

	\begin{abstract}
		We investigate two discrete models of excitable media on a one-dimensional integer lattice $\mathbb{Z}$: the $\kappa$-color Cyclic Cellular Automaton (CCA) and the $\kappa$-color Firefly Cellular Automaton (FCA). In both models, sites are assigned uniformly random colors from $\mathbb{Z}/\kappa\mathbb{Z}$. Neighboring sites with colors within a specified interaction range $r$ tend to synchronize their colors upon a particular local event of 'excitation'. We establish that there are three phases of CCA/FCA on $\mathbb{Z}$ as we vary the interaction range $r$.
		First, if $r$ is too small (undercoupled), there are too many non-interacting pairs of colors, and the whole graph $\mathbb{Z}$ will be partitioned into non-interacting intervals of sites with no excitation within each interval. If $r$ is within a sweet spot (critical), then we show the system clusters into ever-growing monochromatic intervals. For the critical interaction range $r=\lfloor \kappa/2 \rfloor$, we show the density of edges of differing colors at time $t$ is $\sim ct^{-1/2}$ and each site excites about $t^{1/2}$ times up to time $t$. Lastly, if $r$ is too large (overcoupled), then neighboring sites can excite each other and such 'defects' will generate waves of excitation at a constant rate so that each site will get excited at least at a linear rate. For the special case of FCA with $r=\lfloor 2/\kappa \rfloor+1$, we show that every site will become $(\kappa+1)$-periodic eventually. 
	\end{abstract}

	\maketitle


	\section{Introduction}

	An excitable medium is a network of coupled dynamic units whose states get excited upon a particular local event. It has the capacity to propagate waves of excitation, which often self-organize into spiral patterns. Examples of such systems in nature include neural networks, Belousov-Zhabotinsky reaction,  as well as coupled oscillators such as fireflies and pacemaker cells. In a discrete setting, excitable media can be modeled using the framework of generalized cellular automaton (GCA). Given a simple connected graph $G=(V,E)$ and a fixed integer $\kappa\ge 2$,  the microstate of the system at a given discrete time $t\ge 0$ is given by a $\kappa$-coloring of vertices $X_{t}:V\rightarrow\mathbb{Z}_{\kappa}=\mathbb{Z}/\kappa\mathbb{Z}$. A given initial coloring $X_{0}$ evolves in discrete time via iterating a fixed deterministic transition map $\tau: X_{t}\mapsto X_{t+1}$, which depends only on local information at each time step. That is, for each $v\in V$, $X_{t+1}(v)$ is determined by $X_{t}$ restricted on $N(v)\cup\{v\}$, where $N(v)$ is the set of neighbors of $v$ in $G$. This generates a trajectory $(X_{t})_{t\ge 0}$, and its  limiting behavior in relation to the topology of $G$ and structure of $\tau$ is of our interest.

	Greenberg-Hastings Model (GHM) and Cyclic cellular automaton (CCA) are two particular discrete excitable media which have been studied extensively since the 90s. GHM was introduced by Greenberg and Hastings \cite{greenberg1978spatial} to capture the phenomenological essence of neural networks in a discrete setting, whereas CCA was introduced by Bramson and Griffeath \cite{bramson1989flux} as a discrete time analogue of the cyclic particle systems. In GHM, think of each vertex of a given graph as a $\kappa$-state neuron. An excited neuron (i.e., one in state 1) excites neighboring neurons at rest (in state 0) and then needs to wait for a refractory period of time (modeled by the remaining $\kappa-2$ states) to become rested again. In CCA, each vertex of the graph is inhabited by one of $\kappa$ different species in a cyclic food chain. Species of color i are eaten (and thus replaced) by species of color $(i+1)$ mod $\kappa$ in their neighborhood at each time step. More precisely, denoting the $\kappa$-color CCA by $\eta_{t}$, its time evolution is given by 
		\begin{equation}
			\eta_{t+1}(v)=\begin{cases}
				\eta_{t}(v)+1 \,\, (\text{mod $\kappa$})& \text{if $\exists u\in N(v)$ s.t. $\eta_{t}(u)=\eta_{t}(v)+1$ (mod $\kappa$) }\\
				\eta_{t}(v) & \text{otherwise}
			\end{cases}.
		\end{equation}

		In \cite{lyu2015synchronization}, the second author proposed a discrete model for coupled oscillators, and studied criteria for synchronization on some classes of finite graphs. The basic setup is the same as  CCA or GHM. Fix an integer $\kappa\ge 3$ and let $\bZ_{\kappa}=\bZ/\kappa\bZ$ with linear ordering $0<1<2<\cdots<\kappa-1$ be the color space. 
		Consider each vertex as a $\kappa$-state firefly, which blinks whenever it has color $0$. During each iteration, post-blinking fireflies (whose color is in $\{1,\dots,\lfloor \kappa/2 \rfloor\}$) with a blinking neighbor (a neighbor of color $0$) stay at the same color, and all others increase their colors by 1 modulo $\kappa$. The discrete dynamical system $(X_{t})_{t\ge 0}$ generated by the iteration of the above transition map is called the $\kappa$-color \textit{firefly cellular automaton (FCA)} on $G$. More precisely, denoting by $X_{t}$ the $\kappa$-color FCA, its time evolution is given by 
		\begin{align}
			X_{t+1}(x)
			&=\begin{cases}
				X_t(x) & \text{if $\exists w\in N(x)$ with  $X_t(w)=0$ and $X_t(x)\in [1,\lfloor \kappa/2\rfloor]$}\\
				X_t(x)+1 (\text{mod $\kappa$}) & \text{otherwise}\nonumber\end{cases}
		\end{align}

		Both CCA and GHM dynamics have been extensively studied on integer lattices $G=\mathbb{Z}^{d}$ using probabilistic methods, where one takes the initial configuration $X_{0}$ as a random $\kappa$-coloring on sites according to the uniform product probability measure $\mathbb{P}$ on $(\mathbb{Z}_{\kappa})^{\mathbb{Z}^{d}}$. We introduce some terminology for FCA to describe its behavior, which may be defined for CCA and GHM similarly. We say that $X_{t}$ \textit{fixates} if every site is excited only finitely many times $\mathbb{P}$-a.s., and \textit{synchronizes} if for every two vertices $x,y\in V$, there exists $N=N(x,y)\in \mathbb{N}$ such that $X_{t}(x)=X_{t}(y)$ for all $t\ge N$ $\mathbb{P}$-a.s.. It is not hard to see that fixation and synchronization are equivalent notions if and only if any initial coloring on the complete graph with two vertices synchronize, which is the case for GHM and FCA for all $\kappa\ge 3$ and CCA only for $\kappa=3$: CCA for $\kappa\ge 4$ has a pair of distinct but non-interacting colors, resulting in fixation without synchronization.

		In this paper, we study an extended family of CCA and FCA on $\mathbb{Z}$ where the `interaction range', denoted by an integer parameter $r$, is also a variable. Namely, in CCA, note that only nearby colors in the color wheel $\mathbb{Z}_{\kappa}$ interact; in FCA, only colors within distance $\le \lfloor \kappa/2\rfloor$ can interact. Thus, CCA and FCA are now parameterized by the total number $\kappa$ of available colors and the range $r$ of interaction. The extended models read as follows:
		\begin{align}
			\textup{(CCA)}\quad \, \eta_{t+1}(x)
			&=\begin{cases}
				\eta_t(x)+1 \, (\text{mod $\kappa$}) & \text{if $\exists w\in N(x)$ s.t.  $\eta_t(w)\in \{\eta_t(x)+1,\dots, \eta_t(x)+r\}$}\\
				\eta_t(x) & \text{otherwise}\end{cases}
			\label{CCA}  \\
			\textup{(FCA)}\quad X_{t+1}(x)
			&=\begin{cases}
				X_t(x) & \text{if $\exists w\in N(x)$ s.t.  $X_t(w)=0$ and $X_t(x)\in [1,r]$}\\
				X_t(x)+1 \, (\text{mod $\kappa$}) & \text{otherwise}\nonumber\end{cases}
			\label{FCA}	
		\end{align}
		In this paper, we will use $\zeta_{t}$ to denote any one of the two models above on the one-dimensional integer lattice $\mathbb{Z}$ with the initial configuration $\xi_{0}$ distributed according to the uniform product measure on the set $(\mathbb{Z}_{\kappa}^{\Z})$ of all possible $\kappa$-colorings of $\mathbb{Z}$.

			The central notion in the dynamics of the above three discrete models for excitable media is excitation. We say a site $x$ is \textit{excited} at time $t$ if its internal dynamics are affected by its neighbors at time $t$. That is, if $\eta_{t+1}(x)=\eta_{t}(x)+1 (\textup{mod}\, \kappa)$, or if $X_{t+1}(x)=X_{t}(x)$ for FCA. Note that excitation always comes from local disagreements and the sites excite their neighbors to remedy the current disagreement. It is the non-linear aggregation of this mutual effort to synchronize with neighbors that makes studying the global dynamics of excitable media interesting.  
			In particular, one of the main questions in the study of excitable media is if such local efforts for synchronization would actually lead into larger-scale synchronization of the phases. To capture such a phenomenon, we say the system \textit{clusters} if it is overwhelmingly likely to see a single color on any finite subset of sites after a long time, that is, 
			\begin{align}
				\lim_{t\rightarrow \infty} \mathbb{P}(\xi_{t}(x)\ne \xi_{t}(y))=0 
			\end{align}
			for any two sites $x$ and $y$. By union bound, clustering implies that all sites in any finite subset of sites have the same color in probability as $t\rightarrow\infty$. It is known that clustering occurs for the 3-color GHM \cite{durrett1991some} as well as the 3-color critical CCA \cite{fisch1992clustering} and FCA \cite{lyu2017persistence}. 
			
			However, not all excitable media cluster. In fact, we find that clustering is a feature that characterizes `critical excitable media' for which the interaction range $r$ is not too small and not too large. That is, if $r$ is too small, then initially there are too many edges with non-interacting colors at their ends (i.e., color differences $>r$), a positive fraction of which persists and prevent clustering. The classical example for excitable media that does not cluster is the original CCA with $r=1$ for $\kappa\ge 5$. These systems are known to `fixate', where each site excites only finitely many times and hence converges to some limiting color \cite{fisch1990cyclic}. More precisely, for each site $x$ and time $t$, define 
			\begin{equation}
				\mathcal{E}_{t}(x):=\sum_{s=1}^{t} \one(\text{$x$ is excited at time $s$}),
			\end{equation} 
			which is the number of times that a site $x$ excites through time $t$. Since the process $\zeta_{t}$ is translation invariant, the distribution of the random variable above does not depend on $x$. The basic dichotomy of the system is given by the asymptotic behavior of the such excitation counting function. Namely, we say the system $\zeta_{t}$ \textit{fixates} if $\mathcal{E}_t(0)$ is finite a.s., i.e. each the origin is only excited only finitely many times with probability 1 and \textit{fluctuates} otherwise.

			
			

			
			\begin{figure}[h!]
				\centering
				\includegraphics[width=1\textwidth]{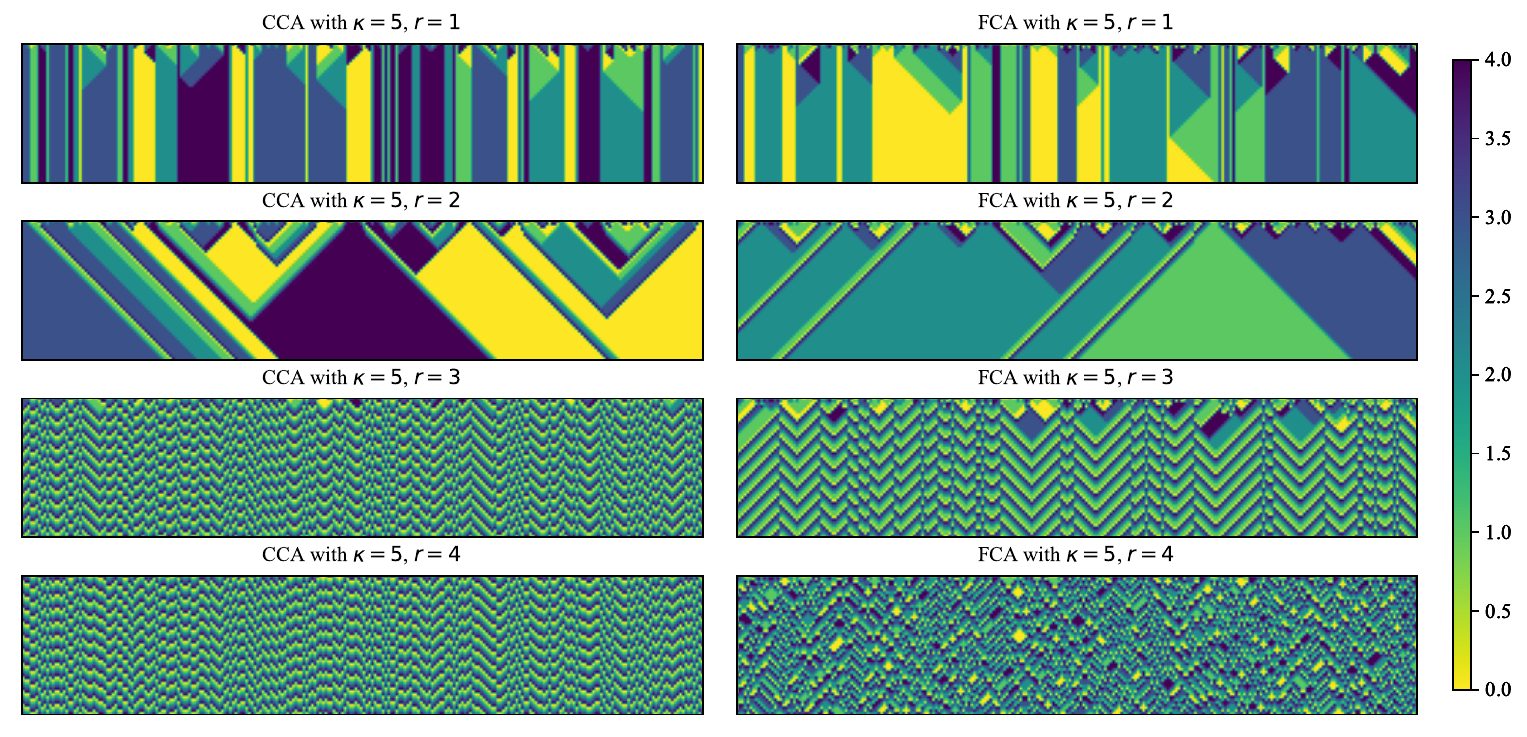}
				\caption{ Simulation of CCA and FCA with $\kappa=5$ and $r\in \{1,2,3,4\}$ on 500 nodes in $\mathbb{Z}$ with 50 and 250 iterations, respectively. Time goes from top to bottom. For CCA configurations at every iteration are shown, while for FCA only the ones at times $\kappa t$ are shown.
				}
				\label{fig:CCA_FCA_sim1}
			\end{figure}

			\begin{figure}[h!]
				\centering
				\includegraphics[width=1\textwidth]{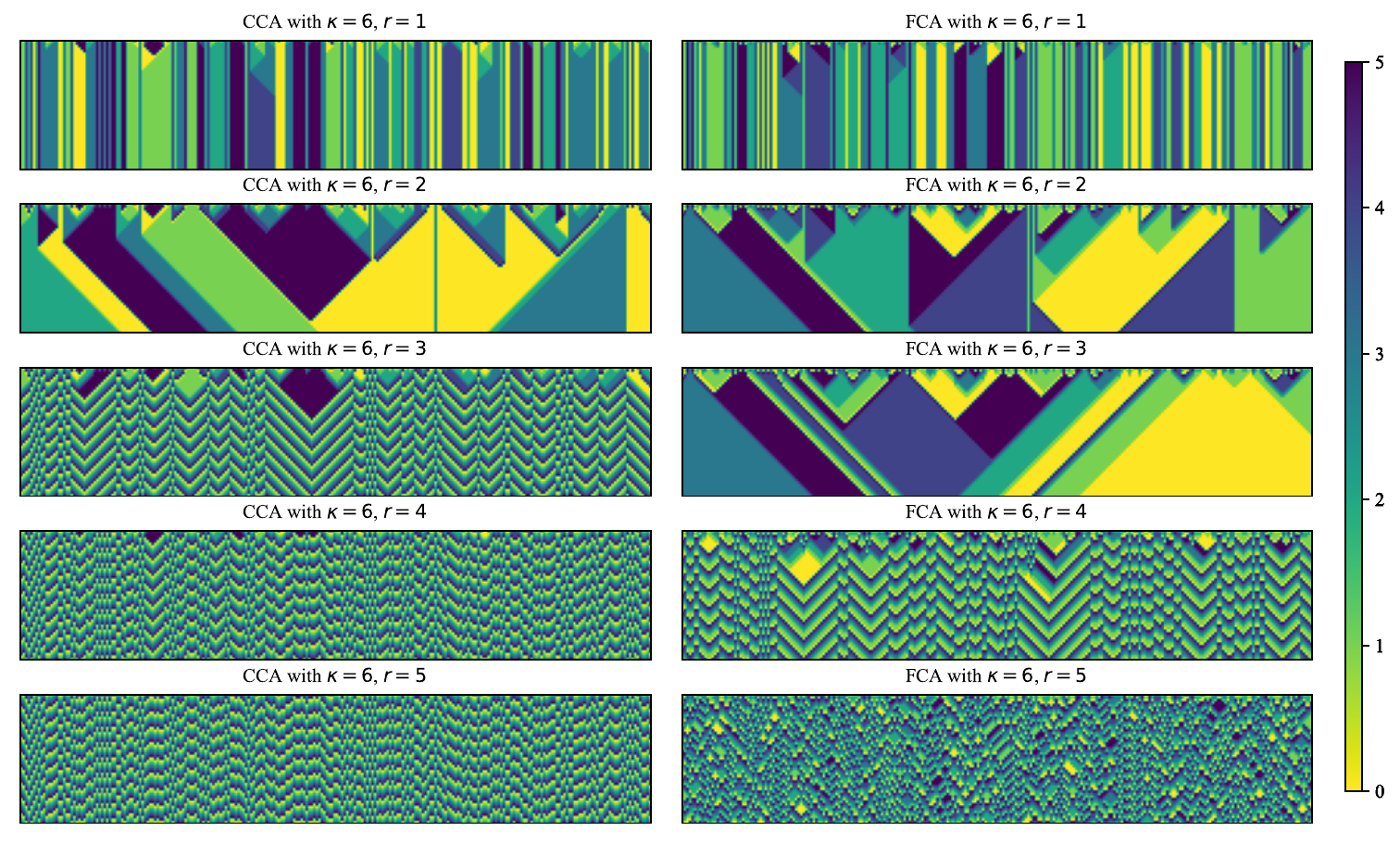}
				\caption{ Simulation of CCA and FCA with $\kappa=6$ and $r\in \{1,2,3,4,5\}$ on 500 nodes in $\mathbb{Z}$ with 50 and 250 iterations, respectively. Time goes from top to bottom. For CCA configurations at every iteration are shown, while for FCA only the ones at times $\kappa t$ are shown.
				}
				\label{fig:CCA_FCA_sim2}
			\end{figure}
			
			Clustering may also not occur when $r$ is too large. In this case, two neighboring sites can excite each other back and forth, and such `defects' will generate waves of excitation at a constant rate. The whole graph $\mathbb{Z}$ will then be partitioned into `cells' that are under the influence of each defect, and each site will get excited at least at a linear rate. This phenomenon was observed for CCA in two-dimension \cite{fisch1991cyclic}. There, a closed loop where each node initially is next to a neighbor of one color larger appear with a positive density at time zero. Each such object, called a stable periodic object (SPO), act as a defect and each constituent site gets excited at a maximal rate (i.e., every single time) and hence it is stable as it cannot be perturbed by any external configuration. SPOs drive the entire graph $\Z^{2}$ into disjoint cells of linear excitation. We find a pair of adjacent sites can act as as SPO in the one-dimensional lattice when $r$ is too large. It is evident that such defects can prevent the system from clustering.


			
			Our main result in this paper makes the three-phase picture we discussed above rigorous. That is, we establish that there are three phases of CCA/FCA on $\mathbb{Z}$ as we vary the interaction range $r$. First, if $r$ is too small (`undercoupled regime'), then there are too many non-interacting pairs of colors and the whole graph $\Z$ will be partitioned into non-interacting intervals of sites and there will be no excitation in each intervals (fixation). If $r$ is too large (`overcoupled regime'), then two neighboring sites can excite each other, and such `defects' will generate waves of excitation at a constant rate. The whole graph $\mathbb{Z}$ will then be partitioned into `cells' that are under the influence of each defect, and each site will get excited at least at a linear rate. Lastly, if $r$ is somewhere in the `sweet spot', including the critical value $\lfloor \kappa/2 \rfloor$, local efforts of resolving immediate phase disagreements can be resolved globally and the system will cluster. Theorem \ref{thm:main1} below summarizes these results.

			

			\begin{theorem}[Phase transition of CCA/FCA in interaction range]
				\label{thm:main1}
				Let $\xi_{t}$ denote the $\kappa$-color CCA or FCA  with interaction range $r$ on the one-dimensional integer lattice $\Z$ for $\kappa\ge 3$, where the initial configuration $\xi_{0}$ is drawn from the uniform product measure $\mathbb{P}$ on $\bZ_\kappa^{\mathbb{Z}}$. The following holds: 
				
				\begin{description}[itemsep=0.1cm]
					\item[(i)] (\textit{Undercoupled regime}) For $r<\lfloor  (\frac{2-\sqrt{2}}{4}) \kappa \rfloor$, every site fixates 
					and clustering does not occur. 
					
					\item[(ii)] (\textit{Critical regime}) Suppose $r=\lfloor \frac{\kappa-1}{2} \rfloor$ for CCA and $r=\lfloor \kappa/2 \rfloor$ for FCA. Then  $\xi_t$ fluctuates and clusters. Furthermore, if $r=\lfloor \kappa/2 \rfloor$ except the case of CCA with $\kappa$ even, 
					\begin{align}\label{eq:clustering_rate}
						\hspace{-2cm} \textup{(Clustering rate)} \qquad    \P(\xi_{t}(0)\ne \xi_{t}(1)) &= \Theta(t^{-1/2}).
					\end{align}

					\item[(iii)] (\textit{Overcoupled regime}) If $\lfloor k/2\rfloor< r<\kappa-1$, then $\xi_t$ fluctuates and $\liminf_{t\rightarrow\infty} \frac{1}{t}\mathcal{E}_{t}(0)>0$. Furthermore,
					\begin{description}[itemsep=0.1cm]
						\item[(a)] For CCA, for all $(\kappa,r)$ such that either 
						$\lfloor k/2\rfloor< r\le \kappa-1$ or $(\kappa,r)=(\textup{even}, \lfloor \kappa/2 \rfloor)$,  every site eventually excites at every iteration. 
						
						\item[(b)] For FCA, there exists a positive density of sites that are eventually $(\kappa+1)$-periodic and no sites can be eventually $\kappa$- or $(\kappa+2)$-periodic. 
						Moreover, if $\kappa \ge 5$ and $r=\lfloor k/2\rfloor+1$, then every site is eventually $(\kappa+1)$-periodic. 
					\end{description}

				\end{description}	
			\end{theorem}

			Next, we obtain asymptotics for the excitation count $\mathcal{E}_{t}(0)$ for some critical CCA and FCA. We show that it scales as $\sqrt{t}$ and identify bounds on its scaling limit. Our result is sharp for the three color systems.

			\begin{theorem}[Excitation rate in critical CCA/FCA]
				\label{thm:main2}
				Let $\xi_{t}$ denote one of the following  critical systems on $\Z$: FCA with $r=\lfloor \kappa/2 \rfloor$ for both $\kappa$ even and odd; and CCA with $\kappa$ even and $r= \kappa/2$. Assume that the initial configuration $\xi_{0}$ is drawn from the uniform product measure $\mathbb{P}$ on $\bZ_\kappa^{\mathbb{Z}}$.  Then the following hold:
				\begin{description}
					\item[(i)] The expected excitation count satisfies 
					\begin{align}\label{eq:excitation_rate}
						\E[\mathcal{E}_t(0)] &= \Theta(t^{1/2}).
					\end{align}
					
					\item[(ii)] Let $\mathtt{M}$ denote a nonnegative random variable satisfying $\P(\mathtt{M}\ge s) = 4\mathbb{P}(Z\ge s) \mathbb{P}(Z\le s)$ for $s\ge 0$, where  $Z\sim N(0,1)$ is a standard normal random variable. Then there exist explicit constants $c,C>0$ such that ($\preceq$ denoting stochastic domination)
					\begin{equation}\label{eq:excitation_rate_gen}
						c\, \mathtt{M}  \preceq  \liminf_{t\rightarrow\infty} \frac{\mathcal{E}_{t}(x)}{ \sqrt{t}}	 \le  \limsup_{t\rightarrow\infty} \frac{\mathcal{E}_{t}(x)}{ \sqrt{t}}	\preceq  C \, \mathtt{M}.
					\end{equation}
					Furthermore, suppose $\kappa=3$ and let $\sigma:=\sqrt{\frac{2}{3}}$ for CCA and $\sqrt{\frac{8}{81}}$ for FCA. Then 
					\begin{equation}\label{eq:3FCA_exact_tournament_eq2}
						\lim_{t\rightarrow\infty} \frac{\mathcal{E}_{t}(x)}{\sigma \sqrt{t}}	\overset{d}{=} \mathtt{M}. 
					\end{equation}
				\end{description}
			\end{theorem}

			The clustering results in Theorem \ref{thm:main1} for the three-color critical regime ($\kappa=3$ and $r=1$) is known due to Fisch \cite{fisch1992clustering} for CCA and Lyu and Sivakoff \cite{lyu2017persistence} for FCA. All other results are new  to this work. Especially, our analysis for FCA in the critical and overcoupled regimes brings several new technical innovations to the literature. 


			\subsection{Related work and discussion}
			\label{sec:related_work}
			Here we briefly summarize known results for discrete excitable media (namely, GHM, CCA, and FCA) in 1-dimension, and some of the main proof techniques. We then discuss our new analysis techniques developed in this work. 
			
			\subsubsection{Fixation vs. Fluctuation.} 
			
			Most results on 1-dimensional models rely on a particle systems analogy where we place ``edge particles" on the boundaries between distinctly colored regions and consider their time evolution. By counting ``live edge'' particles against ``blockade'' particles and with some careful arguments, Fisch \cite{fisch1990one} showed that $\kappa$-color CCA on $\mathbb{Z}$ fixates if and only if $\kappa\ge 5$. We use a similar technique to establish the undercoupled regime (Theorem \ref{thm:main1} \textbf{(i)}). 
			
			\subsubsection{Clustering and persistence of particle counting walks}   
			
			In 1992, Fisch \cite{fisch1992clustering} showed that the critical 3-color CCA on $\mathbb{Z}$ clusters with an exact asymptotic on the clustering rate 
			\begin{equation}\label{eq:3CCA_clustering}
				\mathbb{P}(\xi_{t}(x)\ne \xi_{t}(x+1))\sim \sqrt{\frac{2}{3\pi}}t^{-1/2}. 
			\end{equation}
			Using an embedded particle system, Fisch related the probability of local disagreement at time $t$ as the probability of a certain associated particle counting random walk staying nonnegative for $2t$ steps. Sharp asymptotics for such `persistence probability' of random walks are well known in the literature (see, e.g., Sparre Anderson \cite{andersen1953sums} and Feller \cite[Thm. XII.7.1]{feller1971introduction}).

			The embedded particle system description of the 3-color CCA  is as follows. At time 0, place a right  or left arrow on each edge independently with probability 1/3. Right arrows move to right with constant unit speed and left arrows behave similarly; if opposing particles ever collide or have to occupy the same edge, they annihilate each other and disappear. Now if there is a right arrow on the edge $(0,1)$ at time $t$, this particle must have been on the edge $(-t,-t+1)$ at time $0$ and must travel distance $t$ without being annihilated by an opposing particle. This event is determined by the net counts of right versus left arrows at time 0 starting from the edge $(-t,-t+1)$ and moving rightward. Namely, the excess number of right arrows on successive intervals $[-t,-t+s]$, $1\le s \le 2t+1$, form a random walk. The right arrow moves as long as this random walk survives (stays at positive height), which is an event of probability $\Theta(t^{-1/2})$. A similar technique was incorporated by Durrett and Steif in \cite{durrett1991some} to show similar clustering results for GHM on $\mathbb{Z}$ for $\kappa=3$: The same asymptotic \eqref{eq:3CCA_clustering} with the constant $\sqrt{2/(3\pi)}$ replaced with 
			$\sqrt{2/(27\pi)}$. Later this result was extended to arbitrary $\kappa\ge 3$ by Fisch and Gravner in \cite{fisch1995one}. Clustering for the 4-color critical ($r=1$) CCA was only recently shown by Hellouin de Menibus and Borgne \cite{hellouin2021asymptotic} in 2021. However, the exact clustering rate is still unknown and simulation indicates that the mean cluster size of such system grows at a rate different from $t^{1/2}$  as in the 3-color systems \cite{fisch1992clustering}. This is due to the existence of `blockades' (i.e., edges of non-interacting colors), which can flip the direction of arrows upon collision. (See `stack flipping' in Section \ref{sec:arrow_dynamics} for more details). 
			
			\subsubsection{Particle flipping and random speeds}
			FCA shares a similar embedded annihilating particle system structure, but with additional arrow flipping phenomena at time 0 without blockades. The initial site coloring at time $0$ induces a canonical assignment of edge particles. In the FCA dynamics, the system takes a finite amount of ``burn-in'' period, during which particles may flip their directions and thereafter they stabilize and move in only one direction with annihilation upon collision. This finite burn-in period introduces dependencies between edge particles, so the associated random walk has correlated increments. For example, consider a 3-configuration $\cdots 012 \cdots $ on $\mathbb{Z}$, which corresponds to two consecutive right arrows. Applying the 3-color FCA rule, it evolves to $\cdots 110\cdots $, which has one left particle between 1 and 0, as if the right arrow between 0 to 1 flips the right arrow to its right between 1 and 2. A similar phenomenon occurs for all $\kappa\ge 3$, so an associated random walk has correlated increments. 
			Thus, in order to obtain sharp clustering rate for such systems, one needs to know the sharp asymptotics of the probability that partial sums of correlated increments staying nonnegative. With this motivation, Lyu and Sivakoff \cite{lyu2017persistence} established such a result for Markov additive functionals, which are partial sums where the increments are functionals of an underlying Markov chain. This setting is flexible enough to be applied to the analysis of embedded particle systems for one-dimensional cellular automata, including the critical FCA. Using this general result, Lyu and Sivakoff showed that 3-color critical FCA clusters with rate given by  \eqref{eq:3CCA_clustering}, now with the constant $\sqrt{2/(3\pi)}$ replaced by  
			$\sqrt{8/(9\pi)}$.
			
			Clustering rate for general critical FCA for $\kappa\ge 3$ was not known before. The main issue there was the for larger $\kappa$, there could be multiple arrows on each edge, and depending on that, a given arrow may move with random speed. The same issue occurs with critical CCA with large $\kappa$ (see Figure \ref{fig:9CCA_FCA}). It is relatively easy to uniformly upper and lower bound on the speed of a given particle. However, then due to the ambiguity of particle speed, there could be at least order $t$ many different edges that can potentially send an arrow that is on the edge $(0,1)$ at time $t$, and using union bound on all such candidate edges will ruin the sharp asymptotic on the clustering rate. In this work, we circumvent this issue by using the mass transport principle (Prop. \ref{prop:MTP}). This allows us to convert the event of having an arrow at an edge at time $t$ to the event that a given particle survives until time $t$, thereby removing the need to use a union bound. Mass transport principle has been used by Lyu and Junge to analyize ballistic annihilation \cite{junge2018phase} to relate the probabilities of different collision events.

			\subsubsection{Excitation counts and tournamant expansion}
			
			In Theorem \ref{thm:main2}, we obtain asymptotics for the excitation count $\mathcal{E}_{t}(0)$ for some critical CCA ($\kappa$ odd and $r=\lfloor \kappa/2\rfloor$) and FCA ($r=\lfloor \kappa/2 \rfloor$). 
			The analysis of the excitation count is quite different from the analysis of clustering rate, since one needs to know how many excitation that the origin will get during a time period instead of an arrow surviving for an extended period of time.
			
			For this purpose, we take a novel approach of constructing a monotone comparison process, which was first developed by Gravner, Lyu, and Sivakoff  \cite{gravner2016limiting} to study the 3-color CCA and GHM on arbitrary graphs. We develop a similar technique for critical CCA and FCA on $\mathbb{Z}$ with arbitrary $\kappa$, by which we are able relate the maximum of an associated particle-counting walk to the number of excitations of the origin. We then apply the sharp asymptotic of the expected maxima of Markov additive functionals in \cite{lyu2017persistence} and a functional central limit theorem for Markov additive functionals  \cite[Ch. 17]{meyn2012markov} to deduce Theorem \ref{thm:main2}.

				\subsubsection{Quasi-Stable Periodic Objects in the overcoupled regime}

				In the overcoupled regimes, to sites can excite each other back and forth, so there might be an `echo chamber' of excitation that persist and keep influencing nearby sites. Hence, we expect a roughly linear scaling for the excitation count $\mathcal{E}_t(0)$. It is also worth presenting a more fine-grained view of the behaviour around criticality, where every site starts to become influenced by its neighbors infinitely often. 
				
				For the overcoupled CCA, indeed, any two sites on an edge with color difference $>\lfloor \kappa/2 \rfloor$ excite each other indefinitely, and such local dynamics with internally supported maximal excitation cannot be perturbed by any external dynamics. This is an example of \textit{Stable Periodic Orbits} (SPOs), which was first used by Fisch, Gravner, and Griffeath  to analyze two-dimensional CCA \cite{fisch1991cyclic}. From this, it is easy to deduce that every site will get excited every iteration eventually.  
				
				However, the analysis of overcoupled FCA is significantly more delicate. Since in FCA each site excites its neighbors only when in color 0, such SPOs cannot exist. In FCA on $\Z$, each site can have at most two excitation between consecutive returns to 0, so the lengths between consecutive 0 can be any of $\kappa, \kappa+1$, and $\kappa+2$. In principle, such `blinking gaps' can occur in any combination, so the dynamics of overcoupled FCA can still be very complicated. 
				
				Our analysis of overcoupled FCA proceeds as follows. We first identify a small configuration (color string $rrr0rrr$) of length 7 and show that the dynamics of the internal sites is $(\kappa+1)$-periodic and cannot be perturbed by external configuration, whereas the boundary sites can still be perturbed. While this could be one of many such quasi-SPOs --- a periodic object surrounded by some protective layers --- their existence is enough to deduce that a positive density of sites will be $(\kappa+1)$-periodic in overcoupled FCA. We then show that every single site must become eventually $(\kappa+1)$-periodic when the system is `barely overcoupled' with $r=\lfloor \kappa/2 \rfloor+1$. The main idea is to analyze `defects', which are edges with large color difference that cannot be created spontaneously. We classify all possible local dynamics of such defects (similar local analysis was used to analyze FCA on finite trees \cite{lyu2016time}), and show that defects can be recurrent only if the sites in them are eventually $(\kappa+1)$-periodic. These recurrent defects generate one arrow every $\kappa+1$ iterations, and they travel into intervals of non-defect edges without flipping and collide with opposing arrows. From this we deduce $(\kappa+1)$-periodicity spreads into such intervals of non-defect edges.

				\vspace{9mm}

				\subsection{Organization} In Section \ref{sec:arrow_dynamics},
				we present a dual description for the time evolution of CCA/FCA through their so-called \textit{ arrow dynamics}. Section \ref{sec:CCA_undercoupled} discusses the \textit{undercoupled regime} $(r<\lfloor\frac{\kappa}{2} \rfloor)$ for both models and contains the proof of Theorem \ref{thm:main1} \textbf{(i)}. Section \ref{sec:critical} is devoted to the \textit{critical regime}. We prove the clustering of Theorem \ref{thm:main1} \textbf{(ii)}. In there, we also study clustering (Section \ref{sec:critical_2}) and excitation rates at criticality and thus have the proof of Theorem \ref{thm:main2} (Section \ref{sec:critical_3}). In Section \ref{sec:overcoupled_regime}, we deal with the \textit{overcoupled regime} and the proof of Theorem \ref{thm:main1} \textbf{(iii)}.  For the FCA in the generic overcoupled case, the existence of persistent quasi-stable-periodic-objects ensures the excitation rate of all other sites is linear albeit random (Section \ref{sec:overcoupled_regime_gen}).  In the particular, \textit{weakly overcoupled}, case ($r=\lceil\frac{\kappa}{2}\rceil+1$), we are indeed able to show all lattice sites are eventually $(\kappa+1)$-periodic (Section \ref{sec:FCA_weakly_overcoupled}).


				

				\section{Arrow Dynamics for undercoupled and critical CCA/FCA}
				\label{sec:arrow_dynamics}
				
				In this section, we introduce embedded arrow dynamics for undercoupled and critical CCA/FCA. That is, we will only consider FCA with $r\le \lfloor \kappa/2\rfloor$ and CCA with $r<\lfloor \kappa/2 \rfloor$ or $\kappa$ odd and $r=\lfloor \kappa/2 \rfloor$. In these cases, there is no spontaneous emergence of excitation so the flow of excitation can be understood via certain annihilating particle system. For the overcoupled systems, there could be `echo chambers' of excitation. See Section \ref{sec:overcoupled_regime} for more discussion.    
				
				\subsection{Arrow dynamics for the 3-color critical CCA/FCA} 
\begin{figure}[h!]
					\centering
					\includegraphics[width=1\textwidth]{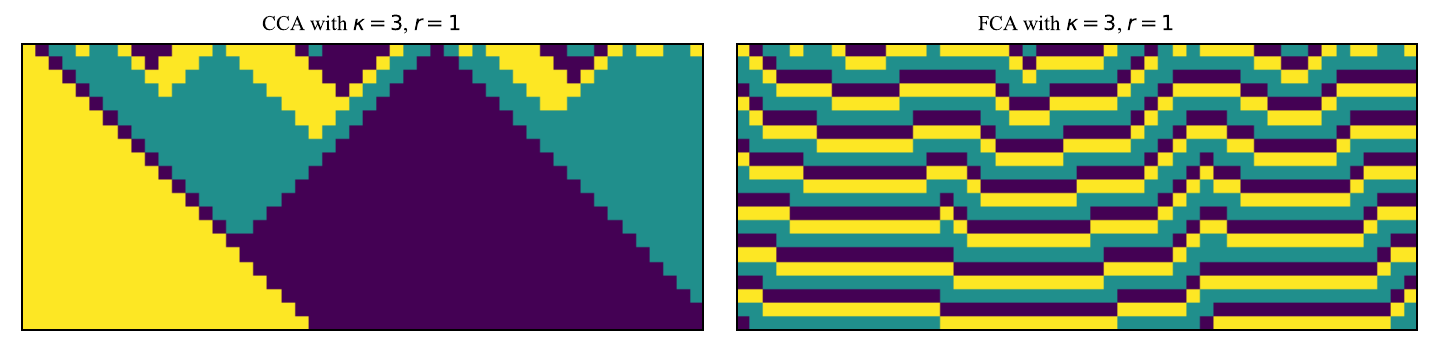}
					\vspace{-0.5cm}
					\caption{ Simulation of 3-color critical CCA and FCA ($r=1$) on 50 nodes in $\mathbb{Z}$ with 20 iterations. Time goes from top to bottom where configurations at every iteration are shown for both systems. 
					}
					\label{fig:3CCA_FCA}
				\end{figure}

One of the classical approaches to the analysis of discrete excitable media in one dimension is to study certain embedded particle systems, where particles keep track of boundaries between intervals of distinct colors \cite{fisch1990one, fisch1991cyclic, durrett1991some}. This technique can be best seen for the 3-color critical CCA ($r=1$), whose sample dynamics is shown in Figure \ref{fig:3CCA_FCA}. Since colors `eat' nearby colors of one less modulo 3, one can place a `right particle' (resp., `left particle') at edges of colors $(i+1,i)$ (resp., $(i,i+1)$) mod 3, which points at the direction of excitation. From each 3-coloring $\xi_{t}$, we may denote such edge particle assignment as $d \xi_{t}$, where $\xi_{t}(x)=+1$ if there is a right (resp., left) particle on the edge $(x,x+1)$ at time $t$. The CCA time evolution $\xi_{t}\mapsto \xi_{t+1}$ induces `arrow dynamics' $d\xi_{t}\mapsto d\xi_{t+1}$ for the edge configurations through the following commutative diagram:
				\begin{align}\label{eq:arrow_diagram}
					\begin{gathered}
						\xymatrix{
							\xi_{t} \ar[rrr]^{\textup{CCA/FCA}} \ar[d]_{d}&&& \xi_{t+1} \ar[d]^{d} \\
							d \xi_{t}
							\ar[rrr]^{\textup{arrow dynamics}}&&& d \xi_{t+1}
						}
					\end{gathered}
				\end{align}
				For the critical 3-color CCA, the induced arrow dynamics is extremely simple. All right and left particles move to the nearest edge that they are pointing at simultaneously; If any two opposing particles must cross or occupy the same edge, then they annihilate each other and are removed from the system. One can clearly see that the boundaries between monochromatic regions in the 3-color CCA dynamics in Figure \ref{fig:3CCA_FCA} behave in such a way. 
				
				Next, we move our attention to the 3-color critical FCA. Clustering for this system with exact asymptotic rate was established in \cite{lyu2017persistence} by using a similar embedded edge particle systems. One can see that, due to the nature of FCA where excitation can only be initiated by sites with blinking color (i.e., 0), boundaries between differently colored region do not move every iteration as in the CCA case; rather, one can see in Figure \ref{fig:3CCA_FCA} that they move in either three or four steps. Nonetheless, it is evident that right and left particles move ballistically and annihilate upon collision.  
				
				Following \cite{lyu2017persistence}, the construction of edge configuration $d\xi_{t}$ from a 3-coloring $\xi_{t}$ for the 3-color critical FCA is the same as the 3-color critical CCA. 
				For instance, the configuration $\xi_{0}=(\cdots00122\cdots)$ is assigned with two right particles on the edges of colors $(0,1)$ and $(1,2)$. The particles move in one step following their direction whenever their tail site has color zero, where opposing particles that cross or coincide mutually annihilate.  
				In \cite{lyu2017persistence}, it was shown that this arrow dynamics commutes with the FCA dynamics for all times $t\ge 1$. The reason that we need to forbid the very first iteration $\xi_{0}\mapsto \xi_{1}$ is that initial particles can flip their direction. For instance, the initial configuration $\xi_{0}=(\cdots00122\cdots)$ evolves to $\xi_{1}=(\cdots11100\cdots)$, which now has a single left particle on the edge between colors $(1,0)$. Such `particle flipping' only occurs for every triple of colors $(0,1,2)$ and $(2,1,0)$. By a simple back-tracking argument, one can show that such configurations cannot occur for all times $t\ge 1$ (for instance, $(0,1,2)$ must have been $(2,0,*)$ with $*=1$ one iteration ago, but $*=1$ is pulled by 0 so it does not evolve to 2 in one iteration, a contradiction).

				\subsection{Arrow dynamics for general CCA and FCA} 
				
				
				
\begin{figure}[h!]
					\centering
					\includegraphics[width=1\textwidth]{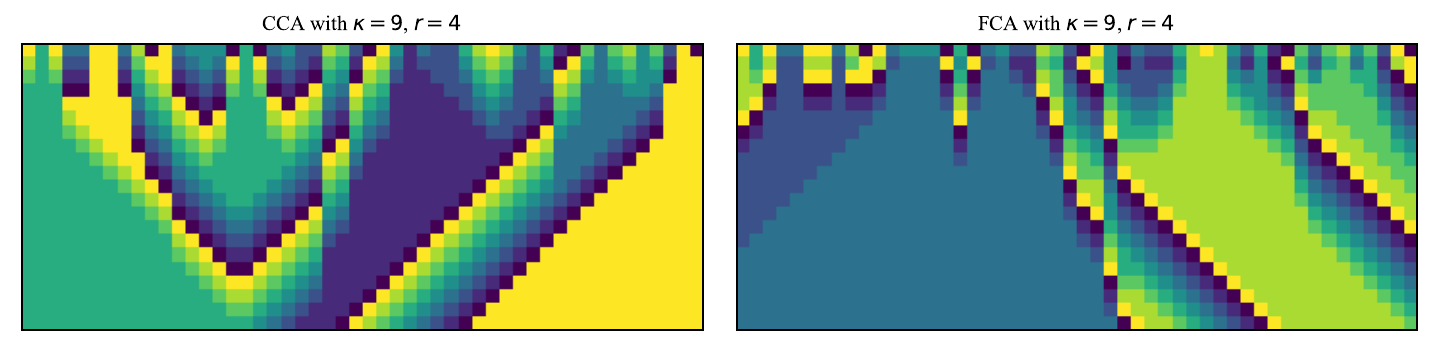}
					\vspace{-0.5cm}
					\caption{ Simulation of 9-color critical CCA and FCA ($r=4$) on 50 nodes in $\mathbb{Z}$ with 20 iterations. Time goes from top to bottom. For CCA configurations at every iteration are shown, while for FCA only the ones at times $9 t$ are shown.
					}
					\label{fig:9CCA_FCA}
				\end{figure}

				For $\kappa \ge 4$, edges may have large color differences. A natural way to interpret the dynamics of the boundaries between monochromatic regions is to think of stacks of right or left particles moving ballistically, annihilating each other when opposing particles meet. See Figure \ref{fig:9CCA_FCA} for examples of 9-color CCA and FCA.

				\subsubsection{Arrow dynamics for CCA}
				We now define the edge configuration $d\eta:\Z\rightarrow [-\lfloor \kappa/2 \rfloor,\lfloor \kappa/2 \rfloor]$ corresponding to a $\kappa$-coloring $\eta:\mathbb{Z}\rightarrow \Z_{\kappa}$ for CCA as 
				\begin{align}\label{eq:CCA_edge_configuration}
					d\eta(x) := \eta(x)-\eta(x+1) \quad (\text{mod $\kappa$}) 
				\end{align}
				where we identify each edge $(x,x+1)$ with the left site $x$. In particular, if $\kappa=3$, then $d\eta$ assigns $-1,0,1$ on the edges of $\Z$, where $+1$ (resp. $-1$) corresponds to right (resp., left) particles on edges where the colors decrease (increase) by one mod 3; 0 is assigned on edges with the same color. For general $\kappa$, we will think of $|d\eta(x)|$ particles of the same direction indicated by the sign of $d\eta(x)$ stacked on the edge $(x,x+1)$, provided $|d\eta(x)|<\kappa/2$. If $\kappa$ is even and $|d\eta(x)|=\kappa/2$, then we assign $\kappa/2$ \textit{bidirectional} arrows on the edge $(x,x+1)$. 
				
				Given CCA dynamics $(\eta_{t})_{t\ge 0}$, \eqref{eq:CCA_edge_configuration} induces a sequence of edge configurations $(d\eta_{t})_{t \ge 0}$. It is possible to evolve these edge configurations directly via certain `arrow dynamics' without evolving the $\kappa$-coloring using CCA dynamics (see the commutative diagram  $\xi_{t}=\eta_{t}$). The detailed arrow dynamics are described below.

\begin{figure}[h!]
					\centering
					\includegraphics[width=0.81\textwidth, height=55mm]{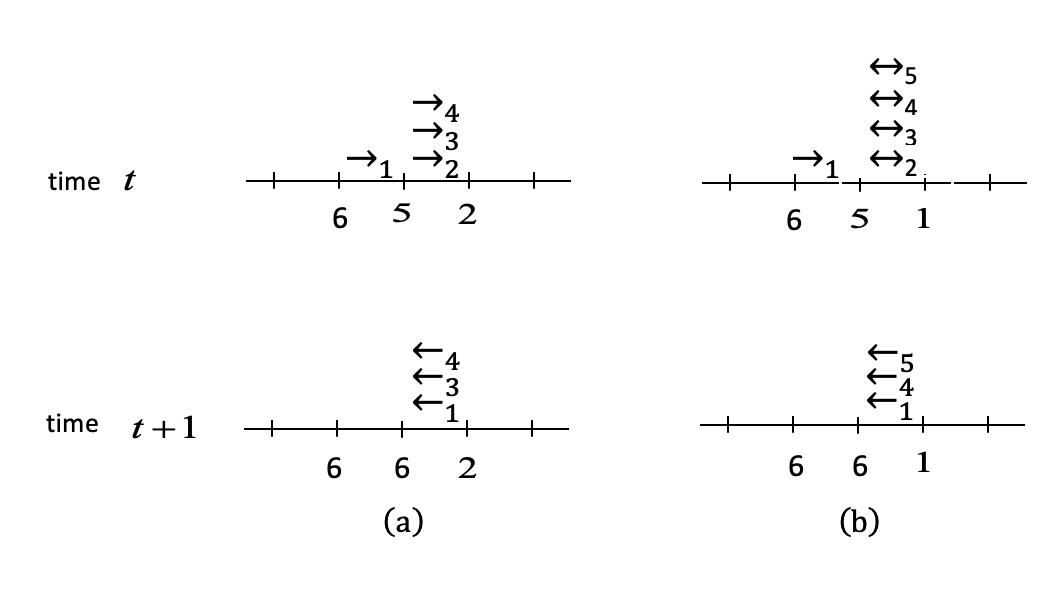}
					\vspace{-0.5cm}
					\caption{ Arrow flipping rule for CCA with (a) $\kappa=7$, $r=2$ and (b) $\kappa=8$, $r=2$. One arrow is lost if $\kappa$ is odd  and two if $\kappa$ is even. Note also the labels of the deleted arrows. Stacks of $\frac{\kappa}{2}$ arrows are bidirectional.
					}
					\label{fig:CCA_arrow3}
				\end{figure}

				\begin{description}[itemsep=0.1cm]
					\item[(i)] \textbf{Labeling}: Arrows are labeled in a pseudo-lexicographic ordering as follows. Labels increase along the lattice bonds from left to right. On a given edge, a stack of right-pointing and bidirectional arrows is labeled in increasing order from the bottom up, and vice-versa for left-pointing stacks. 
					
					\item[(ii)] \textbf{Release}: Nonzero stacks with $\leq r$ arrows are \textit{active}. Stacks with $>r$ arrows are \textit{inactive} and we refer to them as \textit{bloackades}. 
					For all active stacks, top arrows (that is, right (resp., left) arrows with maximum (resp., minimum) labels) are released simultaneously. Released arrows try to jump to the next edge in the direction it is pointing. 
					
					\item[(iii)] \textbf{Collisions}: 
					If two released opposing arrows need to cross each other or occupy the same edge, they mutually annihilate. If a released arrow attempts to jump onto an inactive stack of $<\kappa/2$ opposing arrows, it annihilates the opposing arrow with the closest label. 
					
					If a released arrow attempts to jump onto an inactive stack of $\kappa/2$ bidirectional arrows for $\kappa$ even and if there is another left arrow jumping onto the same stack, then each of the jumping arrows annihilates with a bidirectional arrow in the stack with the closest label.

					If a released arrow attempts to jump onto an inactive stack of $<\lfloor \kappa/2 \rfloor$ arrows of the same direction (including empty stacks), it is added to the stack and no annihilation occurs. 
					Note that the pseudo-lexicographic ordering is preserved after every collision event. 
					

					\item[(iv)] \textbf{Stack flipping}: Suppose $\kappa$ is odd and there is an inactive stack of  $\lfloor \kappa/2 \rfloor$ right (resp., left) arrows and another right (resp., left) arrow with label $\ell$ attempts to jump on it. Suppose that there is no left (resp., right) arrow jumping onto the same inactive stack at the same time.
					Then the entire stack flips to a stack of left (resp., right) arrows, including a left (resp., right) arrow of label $\ell$. In the flipped stack, there are  $\lfloor \kappa/2 \rfloor$ left (resp., right) arrows. In this interaction, one right arrow in the stack with labels closest to $\ell$ are annihilated. 
					
					Suppose $\kappa$ is even and there is an inactive stack of $\kappa/2$ bidirectional arrows and a right arrow with label $\ell$ attempts to jump on it and there is no left arrow jumping onto the same stack. Then the bidirectional arrows become $(\kappa/2) -1$ left arrows, including one with label $\ell$. In this interaction, two bidirectional arrows in the stack with labels closest to $\ell$ are annihilated. 
					See 
					Figure \ref{fig:CCA_arrow3}.

				\end{description}
				
				We remark that stack flipping for CCA may occur only if $r<\lfloor \kappa/2 \rfloor$, since otherwise stacks of size $\lfloor \kappa/2 \rfloor$ is always active. 
				
				\begin{figure}[h!]
					\centering
					\includegraphics[width=0.85\textwidth , height=83mm]{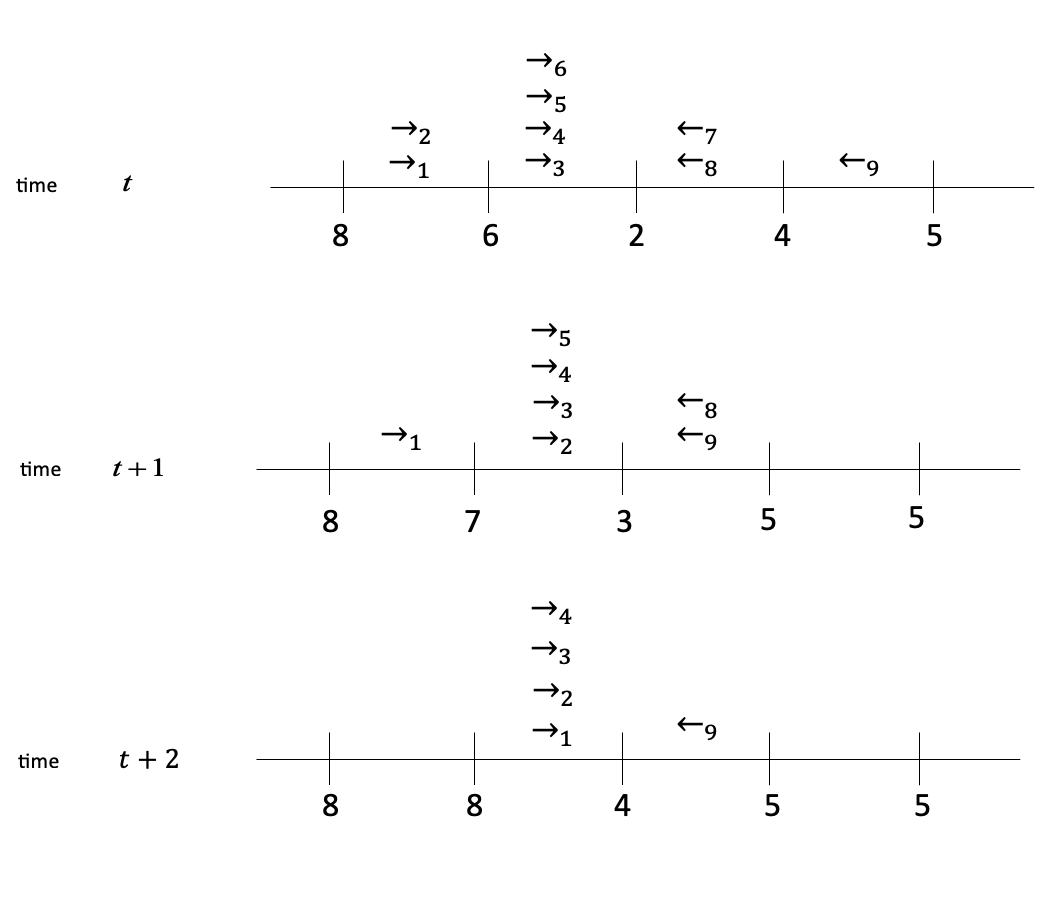}
					\vspace{-0.5cm}
					\caption{ Arrow collisions for CCA with $\kappa=11$ and $r=4$. Pseudo-Lexicographic labeling of arrows ensures proper 
						monotonic ordering.
					}
					\label{fig:CCA_arrow1}
				\end{figure}

				Below we show that the arrow dynamics for CCA is indeed compatible with the CCA dynamics. 
				
				\begin{prop}\label{prop:arrow_CCA}
					Suppose $r<\lfloor \kappa/2 \rfloor$ or $\kappa$ is odd and $r=\lfloor \kappa/2 \rfloor$. Then 
					The diagram in \eqref{eq:arrow_diagram} commutes for CCA with the arrow dynamics described above. 
				\end{prop}
				
				\begin{proof}
					Consider the edge $(0,1)$ at time $t$. Suppose $d\eta_{t}(0)=a\in \{-\lfloor \kappa/2 \rfloor,\, \lfloor \kappa/2 \rfloor\}$. By symmetry, without loss of generality, we may assume $a>0$, $\eta_{t}(0)=0$ and $\eta_{t}(1)=-a$. In terms of the arrow dynamics, there is a stack of $a$ right arrow on the edge $(0,1)$. We need to show the value of $d\eta_{t+1}(0)$ computed using the arrow dynamics coincides with the value computed by first evolving the coloring by the CCA dynamics one iteration and then computing the edge configuration. 
					
					We first consider the case $a\le r$ so that the stack on the edge $(0,1)$ is active at time $t$, so a single right arrow from this stack is released to the edge $(1,2)$. Since $a\le r$, the stack on the edge $(0,1)$ does not flip, so we should have 
					$a$ or $a-1$ right arrows on the edge $(0,1)$ at time $t+1$ depending on whether another right arrow jumps from the edge $(-1,0)$ or not. In terms of CCA dynamics, note that $1$ excites $0$ (since $d\eta_{t}(0)>0$) so that $\eta_{t+1}(0)=1$. If 0 is also excited (necessarily by $-1$), then $\eta_{t+1}(0)=1$ and hence $d\eta_{t+1}(0)=a$; otherwise $\eta_{t+1}(0)=0$ and $d\eta_{t+1}(0,1)=a$. We have shown that 
					\begin{align*}
						&\textup{($\#$ of right arrows on the edge $(0,1)$ at time $t+1$)} \\
						&\qquad = (a-1) + \mathbf{1}(\textup{a right arrow jumps from the edge $(-1,0)$ at time $t$}) \\
						&\qquad = (a-1) + \mathbf{1}(\textup{$0$ is excited by $-1$ at time $t$}) \\
						&\qquad = d\eta_{t+1}(0),
					\end{align*}
					which is the desired conclusion. 
					
					Second, consider the case $a> r$ so that the stack on the edge $(0,1)$ is inactive and hence is a blockade at time $t$. We will consider two cases depending on whether this stack flips or not. First, suppose it does not flip. Then note that 
					\begin{align*}
						&\textup{($\#$ of right arrows on the edge $(0,1)$ at time $t+1$)} \\
						&\qquad = a + \mathbf{1}(\textup{a right arrow jumps from the edge $(-1,0)$ at time $t$}) \\
						&\hspace{2cm} - \mathbf{1}(\textup{a left arrow jumps from the edge $(1,2)$ at time $t$})  \\
						&\qquad = a + \mathbf{1}(\textup{$0$ is excited by $-1$ at time $t$}) - \mathbf{1}(\textup{$1$ is excited by $2$ at time $t$}) \\
						&\qquad = d\eta_{t+1}(0),
					\end{align*}
					as desired. Lastly, suppose the inactive stack on the edge $(0,1)$ flips. Under the assumption that $\eta_{t}(0)=0$ and $\eta_{t}(1)=-a<0$, 
					this is precisely when $a=\lfloor \kappa/2 \rfloor>r$ and the larger color of $\eta_{t}(0)$ and $\eta_{t}(1)$ increments by one, so $\eta_{t+1}(0)=1$ and $\eta_{t+1}(1)=-a$. It follows that 
					\begin{align}
						d\eta_{t+1}(0) = -\lfloor \frac{\kappa-1}{2} \rfloor. 
					\end{align}
					The arrow dynamics for stack flipping accounts precisely the above. 
				\end{proof}

				\subsubsection{Arrow dynamics for FCA}
				\label{sec:FCA_arrow}
				
				For FCA, we define the edge configuration corresponding to a $\kappa$-coloring $X:\Z\rightarrow \Z_{\kappa}$ as a map $dX:\Z\rightarrow [-m,m]$ for $m:=\lfloor \kappa/2 \rfloor$, where 
				\begin{align}\label{eq:FCA_edge_config}
					dX(x) =
					\begin{cases}
						X(x+1)-X(x)  & (\textup{mod}\,\, \kappa)\\
						m & \text{if $\kappa=2m$, $|X(x+1)-X(x)|=m$, and $X(x)\in [r+1, \kappa]$} \\
						-m & \text{if $\kappa=2m$, $|X(x+1)-X(x)|=m$, and $X(x+1)\in [r+1,\kappa]$}.
					\end{cases}
				\end{align}
				If $\kappa$ is odd, the definition of edge configuration $dX$ is similar to that for CCA except the sign is the opposite, since excitation goes from smaller colors (in fact, only from the blinking color 0) to larger ones in FCA but the other way around in CCA. When $\kappa$ is even, there could be edges with maximum color difference $m$, and care has to be taken to assign the direction of the particles on such edges. For CCA, we assigned a stack of $m$ bidirectional arrows on such edges. For FCA, we do not need to introduce bidirectional edges. Recalling that excitation in FCA goes from color 0 to colors in $[1,r]$, we assign the direction of the particles according to which of the two sites on the edge $(x,x+1)$ becomes the blinking color 0 first. Since colors in $[r+1,\kappa]$ advance to color zero without interruption, $x$ achieves color 0 before $x+1$ given $dX(x)=m$ if and only if $x$ currently has a color in $[r+1,\kappa]$.

				
				Given FCA dynamics $(X_{t})_{t\ge 0}$, \eqref{eq:FCA_edge_config} induces a sequence of edge configurations $(dX_{t})_{t \ge 0}$, which we can interpret as an embedded dynamical system with stacks of arrows. Unlike CCA, it would be convenient to use the FCA dynamics $(X_{t})_{t\ge 0}$ to define such embedded arrow dynamics since only arrows starting from sites with color 0 are `active' in the arrow dynamics. 
				For each edge $(x,x+1)$ with $dX_{t}(x)=k>0$, we imagine a stack of $k$ unit particles residing on the edge $(x,x+1)$ at time $t$, which we will represent with arrows,  directed from $x$ to $x+1$. Similarly, negative values of $dX_t$ correspond to stacks of arrows moving left. Note also that $\lVert dX_{t} \rVert_{\infty}\le \lfloor \kappa/2 \rfloor$, so the height of the stacks on the edges is uniformly bounded by $\lfloor \kappa/2 \rfloor$.  The arrow dynamics for FCA is almost identical to that of FCA, but it crucially differs in that stacks on the edges where arrows are pointing away from  color 0 are active, and particles on such active stacks are released. Hence, unlike in the arrow dynamics for CCA, stacks with $<r$ arrows can be inactive unless the site at `tail' have the blinking color 0.



\begin{figure*}[ht]
					\centering
					\includegraphics[width=1 \linewidth]{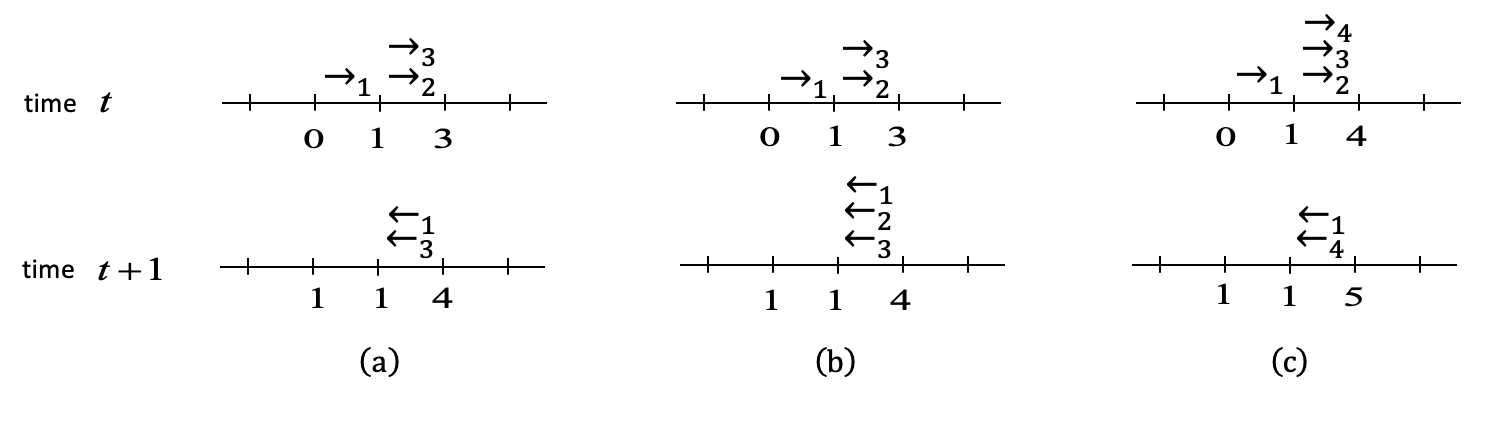}
					\caption{\footnotesize Examples of arrow dynamics for $\kappa$-color FCA with (a) $\kappa=5$ and (b), (c) $\kappa=6$ with $r=2$ for all three cases. (a) and (b) show instances of stack flipping but (c) is the usual arrow-arrow annihilation. 
					}
					\label{fig:FCA_flipping}
				\end{figure*}

				\begin{description}[itemsep=0.1cm]
					\item[(i)] \textbf{Labeling}: Same as in CCA but note that there are no bidirectional arrows for FCA. 
					
					
					
					\item[(ii)] \textbf{Release}: Nonzero stacks with $\leq r$ arrows on edges where the sites toward the tails of the arrows have color 0 are called 
					\textit{active}; stacks of $>r$ arrows are  called \textit{blockades}; all other stacks with nonzero arrows are \textit{inactive}. 
					For all active stacks, top arrows (that is, right (resp., left) arrows with maximum (resp., minimum) labels) are released simultaneously. The released arrows try to jump to the next edge in the direction it is pointing.

					\item[(iii)] \textbf{Collisions}:  Same as in CCA.
					

					\item[(iv)] \textbf{Stack Flips}: Suppose there is an inactive stack of  $\lfloor \frac{\kappa-1}{2} \rfloor$ right (resp., left) arrows and another right (resp., left) arrow with label $\ell$ attempts to jump on it. 
					Then the entire stack flips to a stack of $\lfloor \kappa/2 \rfloor$ left (resp., right) arrows, including a left (resp., right)  arrow of label $\ell$. Note that when $\kappa$ is odd, then one arrow in the inactive stack with the closest label to $\ell$ is annihilated during flipping. When $\kappa$ is even, no annihilation occurs due to flipping. See Figure \ref{fig:FCA_flipping} for illustration. 
				\end{description}
				

				\begin{figure*}[ht]
					\centering
					\includegraphics[width=1.06\linewidth, height=55mm]{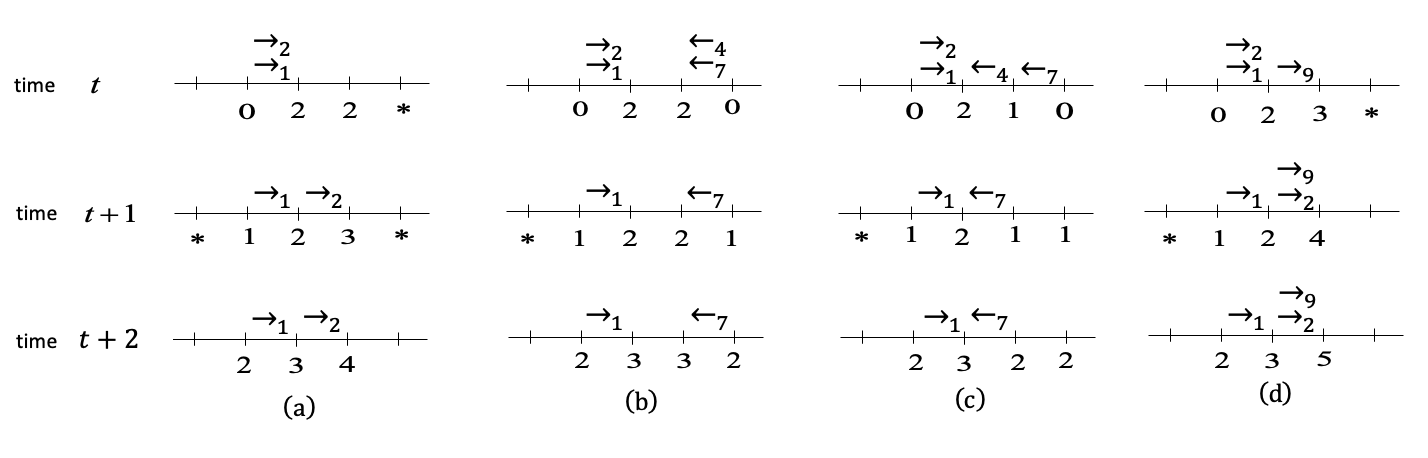}
					\caption{\footnotesize Queuing and annihilation rules for the particle system expansion of the $6$-color FCA on $\mathbb{Z}$ with $r=3$. We denote particle $(i,\R)$ by $\rightarrow_{i}$ and similarly for left arrows. $*\in \mathbb{Z}_{6}\setminus \{0\}$. Note deleted arrows. 
					}
					\label{fig:aqps}
				\end{figure*}

				To motivate the above construction of the arrow dynamics for FCA, we walk through some examples when $\kappa=6$ and $r=3$.  Recall that $0$ is the blinking color and there are three post-blinking colors $1,2$ and $3$ whose update to the next color is inhibited when in contact with color $0$. Suppose $X_{t}(x)=0$ and $X_{t}(x+1)=2$ so that there are two $\mathbf{r}$ particle on the edge $(x,x+1)$ and site $x$ blinks at time $t$. Suppose that the edge $(x+1,x+2)$ is vacant and site $x+3$ does not blink at time $t$; so $X_{t}(x+2)=2$ and $X_{t}(x+3)=*\ne 0$ (see Figure \ref{fig:aqps} $(a)$). Then at time $t+1$ sites $x,x+1$, and $x+2$ have colors $1,2,$ and $3$ respectively, so there is a single $\mathbf{r}$ particle on each of the edges $(x,x+1)$ and $(x+1,x+2)$ at time $t+1$. We view this as the bottom $\mathbf{r}$ particle on the edge $(x,x+1)$ having moved onto the vacant edge $(x+1,x+2)$. If $X_{t}(x+3)=0$, then both the bottom particles on edges $(x,x+1)$ and $(x+2,x+3)$ try to move into the vacant edge $(x+1,x+2)$, resulting in their annihilation at time $t+1$ (see Figure \ref{fig:aqps} (b)). Similar annihilation occurs when these particles are closer to each other, i.e., the edge $(x+1,x+2)$ also has an $\mathbf{l}$ particle (see Figure \ref{fig:aqps} (c)). Lastly, consider the case when there are right particles on two consecutive edges, namely, $(x,x+1)$ and $(x+1,x+2)$, and suppose site $x$ blinks at time $t$. If $dX_{t}(x+1)=1$, then the bottom $\mathbf{r}$ particle on the edge $(x,x+1)$ migrates to the top of the queue on the edge $(x+1,x+2)$ (see Figure \ref{fig:aqps} (d)).

				Below we show that the arrow dynamics for FCA is indeed compatible with the FCA dynamics. 
				
				\begin{prop}\label{prop:arrow_FCA}
					Suppose $r\le \lfloor \kappa/2 \rfloor$. Then 
					The diagram in \eqref{eq:arrow_diagram} commutes for FCA with the arrow dynamics described above. 
				\end{prop}
				
				\begin{proof}
					The argument is similar to that for the proof of Prop. \ref{prop:arrow_CCA}. Details are omitted.
				\end{proof}

				\section{CCA and FCA in the undercoupled regime}
				\label{sec:CCA_undercoupled}



				Without loss of generality, we shall prove that the origin fixates. If it didn't,  it would need to change colors infinitely often. However, each color change corresponds to an arrow crossing the origin. The following lemma shows that if there are infinitely many  crossings of the origin then arrows must eventually have an arbitrarily large label (i.e. they come from arbitrarily far away).

				\begin{lemma}\label{lem:CCA_fixation_suff}
					Assume that as $t\rightarrow \infty$, there are infinitely many arrow-crossings of the origin and denote by $A_t$ the set of  arrows that have crossed the origin at time $t$. Then as $t\rightarrow\infty$ we have that $|A_t|\rightarrow \infty$.
				\end{lemma}

				\begin{proof} Suppose that $|A_{\infty}|$ is finite. Then there would be at least one arrow  that crosses the origin infinitely many times and therefore switches directions infinitely often. Without loss of generality suppose this arrow is initially to the left of the origin. Let $k$ be the minimal label of such arrow. Observe that for $k$ to cross the origin infinitely many times, it must switch directions from right to left at locations to the right of the origin infinitely many times. In each direction flip, $k$ must be involved in an arrow interaction of the form stack flipping, which means that at these times the location of $k$ is entered by an arrow from the left that has strictly smaller label than $k$ and which therefore also started to the left of the origin and passed the origin. Let us keep track of the set of arrows that induce the direction flips in $k$ by naming the arrow involved in the $m$-th direction flip $k_m$. Denote $K=\cup_{i=1}^{\infty}k_i$ and observe that $K\subset A_{\infty}$. Since $K$ is finite and arrows cannot cross, then there must be another arrow $k'\in K$ that also crosses the origin infinitely many times. This contradicts the minimality of $k$.
				\end{proof}



				Thus, by Lemma \ref{lem:CCA_fixation_suff}, 
				if the origin is crossed by arrows infinitely many times, then there are infinitely many distinct labels of all arrows that cross the origin, and hence, for every $n$ there exists a site $z$ with $|z|>n$ such that the arrow located at, say  $(z-1,z)$, crosses the origin. To proceed with the argument, we introduce the edge elimination function $\varphi$. This is effectively a book-keeping device to count the net number of live arrows vs blockades in an interval. Here we give an edge elimination function that suffices for our purposes but we remark that it may be improved for a more optimal result. For a $\kappa$-coloring $\xi:\mathbb{Z}\rightarrow \Z_{\kappa}$ for CCA or FCA,
				\begin{equation}
					\label{eq1} 
					\varphi\of{d\xi(x)}=\begin{cases}
						-|d\xi(x)| \quad \quad \quad &\text{if}\quad  0\leq |d\xi(x)|\leq r\\
						|d\xi(x)|-2r\quad \quad \quad & \text{if}\quad  r\leq |d\xi(x)|\leq \lfloor\frac{\kappa}{2}\rfloor.
					\end{cases}
				\end{equation}
				The function $\varphi$ above gives a conservative estimate of the total blockade capacities that all arrows on the edge $(x,x+1)$ can reduce, if CCA or FCA is initialized with $\xi$. It is conservative in that we do not consider the effect of arrows annihilating each other and always counted against blockade capacities. It is possible, however, to consider such finite-time effects (e.g., arrow-arrow annihilation, immediate blockade formation) to improve the resulting sufficient condition on fixation as done by Fisch \cite{fisch1991cyclic} for $\kappa$-color CCA with $r=1$.

				Lemma \ref{lem:fixation_suff} below gives a necessary condition for an arrow to eventually reach the origin in terms of $\varphi$. For any $z\in \Z$ 
				define the first time $T(z)$ from time that size $z$ influences the origin: 
				\begin{align}\label{eq:1}
					T(z)=\inf \left\{t \, \bigg| \,   \begin{matrix}\text{time that there is an arrow visiting the origin at time $t$,} \\ \textup{which initially was at an edge containing $z$} \end{matrix} \right\}.
				\end{align}


				\begin{lemma}\label{lem:fixation_suff}
					Let $\xi_{t}$ denote the $\kappa$-color CCA or FCA on $\Z$ with  $r< \lfloor \kappa/2 \rfloor$. The following two statements combined give us sufficient conditions for fixation: 
					\begin{description} 
						\item[(i)] $\xi_t$ fixates if $ \lim_{n\rightarrow\infty}\P(T(z)<\infty\,\, \text{for some}\,\,  z<-n)=0$. 
						
						\vspace{0.1cm}
						\item[(ii)]  For an arrow coming from  $z\in \Z^{-}$  to reach the origin we must have an interval $I$ containing $[z,0]$ such that $\sum_{I}\varphi(d\xi_{0}(x))<0$. More precisely, we have the following event inclusion
						\begin{align*}
							\{T(z)<\infty \quad \text{for some } z<-n\}\subset \left\{\sum_{I}\varphi(d\xi_{0}(x))< 0\right\}.
						\end{align*}
					\end{description}
				\end{lemma}

				\begin{proof} See  \hyperlink{appendix}{Appendix A} (adapted from \cite{bramson1989flux, fisch1990one}). \end{proof}

				We are now ready to prove of Theorem \ref{thm:main1} \textbf{\textup{(i)}}. 
				
				\begin{proof}[\textbf{Proof of Theorem \ref{thm:main1} \textbf{\textup{(i)}}}]

					To conclude the argument, all we need, by Lemma \ref{lem:fixation_suff}, is to compute the asymptotics of the following expectation 
					\begin{align}\label{bcount}
						\E\left[ \sum_{-n}^n\varphi(d\xi_{0}(x)) \right]
					\end{align}
					in terms of $\kappa(r)$. When this expectation is  non-positive, the two parts of \textbf{Lemma} \ref{lem:fixation_suff} give us the required vanishing probability for fixation. Without loss of generality, we consider  we are in the case of odd $\kappa$ (terms don't change much when $\kappa$ is even). Then using (\ref{eq1}) and the translational symmetry of the original random configuration one has:
					
					\begin{align}(\ref{bcount})=&\sum_{m=1}^r- \frac{2m}{\kappa}+\sum_{m=r+1}^{\lfloor\frac{\kappa}{2}\rfloor}\frac{2(m-2r)}{\kappa}\nonumber\\
						=&\frac{-r(r+1)}{\kappa}+\frac{1}{\kappa}\of{\lfloor\frac{\kappa}{2}\rfloor\of{\lfloor\frac{\kappa}{2}\rfloor+1}-(r+1)(r+2)}-\frac{4r}{\kappa}\of{\lfloor\frac{\kappa}{2}\rfloor-(r+1)}\label{approx1}
					\end{align}
					Asymptotically in large $\kappa$, this gives us:
					\begin{align*}\text{RHS}(\ref{approx1})\approx-\frac{r^2}{\kappa}+\frac{\kappa}{4}-\frac{r^2}{\kappa}-2r+\frac{4r^2}{\kappa}.
					\end{align*} 
					Combining like terms, for large enough $\kappa$, one has:
					\begin{align}\label{lim}\lim_{n\rightarrow \infty}\E\of{\sum_{x=-n}^n\varphi(d\eta(x))}\longrightarrow\frac{2r^2}{\kappa}-2r+\frac{\kappa}{4}.
					\end{align}
					Hence, we have fixation if $r<\frac{2-\sqrt{2}}{4}\kappa$. ($\approx 0.14644 \kappa $)
				\end{proof}

				\section{Critical CCA and FCA}
				\label{sec:critical}
				
				In this section, we prove Theorem \ref{thm:main1} \textbf{(ii)} concerning CCA and FCA in the critical regime. Suppose $r=\lfloor \frac{\kappa-1}{2} \rfloor$ for CCA and $r=\lfloor \kappa/2 \rfloor$ for FCA.  We will first show clustering for these systems using a soft ergodic argument. Then, we will prove the asymptotic rate of clustering in \eqref{eq:clustering_rate} when $r=\lfloor \kappa/2 \rfloor$ except even $\kappa$ for CCA.

				\subsection{Proof of clustering with an ergodic argument}
				
				Here we show that CCA and FCA in the critical regime clusters using an ergodic argument. We first show that there is no spontaneous blockade formation.

				\begin{lemma}[No spontaneous emergence of blockades]\label{lem:FCA_undercoupled_flipping_transient1}
					Suppose $\kappa=2m$ for some integer $m\ge 2$ and let $r=m-1$. 
					For each site $x\in \Z$, if $|d\xi_t(x)|=m$ for some $t\ge m$, then $|d\xi_{s}(x)|=m$ for some $0 \le s < t$. In particular, blockades are not created spontaneously after time $m$. 
				\end{lemma}

				\begin{proof}
					Suppose $|d\xi_{t}(x)|=m$ for some times $t\ge m$. Further assume that $|d\xi_{s}(x)|<m$ for all $t-m \le s < t$. We will show that there is a contradiction by backtracking the dynamics $\kappa$ steps backward in time. First, consider the CCA dynamics. By symmetry, suppose $\xi_{t}(x)=0$ and $\xi_{t}(x+1)=m$ so that there is a blockade $d\xi_{t}(x)=m$ at time $t$ on the edge $(x,x+1)$. Then there is the last time $s$ before time $t$ that there is no blockade on the edge $(x,x+1)$, i.e., $d\xi_{s}(x)<m$. By symmetry, we may assume $\xi_{s}(x)=0$ and $\xi_{s}(x+1)=m-1$. Then $x$ must be excited at time $s$ so that $\xi_{s+1}(x+1)=m$. However, we must also have $\xi_{t+1}(x)=1$ since $x$ is excited by $x+1$ at time $s$ for $d\xi_{s}(x)=m-1=r$ (see \eqref{eq:CCA_blockade_no_emergence}). This is a contradiction. 
					\begin{align}\label{eq:CCA_blockade_no_emergence}
						\begin{matrix}
							\textup{time}: & t & t-1 & &  \\
							\hline
							x+2 :  &  & &  & &  \ge m  \\
							x+1 :& m & m  &\hdots & m & m-1 \\
							x :   &  0 & 0  & \hdots & 0 & 0  
						\end{matrix}
					\end{align}
					
					
					Next, consider the FCA dynamics. Suppose there is a blockade at time $t\ge \kappa$ at edge $(x,x+1)$ but there is not at times in $\{ t-\kappa,\dots,t-2,t-1 \}$. By symmetry, assume $\xi_{t}(x)=a\in \{0,1,\dots,m\}$  and $\xi_{t}(x+1)=a+m$. Then either of the following holds: (1) $\xi_{t-a}(x)=0$ and $x$ is not excited during time interval $[t-a,t)$; or (2) $\xi_{t-a-1}(x)=0$ and $x$ is excited exactly once (by $x-1$) during time interval $[t-a,t)$. There are no other cases since $x$ cannot be excited more than once during this period (since $x-1$ is the only neighbor of $x$ that can excite it during this period and every node in FCA returns to zero at most once in every $\kappa$ iterations).  
					
					First consider the case (1). Note that $x+1$ has colors in $\{m,m+1,\dots,\kappa-1\}$ during time interval $[t-a,t)$ so it cannot get excitation during this period. Thus the colors of $x$ and $x+1$ decrease by one when we backtrack the dynamics from time $t$ to $t-a+1$, giving $\xi_{t-a+1}(x)=1$ and $\xi_{t-a+1}(x+1)=m+1$. Then $\xi_{t-a}(x)=0$ and $\xi_{t-a}(x+1)\in \{m,m+1\}$, but we get contradiction, for if $\xi_{t-a}(x+1)=m$, then $x+1$ is excited at time $t-a$ by $x$ so that $\xi_{t-a+1}(x+1)=m$ contrary to $\xi_{t-a+1}(x+1)=m$, and if $\xi_{t-a}(x+1)=m+1>r$, then $x+1$ is not excited at time $t$ so $\xi_{t-a}(x+1)=m+2$, again contrary to $\xi_{t-a+1}(x+1)=m$. See \eqref{eq:FCA_blockade_no_emergence1}. 
					\begin{align}\label{eq:FCA_blockade_no_emergence1}
						\begin{matrix}
							\textup{time}: & t & t-1 & & & t-a+1 &  t-a  \\
							\hline
							x+1 : & a+m & a+m-1 & a+m-2  &\cdots & m+1 & m \\
							x :  & a  &  a-1 & a-2  & \cdots & 1 & \mathbf{0}  
						\end{matrix}
					\end{align}
					
					Lastly, consider the case (2). In this case, $\xi_{t-a}(x+1)=m$ as before, but $\xi_{t-a}(x)=1$ and $\xi_{t-a-1}(x)=0$ since $x$ is excited by $x-1$ at some time between $t-a-1$ and $t-1$. Below in \eqref{eq:FCA_blockade_no_emergence}, we depict such a situation where $x$ is excited by $x-1$ at time $t-1$ for illustration purpose. 
					\begin{align}\label{eq:FCA_blockade_no_emergence}
						\begin{matrix}
							\textup{time}: & t & t-1 & & &t-a+1 & t-a &  t-a-1  \\
							\hline
							x+1 : & a+m & a+m-1 & a+m-2  &\cdots & m +1 & m & m \\
							x :  & a  &  a & a-1  & \cdots &  & 1 & \mathbf{0}  \\
							x-1 : & 1 & \mathbf{0} 
						\end{matrix}
					\end{align}
					Since $t\ge m$, $t-a-1\ge 0$. Since $\xi_{t-a-1}(x)=0$ and $\xi_{t-a}(x+1)=m = r-1$, we cannot have $\xi_{t-a-1}(x+1)=m-1$ since then $x+1$ is excited by $x$ at time $t-a-1$ so that  $\xi_{t-a}(x+1)=m-1$, contrary to $\xi_{t-a}(x+1)=m$. Thus we must have $\xi_{t-a-1}(x+1)=m$. But then we have a blockade on the edge $(x,x+1)$ at time $t-a-1\ge 0$, contrary to the hypothesis that there was no blockade on this edge at prior times. This finishes the proof.
				\end{proof}

				In the following proposition, we show that critical CCA and FCA must either fixate or cluster. We include the case of FCA with $\kappa$ even and $r=\kappa/2$ in the statement, which seems to cluster according to simulation (see Figures \ref{fig:CCA_FCA_sim1}, \ref{fig:CCA_FCA_sim2}, and \ref{fig:4CCA_FCA}) but is excluded from Theorem \ref{thm:main1} \textbf{(ii)} since our current argument does not rule out the possibility of fixation into non-interacting colors.
				
				\begin{prop}\label{prop:clustering_blockade}
					Let $\xi_{t}$ denote $\kappa$-color CCA or FCA, where 
					$r\in \{\lfloor \kappa/2 \rfloor-1, \lfloor \kappa/2 \rfloor\}$ except the case of CCA with $\kappa$ even and $r=\kappa/2$. Then either the system clusters or a positive density of blockades persists. 
				\end{prop}
				
				\begin{proof}
					
					By symmetry, the finite time densities of right and left arrows must always be equal. We first consider the case when  $\kappa$ is odd for both CCA and FCA or $\kappa$ is even with $r=\lfloor \kappa/2 \rfloor$ for FCA. In this case, there are no blockades for all times. Suppose that fixation occurs. Then we claim that that all sites synchronize eventually, and hence the system clusters. To see this, suppose a site $x$ fixates. Then $x$ is not excited after some time, say $t_{0}$. If $x+1$ does not have the same color as $x$ at time $t_{0}$, then we must have $d\xi_{t_{0}}(x)< 0$ since otherwise $x$ is excited by $x+1$ at time $t_{0}$, a contradiction. Then 
					$x+1$ becomes the same color as $x$ after a finite time getting excited by $x$ and possibly by $x+2$ as well. Once $x+1$ is synchronized with $x$, $x+1$ cannot get any more excitation since otherwise $x$ will also get excited. This shows that $x+1$ fixates as well and will eventually be synchronized with $x$. By induction, it follows that all sites fixate and synchronize with $x$, as desired. 
					
					Next, keep the assumption that $\kappa$ is odd and suppose fluctuation occurs. Then by symmetry, the density of right and left arrows are always the same and positive (if it ever equals zero, we have fixated and it goes back to the previous case). Then every arrow must eventually be destroyed, since it will eventually encounter an opposing arrow and arrows do not cross each other without mutual annihilation. Thus arrow densities tend to zero and we deduce clustering. 
					
					Lastly, we assume $\kappa=2m$  for some integer $m\ge 2$ and  $r=m-1$ for CCA and FCA. 
					This case requires more care since there could be blockades, which are precisely the edges with color difference $m$ in this case. Note that in the arrow dynamics in both systems, blockades turn into stacks of $m-1$ opposing arrows when collided with a released arrow. By Lemma \ref{lem:FCA_undercoupled_flipping_transient1} there is no spontaneous formation of blockades after time $\kappa/2$, meaning that once a blockade at an edge is destroyed at time time $\ge m$, then that edge will never have a blockade thereafter. See the simulation in Figure \ref{fig:4CCA_FCA}. 
					\begin{figure}[h!]
						\centering
						\includegraphics[width=1\textwidth]{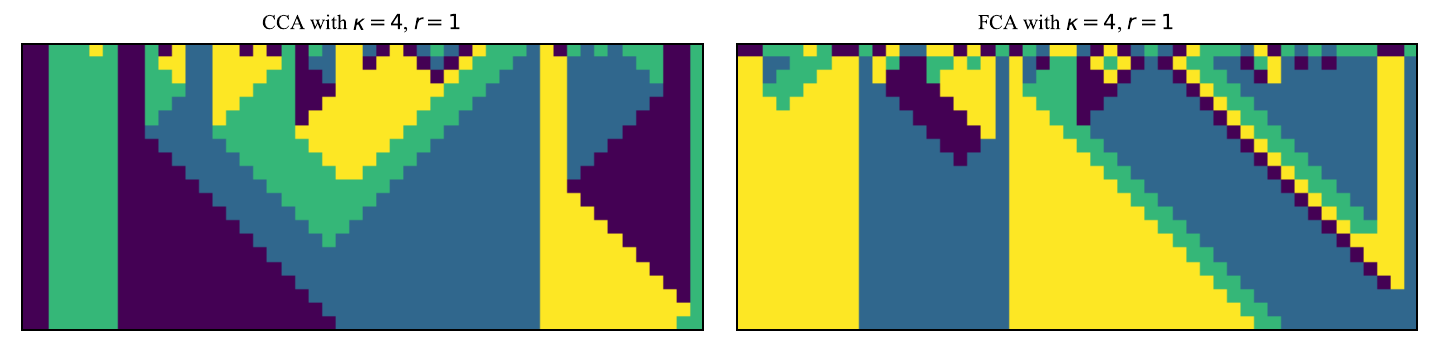}
						\vspace{-0.5cm}
						\caption{ Simulation of critical 4-color critical CCA and FCA ($r=1$) on 50 nodes in $\mathbb{Z}$ with 20 iterations. Blockades (edges with color difference 2) reflect arrows and disappear upon collision. Time goes from top to bottom. For CCA configurations at every iteration are shown, while for FCA only the ones at times $4 t$ are shown.
						}
						\label{fig:4CCA_FCA}
					\end{figure}
					

					It suffices to show that, assuming blockades are almost surely annihilated, that arrows are also almost surely annihilated. Suppose for contradiction that there is a positive density of arrows surviving forever. By ergodicity, there exist infinitely many arrows that survive forever. There are two cases:
					Suppose some of the surviving arrows hit finitely many blockades (hence flip their direction finitely many times). By ergodicity, then, there is a positive density of these and equal densities of those that eventually move right and eventually move left. But then there must be some opposing pair that meet, a contradiction. Therefore, any arrow that survives forever must hit infinitely many blockades.

					Consider two consecutive arrows that survive forever and hit infinitely many blockades. After some long finite time, all of the blockades between them will be removed. If the arrows are pointing away from each other, wait until the next time one of them turns around. Now they are pointing in the same direction or towards each other. If they are headed in the same direction, there are no blockades between them, so the next time the other arrow hits a blockade they will be pointing towards each other. Their next interaction is annihilation, a contradiction.
					Therefore, all arrows should be destroyed. This argument assumes that arrows never cross each other without annihilation and that an arrow keeps its identity when it is flipped after colliding with a blockade. The arrow dynamics for CCA was constructed to guarantee these properties. 
				\end{proof}
				
				We now deduce Theorem \ref{thm:main1} \textbf{(ii)} from the result above. 
				
				\begin{proof}[\textbf{Proof of clustering for critical CCA and FCA; Theorem \ref{thm:main1}}]
					By Proposition \ref{prop:clustering_blockade} it suffices to show that all blockades are eliminated eventually. The only case that there could be blockades in the statement is when $\kappa$ is even, $r=(\kappa/2)-1$, and the dynamics is CCA.

					We first show that the system must fluctuate. Suppose for a contradiction that the system fixates. We first consider CCA. Then the system must fixate into any set $\{a,a+m\}$ of non-interacting colors for $a=0,1,\dots,m-1$. By symmetry under color shift, these fixation events have the same probability and they are also invariant under spatial shift. Using the fact that initial colors of sites are independent,  by Kolmogorov's 0-1 law, we can deduce that these fixation events have probability 0 or 1 simultaneously. It follows that they must all have probability zero, so fluctuation follows. (The same argument does not rule out fixation for FCA as the system is not invariant under color shift.) 
					
					It remains to show that the blockade density must tend to zero. Indeed, if there exists a positive blockade density in the limit, then there are infinitely many blockades that survive forever by ergodicity. The arrows in between these forever-surviving blockades must eventually all be eliminated since otherwise a surviving arrow will eventually kill such surviving blockade, a contradiction. Thus arrows are annihilated almost surely but blockades survive with a positive probability, implying fixation. This contradicts that the system must fluctuate, as we have just established above. 
				\end{proof}
				
				We remark on a particular consequence of the clustering result above. It shows that the 4-color CCA with $r=1$ clusters:
				\begin{corollary}
					4-color CCA with $r=1$ clusters. 
				\end{corollary}
				
				\noindent The original $4$-color CCA with $r=1$ was known to fluctuate by Fisch \cite{fisch1990cyclic} in 1991 and was conjectured to cluster. Surprisingly, this result was proven only in 2021 by Hellouin de Menibus and Borgne \cite{hellouin2021asymptotic}, where the authors characterized the limiting measure with general product initial measure (not necessarily uniform). Analogous system with asynchronous update, known as the 4-color cyclic particle system, was conjectured to cluster by Bramson and Griffeath in 1989 \cite{bramson1989flux}. This was confirmed affirmatively by Foxall and Lyu \cite{foxall2018clustering} in 2018.

				\subsection{Transience of flipping and rate of clustering}
				\label{sec:critical_2}
				In the remainder of this section, we consider one of the three critical systems: FCA with $r=\lfloor \kappa/2 \rfloor$ for both $\kappa$ even and odd; and CCA with $\kappa$ even and $r= \kappa/2$. In these cases, we would like to establish the clustering rate 
				\begin{align}
					\P(\xi_{t}(0)\ne \xi_{t}(1)) = \Theta(t^{-1/2}). 
				\end{align}
				We use two different arguments to establish the upper and lower bounds of the asymptotic above. In this section, we establish the upper bound by relating it to the persistence probability of a random walk with correlated increments. 

				\vspace{2mm}

				
				First, we show that stack flipping is a transient phenomenon that does not occur after time $5\kappa$. This is stated in Proposition \ref{prop:finiteflipping} below. 
				\begin{prop}\label{prop:finiteflipping}
					Let $\xi_{t}$ denote one of the three critical systems: FCA with $r=\lfloor \kappa/2 \rfloor$ for both $\kappa$ even and odd; and CCA with $\kappa$ even and $r= \kappa/2$. Suppose $\xi_{0}$ is arbitrary. Then there exists $0\le t_{0}\le 5\kappa$ such that there is no stack flipping after time $t_{0}$. In particular, we can take $t_{0}=0$ for CCA. 
				\end{prop}

				\begin{proof}
					
					Suppose $\kappa=2m+1$ for some $m\in \mathbb{N}$. Note that $r=m$ by the hypothesis. The assertion \textbf{(i)} holds trivially for CCA since for $\kappa$ odd and $r=\lfloor \kappa/2\rfloor$, there is no stack flipping. To show \textbf{(i)} for the FCA, suppose for the sake of contradiction that the stack at an edge $(x,x+1)$ flips at time $t\ge m$. We may assume without loss of generality that $d\xi_{t}(x)=m$ and $x$ is excited at time $t$. This requires $\xi_{t}(x)=a$ for some $a\in [1,r]$ and $x$ has another neighbor $x-1$ such that $\xi_{t}(x-1)=0$. Backtracking $a$ steps, and noting the fact that $x$ cannot get excited by from either side from times $t-a+1$ through $t-1$, we deduce $\xi_{t-a+1}(x+1)=m+1$ and $\xi_{t-a}(x)=0$. It follows that $\xi_{t-a}(x+1) = m$ or $m+1$, but both choices are impossible; If $\xi_{t-a}(x+1) = m$, then $x+1$ is excited by $x$ at time $t-a$ (recall $m=r$) so that $\xi_{t-a+1}(x+1) = m$, contradicting $\xi_{t-a+1}(x+1)=m+1$; If $\xi_{t-a}(x+1) = m+1$, then $\xi_{t-a}(x+1)=m+2\ne m+1$, a contradiction. See \eqref{eq:FCA_odd_transient_flipping} below.
					\begin{align}\label{eq:FCA_odd_transient_flipping}
						\begin{matrix}
							\textup{time}: & t & t-1 & & & t-a+1 &  t-a  \\
							\hline
							x+1 : & a+m & a+m-1 & a+m-2  &\cdots & m+1 & b \\
							x :  & a  &  a-1 & a-2  & \cdots & 1 & \mathbf{0}   \\
							x-1 : &  \mathbf{0} & \kappa-1 & \cdots
						\end{matrix}
					\end{align}

					Now suppose $\kappa=2m$ for some integer $m\ge 2$. Recall that $r=m$ is the interaction range in this case. Suppose, for the sake of contradiction, that we have the stack at the edge $(x,x+1)$ flips at time $t=N$ for some integer $N> 10m$. Then $|d\xi_{t}(x)|=m-1$. We first claim that $\{\xi_{t_{1}}(x),\xi_{t_{1}}(x+1)\}=\{0,m\}$ for some $t_{1}\ge 9m$.

					Without loss of generality, we may assume $d\xi_{N}(x)=m-1$ and $x$ is excited by a blinking neighbor $x-1$ at time $t=N$. So $\xi_{N}(x-1)=0$, $\xi_{N}(x)=a$ for some $a\in [1,m]$, and $\xi_{N}(x+1)=a+m-1$. As before, none of $x\pm 1$ blink during time interval $(N-a,N]$, so backtracking $a-1$ iterations gives $\xi_{N-(a-1)}(x)=1$ and $\xi_{N-(a-1)}(x+1)=m$. This yields $\xi_{N-a}(x)=0$ and $\xi_{N-a}(x+1)=m$ since $x$ excites $x+1$ at time $t_{1}:=N-a$. This shows the claim.

					Next, we show that as we further backtrack the dynamics, the same ``opposite'' local configuration $\{\xi_{t}(x),\xi_{t}(x+1)\}=\{0,m\}$ appears with period $m+1$. Since $t_{1}\ge 9m\ge 6(m+1)$, this yields that the same pattern appears at least seven times during $[0,t_{1}]$. We will then obtain a contradiction. For simplicity, we work with $\kappa=6$ (so $m=3$) case but generalizing to any even $\kappa$ is immediate.

					In the following backtracking tables, time increases from right to left and each column gives a local configuration on sites $x-1,x,x+1,x+2$. The first column shows the opposite local configuration at $t=t_{1}$: 
					\begin{equation}
						\begin{matrix}
							\textup{time}: & t_{1} & && t_{1}-3& \\
							\hline
							x+2	:&	b_{1} & b_{2} & b_{3} & b_{4} \\
							x+1	:&	3 & & & *&  \\		
							x	:&	\mathbf{0} & 5 & 4 & 3 \\
							x-1	:&	a_{1} & a_{2} & a_{3} & a_{4} 
						\end{matrix}
					\end{equation}
					
					Observe that $*=0$ in the above table contradicts the transition rule since $3=m=r$. 
					Hence we can further backtrack the dynamics as follows:
					\begin{equation}
						\begin{matrix}
							\textup{time}: & t_{1} & &&&  t_{1}-4 & &&  t_{1}-7\\
							\hline			
							x+2	:&	b_{1} & b_{2} & b_{3} & b_{4} & b_{5} & b_{6} & b_{7} & b_{8}\\
							x+1	:&	3 & & && \textbf{0}& 5 & 4 & 3  \\
							x	:&	\mathbf{0} & 5 & 4 & 3 & 3 \\
							x-1	:&	a_{1} & a_{2} & a_{3} & a_{4} & a_{5} & a_{6} & a_{7} & a_{8}
						\end{matrix}
					\end{equation}
					with $0\in \{b_{2},b_{3},b_{4}\}$.  But note that $b_{2}\ne 0$ since otherwise we have $b_{4}=4$ and $b_{5}=3$, which contradict the transition rule since then $\xi_{2}$ would excite $\xi_{3}$ at the column containing $b_{5}$ so that $b_{4}=3$. Notice that the column containing $a_{5}$ is identical to the first one with the role of $\xi_{1}$ and $\xi_{2}$ exchanged. Hence the opposite local configuration appears with period $4=m+1$.

					Since $t_{1}\ge 6(m+1)$, we further extend our backtracking table as follow:
					\begin{equation}
						\begin{matrix}
							\textup{time}: & t_{1} & && &  t_{1}-4 &&&& t_{1}-8 &&&& t_{1}-12\\
							\hline			
							x+2:	&	b_{1} & b_{2} & \mathbf{b_{3}} & \mathbf{b_{4}} & b_{5} & b_{6} & b_{7} & b_{8} & b_{9} & b_{10} & \mathbf{b_{11}} & \mathbf{b_{12}} & b_{13}\\
							x+1:	&	3 & & && \textbf{0}& 5 & 4 & 3  & 3 &&&& \mathbf{2}\\
							x:	&	\mathbf{0} & 5 & 4& 3 & 3 & & & & \textbf{0} & 5 & 4 & 3 & 3\\
							x-1	:&	a_{1} & a_{2} & a_{3} & a_{4} & a_{5} & a_{6} & a_{7} & a_{8} & a_{9} & a_{10} & a_{11} & a_{12} & a_{13}
						\end{matrix}
						\label{oppositebacktracking}
					\end{equation}
					where $0\in \{b_{3},b_{4}\}$ and $0\in \{b_{11},b_{12}\}$ by repeating the same argument. In fact, none of the four combinations is possible. Suppose $b_{3}=0$. Then $b_{5}=4$, so $\xi_{2}$ does not excite $\xi_{3}$ at the column containing $b_{5}$. But then by the degree condition, $\xi_{3}$ must have been excited by at most once during the transition $b_{3}\leftarrow b_{11}$, which is impossible as this transition requires two excitations. Thus $b_{4}=0$. Since we  have enough history to apply the same argument to $b_{11}$ and $b_{12}$, we conclude $b_{12}=0$. But then by a similar argument the transition $0=b_{4}\leftarrow b_{12}=0$ is impossible. This shows the assertion.  
				\end{proof}

				Recall the classical mass transport principle. 
				
				\begin{prop}[Mass transport principle]\label{prop:MTP}
					Define a non-negative random variable $Z(a,b)$ for integers $a,b\in \mathbb{Z}$ such that its distribution is diagonally invariant under translation, i.e., for any integer $d$,  $Z(a+d,b+d)$ has the same distribution as $Z(a,b)$. Then for each $a\in \mathbb{Z}$,
					\begin{equation}
						\E \left[ \sum_{b\in \mathbb{Z}} Z(a,b) \right] = \E \left[ \sum_{b\in \mathbb{Z}} Z(b,a)\right].
					\end{equation}
				\end{prop}
				
				\begin{proof}
					Using  linearity of expectation and translation invariance of $\E[Z(a,b)]$, we get 
					\begin{equation}
						\E \left[ \sum_{b\in \mathbb{Z}} Z(a,b) \right] =  \sum_{b\in \mathbb{Z}} \E [Z(a,b)]  = \sum_{b\in \mathbb{Z}} \E [Z(2a-b,a)] = \sum_{b\in \mathbb{Z}} \E [Z(b,a)] = \E  \left[ \sum_{b\in \mathbb{Z}} Z(b,a) \right].
						\nonumber\end{equation}  
				\end{proof}
				
				Let $t_{0}$ be as in Proposition \ref{prop:finiteflipping}, the time after which there is no stack flipping. In terms of the arrow process, $d\xi_{t_{0}+\tau}(x)>0$ for some $\tau\ge 0$ means that there is a right arrow at time $t_{0}+\tau$ on the edge $(x,x+1)$. Since there is no spontaneous creation of arrows and each arrow maintains its identity until annihilation, there must be some arrow at a prior time $t_{0}$ somewhere that made it to the edge $(x,x+1)$ at time $t_{0}+\tau$. 
				In this way, we can relate the probability of local disagreement with the survival probability of arrows at time $t_{0}$. We make this precise in Lemma \ref{lem:disagreement_survival} below.

				\begin{lemma}[Local disagreement and particle survival]\label{lem:disagreement_survival}
					Let $(\xi_{t})_{t\ge 0}$ and $0\le t_{0} \le 5\kappa$ be as in Proposition \ref{prop:finiteflipping}. Let $\vec{N}_{t_{0}}(\tau)$ denote the number of right arrows on the edge $(0,1)$ at time $t_{0}$ that are alive at time $t_{0}+\tau$.  Then for any $\tau\ge 0$, 
					\begin{align}\label{eq:disagreement_persistence}
						\frac{2}{\lfloor \kappa/2 \rfloor} \,\E\left[ \vec{N}_{t_{0}}(\tau) \right]  \le   \P\left( \xi_{t_{0}+\tau}(0)\ne \xi_{t_{0}+\tau}(1)  \right) \le 2 \,\E\left[ \vec{N}_{t_{0}}(\tau) \right].
					\end{align}
				\end{lemma}
				
				\begin{proof}
					For each integers $a,b\in \Z$ and $t\ge 0$, define a nonnegative random variable $Z(a,b)$ as 
					\begin{align}
						Z_{t}(a,b) := \#\left( \begin{matrix} \textup{right arrows on the edge $(a,a+1)$ at time $t_{0}$} \\ \textup{that are on $(b,b+1)$ at time $t_{0}+\tau$} \end{matrix} \right). 
					\end{align}
					By the translation invariance of the edge particle system with initial particles given by $d\xi_{t_{0}}$, $Z_{t}(a+c,b+c)$ has the same distribution as $Z_{t}(a,b)$ for any integer $c$. Hence by Proposition \ref{prop:MTP}, we have 
					\begin{align*}
						\E\left[ \vec{N}_{t_{0}}(\tau) \right] &= \E\left[\sum_{b\in Z} Z_{\tau}(0,b)\right]  =\E\left[\sum_{b\in Z} Z_{\tau}(b,0)\right] \nonumber \\
						& = \E\left[\#\left( \textup{right arrows on $(0,1)$ at time $t_{0}+\tau$} \right) \right]\nonumber\\
						&\ge \P\left( \textup{$\exists$ an right arrow on $(0,1)$ at time $t_{0}+\tau$}  \right)\nonumber\\
						& = \P\left(|d\xi_{t_{0}+\tau}(0,1)|>0 \right)/2,
					\end{align*}
					where the third equality uses the fact that there is no stack flipping after time $t_{0}$ (Prop. \ref{prop:finiteflipping}). This establishes the second inequality in \eqref{eq:disagreement_persistence}. For its first inequality, simply note that 
					\begin{align*}
						& \E\left[\#\left( \textup{right arrows on $(0,1)$ at time $t_{0}+\tau$} \right) \right]  \\
						&\qquad \le \lfloor\kappa/2 \rfloor \P\left( \textup{$\exists$ an right arrow on $(0,1)$ at time $t_{0}+\tau$} \right)  \\
						&\qquad =  \frac{\lfloor \kappa/2 \rfloor}{2}  \P\left(|d\xi_{t_{0}+\tau}(0,1)|>0 \right),
					\end{align*}
					where the inequality follows since every stack can hold at most $\lfloor\kappa/2 \rfloor$ arrows. This is enough to conclude. 
				\end{proof}

				To proceed further, we define the \textit{associated additive process} $\{S_{n}\}_{n\ge 0}$ for the critical excitable media we consider in this section by $S_{0}=0$ and 
				\begin{align}\label{eq:def_associated_RW}
					S_{n} = \sum_{i=0}^{n-1}d\xi_{t_{0}}(i,i+1).
				\end{align}
				Note that $S_{n}$ counts the number of right arrows at time $t_{0}$ in the interval $[0,n]$ against the number of left arrows at time $t_{0}$ in the same interval. As we will prove 
				
				The following lemma translates the temporal event of seeing a particle on a particular edge at time $t$ into the spatial event of having positive partial sums at prior times; this is analogous to a duality relation in the graphical construction of an interacting particle system, in which one tracks the origin of a particle backwards in time. 
				\begin{lemma}[Particle survival and persistence of partial sums]\label{lem:RWcondition}
					Let $(\xi_{t})_{t\ge 0}$ and $0\le t_{0} \le 5\kappa$ be as in Proposition \ref{prop:finiteflipping}. Let $\{S_{n}\}_{n\ge 0}$ be as in \eqref{eq:def_associated_RW}. For CCA (resp., FCA), every released arrow moves in its prescribed direction until annihilation with speed between $v^{-}=1/\lfloor \kappa/2 \rfloor$ (resp., $1/\lfloor \kappa/2\rfloor (\kappa+1)$) and $v^{+}=1$ (resp., $1/\kappa$). 
					Furthermore, suppose a right arrow is on the edge $(0,1)$ at time $t=t_{0}$ and it is alive at time $t_{0}+\tau$ for some $\tau\ge 0$. Then 
					\begin{align}\label{eq:CCA_partial_sum}
						S_{z} \ge -\floor{\kappa/2} + 2 \qquad \text{for all $1\le z\le 2\lfloor \tau v^{-}\rfloor$}. 
					\end{align}
					Conversely, if 
					\begin{align}\label{eq:CCA_partial_sum2}
						S_{z} \ge 1 \qquad \text{for all $1\le z\le 2\lceil \tau v^{+} \rceil  $}, 
					\end{align}
					then all right arrows on the edge $(0,1)$ at time $t_{0}$ are alive at time $t_{0}+\tau$.

				\end{lemma}

				\begin{proof}
					We first consider the critical CCA. Note that in the arrow dynamics for the CCA with odd $\kappa$ and $r=\lfloor \kappa/2 \rfloor$, there is no blockade (and hence no stack flipping) and every nonzero stack of arrows is active at all times. The size of the stacks is at most $\lfloor \kappa/2\rfloor$. 
					Thus each right arrow is released again at least once in every $\lfloor \kappa/2 \rfloor$ iterations at at most once every iteration. This shows that each right particle moves ballistically to the right until annihilation with speed between $1/\lfloor \kappa/2 \rfloor$ and 1. 
					
					Furthermore, suppose a right arrow with label $\ell$ is on the edge $(0,1)$ at time $t=t_{0}$ and it is alive at time $t_{0}+\tau$ for some $\tau\ge 0$. Since this arrow moves to the right with speed at least $1/\lfloor \kappa/2 \rfloor$, it must be on the edge $(x,x+1)$ at time $t_{0}+\tau$ for some $x\ge \tau/\lfloor \kappa/2\rfloor$. If the right arrow with label $\ell$ is at height 0 in the stack on the edge $(0,1)$ at time $t_{0}$, then in order for it to survive until time $t_{0}+\tau$, the number of right arrows in any interval $[0,z]\subset [0,2 \lfloor \tau/v^{-} \rfloor]$ at time $t_{0}$ must exceed that of the left arrows in the same interval at the same time, because otherwise the right arrow on the edge $(0,1)$ at time $t=t_{0}$ will be annihilated before time $\tau$ since left arrows also move with speed at least $v^{-}$, a contradiction. Thus the partial sums in \eqref{eq:CCA_partial_sum} in this case must always be at least 1 for all $1\le z\le 2\lfloor \tau/v^{-} \rfloor$.  In general, the height of the right arrow with label $\ell$ at time $t_{0}$ is at most $\lfloor \kappa/2 \rfloor$, so the partial sums must be at least $-\lfloor \kappa/2\rfloor$ for all $1\le z \le 2\lfloor \tau/\lfloor \kappa/2 \rfloor \rfloor$. 
					
					Conversely, suppose \eqref{eq:CCA_partial_sum2} holds. Pick a right arrow on the edge $(0,1)$ at time $t_{0}$ and let $\ell$ be its label. Suppose for a contradiction that this right arrow is annihilated at time  before $t_{0}+\tau$. This annihilation must occur on some edge $(x,x+1)$ for some $x \in [ \lfloor \tau v^{-} \rfloor, \lceil \tau v^{+} \rceil ]$ since right arrows move with speed between $v^{-}$ and $v^{+}$ without flipping. Since the same holds for the left arrows, the unique left arrow that annihilates the right arrow with label $\ell$ must be on some edge $(y,y+1)$ for some $y\in [x+\lfloor \tau v^{-} \rfloor, x+\lceil \tau v^{+} \rceil]\subset [0,2\lceil \tau v^{+} \rceil]$. Since arrows do not flip after time $t_{0}$ and no two arrows cross without mutual annihilation, it follows that $S_{y}<0$, which contradicts \eqref{eq:CCA_partial_sum2}. This show the assertion for the critical CCA. 
					
					Next, we show the assertion for the critical FCA.  Recall that we assume $r=\lfloor \kappa/2 \rfloor$ for critical FCA. Suppose a right arrow with label $\ell$ is released from an edge $(y,y+1)$ at some time $t$ and eventually released again from the edge $(y+1,y+2)$ without being annihilated. Note that the stack on the edge $(y+1,y+2)$ must have either no particles or some right arrows at time $t$. Suppose the former. We then have $\xi_{t}(y)=0$, $\xi_{t}(y+1)=\xi_{t}(y+2)\in [1, r]$, and $\xi_{t+1}(y+1)=\xi_{t}(y)$ since $y+1$ is excited by $y$ at time $t$. Hence the next time that $y+1$ blinks (i.e., has color 0) 
					is between $t+r+2$ and $t+\kappa$. 
					Since by Proposition \ref{prop:finiteflipping} there is no stack flipping after time $t_{0}$, there could be at most $\lfloor \kappa/2 \rfloor-1$ right arrows  on the edge $(y+1,y+2)$ at time $t$. Note that $y+1$ can be excited only by $y$ until all right arrows  on $(y+1,y+2)$ have been released, so $y+1$ blinks at least once in every $\kappa+1$ iterations until then. Hence, it takes at most $\lfloor \kappa/2 \rfloor(\kappa+1)$ iterations for $(y+1,y+2)$ to release the right particle with label $\ell$. This shows that right particles travel ballistically with speed between $1/\lfloor \kappa/2\rfloor (\kappa+1)$ and $1/\kappa$. The rest of the statement for the critical FCA can be shown similarly as before. 
				\end{proof}
				
				We are now ready to prove the clustering rate in Theorem \ref{thm:main1} \textbf{(ii)}. The key idea of the proof is to view the additive process  $S_{n}$ in \eqref{eq:def_associated_RW} as a `Markov additive process' of a certain underlying Markov chain. Namely, note that the increments $d\xi_{t_{0}}(x)$ of $S_{n}$ can be thought of as a functional of an underlying spatial Markov chain, whose state at `site $x$' is the color tuple 
				\begin{align}\label{eq:def_color_tuple}
					\vec{Y}_{x}:= (\xi_{0}(y); y\in [x-q, x+1 +q])\in \Z_{\kappa}^{2t_{0}+2}. 
				\end{align}
				Note that $(\vec{Y}_{x})_{x\in \Z}$ is ergodic,  being irreducible on finite state space, and it has a finite correlation length: $\vec{Y}_{x}$ and $\vec{Y}_{y}$ are independent as soon as $|x-y|>2t_{0}+1$. Then $d\xi_{t_{0}}(x) = g(\vec{Y}_{x})$ for some functional $g$, since $d\xi_{t_{0}}(x)$ is determined by $\xi_{t_{0}}(x)$ and $\xi_{t_{0}}(x+1)$, which are in turn determined by the initial colors in $\vec{Y}_{x}$ through $t_{0}$ iterations of the CCA/FCA transition map. Partial sums with increments of  a functional of an underlying Markov chain are called \textit{Markov additive functional}. Hence the associated additive process $S_{n}$ in \eqref{eq:def_associated_RW} is indeed a Markov additive process with underlying Markov chain $(\vec{Y}_{x})_{x\ge 0}$.

				A crucial quantity in the analysis of the Markov additive process $S_{n}$ is the limiting variance  $\gamma_{g}^{2}$ defined by 
				\begin{align}\label{eq:limiting_variance}
					\gamma_{g}^{2}:= \Var( d\xi_{t_{0}}(0) ) + 2\sum_{j=1}^{\infty} \textup{Cov}(d\xi_{t_{0}}(0) ,\, d\xi_{t_{0}}(j)). 
				\end{align}
				According to \cite[Thm. 17.5.3]{meyn2012markov}, we can write it as 
				\begin{align}
					\gamma_{g}^{2}=\lim_{n\rightarrow \infty} n^{-1}\E[S_{n}^{2}]. 
				\end{align}
				Recalling that $d\xi_{t_{0}}(x)$ and $d\xi_{t_{0}}(y)$ are independent as soon as $|x-y|>2t_{0}+1$ and observing that $S_{2t_{0}+1}$ has positive finite variance, it follows that  $\gamma_{g}^{2}$ is also positive and finite. We denote $\gamma_{g}:=\sqrt{\gamma_{g}^{2}}$.

				\begin{proof}[\textbf{Proof of the clustering rate in Theorem \ref{thm:main1} \textbf{(ii)}}]
					
					Let $t_{0}$ be as in Proposition \ref{prop:finiteflipping}. Let $v^{\pm}$ be as in Lemma \ref{lem:RWcondition}. Then by Lemmas \ref{lem:disagreement_survival} and \ref{lem:RWcondition}, we have 
					\begin{align*}
						\P\left( \xi_{t}(0)\ne \xi_{t}(1)  \right) 
						&\le 2 \lfloor \kappa/2 \rfloor  \P(d\xi_{t_{0}}(0)>0)\, \P\left( \begin{matrix}\textup{right arrow on $(0,1)$ at time $t_{0}$} \\ \textup{is alive at time $t$} \end{matrix} \,\bigg|\, d\xi_{t_{0}}(0)>0\right)   \\
						&\le 2 \lfloor \kappa/2 \rfloor  \P( S_{1}\ge -\lfloor \kappa/2\rfloor+2,\dots, S_{\lfloor t/v^{-} \rfloor} \ge -\lfloor \kappa/2\rfloor+2),
					\end{align*}
					where the first inequality uses the fact that $\lVert d\xi_{0} \rVert_{\infty}\le \lfloor \kappa/2 \rfloor$. Similarly, we have 
					\begin{align*}
						\P\left( \xi_{t}(0)\ne \xi_{t}(1)  \right) 
						&\ge 2  \P(d\xi_{t_{0}}(0)>0)\, \P\left( \begin{matrix}\textup{right arrow on $(0,1)$ at time $t_{0}$} \\ \textup{is alive at time $t$} \end{matrix} \,\bigg|\, d\xi_{t_{0}}(0)>0\right)   \\
						&\ge 2  \P(d\xi_{t_{0}}(0)>0)\, \P( S_{1}\ge1 ,\dots, S_{2 \lceil t/v^{+} \rceil } \ge 1). 
					\end{align*}
					Thus for the desired asymptotic for the clustering rate, it suffices to show that 
					\begin{equation}\label{eq:survivalprob_for_thm2}
						\mathbb{P}\left( S_{1}\ge b,\cdots,S_{t}\ge b \right) = \Theta(t^{-1/2})
					\end{equation}
					for any fixed constant $b\in \mathbb{R}$. 
					
					For the critical CCA ($\kappa$ odd and $r=\lfloor \kappa/2 \rfloor$), we can take $t_{0}=0$ in Proposition \ref{prop:finiteflipping} so in this case the increments $(d\xi_{0}(x))_{x\ge 0}$ for the associated additive process $S_{n}$ are i.i.d. with mean zero and uniformly bounded values (by $\lfloor \kappa/2\rfloor$). In this case, the probability in \eqref{eq:survivalprob_for_thm2} is called a \textit{persistence probability} of a random walk (see, e.g., \cite{bray2013persistence}) and it is in fact known that 
					\begin{align}
						\mathbb{P}\left( S_{1}\ge b,\cdots,S_{t}\ge b \right) \sim C t^{-1/2}
					\end{align}
					for some constant $C>0$. Such an exact asymptotics of persistence probabilities for simple random walk was first obtained by Sparre Anderson \cite{andersen1953sums} in 1953 and later it was generalized by Feller \cite[Thm. XII.7.1]{feller1971introduction} to arbitrary random walks with i.i.d. increments satisfying some moment conditions. In particular, this justifies the  asymptotic in \eqref{eq:survivalprob_for_thm2}. 
					
					For the critical FCA, the situation is more delicate. Recall that in this case $t_{0}$ is not necessarily zero but can be taken to be at most $5\kappa$ (see Prop. \ref{prop:finiteflipping}). Hence now the increments $(d\xi_{t_{0}}(x))_{x\ge 0}$ for the associated random walk $S_{n}$ are (still identically distributed due to the translation invariance of the process but) are 
					not independent, since $d\xi_{t_{0}}(x)$ depends not only on the initial colors $\xi_{0}(x)$ and $\xi_{0}(x+1)$, but also on all initial colors $\xi_{0}(y)$ for $y\in [x-t_{0}, x+1 +t_{0}]$. This is because during the first $t_{0}$ iterations of `burn-in' period, sites interact with their neighbors to update their colors. Thus classical results on random walk persistence probabilities cannot be used to deduce \eqref{eq:survivalprob_for_thm2} for the critical FCA. 
					
					Instead, we have noted above the proof that the increments $d\xi_{t_{0}}(x)$ for the additive process $S_{n}$ can be thought of as a functional of underlying Markov chain $(\vec{Y}_{x})_{x\ge 0}$. The asymptotics of the persistence probabilities of Markov additive functionals are obtained by Lyu and Sivakoff \cite{lyu2017persistence}. In particular, \cite[Thm. 2]{lyu2017persistence} yields 
					\begin{align}\label{eq:MAF_persistence1}
						\sum_{b\in \Z: |b|\le \lfloor \kappa/2 \rfloor}  \mathbb{P}\left( S_{1}\ge b,\cdots,S_{t}\ge b \right) \sim \frac{\gamma_{g}}{\sqrt{2\pi}} t^{-1/2}, 
					\end{align}
					where the quantity $\gamma_{g}^{2}$ is defined in \eqref{eq:limiting_variance}. To conclude, note that the finitely many persistence probabilities in \eqref{eq:MAF_persistence1} for various values of $b$'s satisfy linear relations which can be derived by standard first-step analysis (see \cite[Prop. 3.4]{lyu2017persistence}). Thus \eqref{eq:MAF_persistence1} implies \eqref{eq:survivalprob_for_thm2}, as desired. 
				\end{proof}


				\subsection{Tournament expansion and the excitation rate} 
				\label{sec:critical_3}
				In this section, we establish the $\sqrt{t}$ excitation rate of the critical excitable media stated in Theorem \ref{thm:main1} \textbf{(ii)} \eqref{eq:clustering_rate}. While we were able to use the the embedded arrow dynamics to obtain the clustering rate, we need a different technique to analyze the excitation rate.  To this effect, we introduce another comparison process for the critical excitable media from a different perspective, which we call the \textit{tournament expansion}. This technique was originally developed by Gravner, Lyu, and Sivakoff  \cite{gravner2016limiting} in order to better understand the 3-color CCA and GHM dynamics on arbitrary underlying graphs.

				We begin with an instructive toy model called a \textit{tournament process}. Instead of coloring nodes of a graph $G=(V,E)$  by using mod $\kappa$ colors, consider a map $\mathtt{Rk}_{0}:V\rightarrow \mathbb{Z}$ called a \textit{ranking} on $G$. The transition map from time $t$ to $t+1$ is given by 
				\begin{equation}
					\mathtt{Rk}_{t+1}(x) = \max \{\mathtt{Rk}_{t}(y)\,|\, y\in N(x)\cup\{x\}\}. 
				\end{equation}  
				In words, in each transition $\mathtt{Rk}_{t}\mapsto \mathtt{Rk}_{t+1}$, each node simultaneously copies the maximum rank among it and its neighbors. Observe that if $G$ is finite, then for any initial ranking $\text{Rk}_{0}$ on $G$ there is a global maximum, and each node will eventually adopt the global maximum. In general, locally maximum rank propagates with unit speed across the graph until it is overcome by a higher ranker. 
				
				The basic idea of relating a similar process to the FCA dynamics on $\mathbb{Z}$ is the following: we construct an accompanying tournament-like system where each site increases its rank if and only if it gets excited. To give a detailed construction, we denote
				\begin{equation}
					\int_{x}^{y} d\xi_{t}:= \sum_{x \le z < y} d\xi_{t}(z) .
				\end{equation} 
				for sites $x\le y$ in $\Z$ and times $t\ge 0$. If $x>y$, then we define $\int_{x}^{y} d\xi_{t} = -\int_{y}^{x} d\xi_{t}$. 

				
				Let $(\xi_{t})_{t\ge 0}$ and $0\le t_{0} \le 5\kappa$ be as in Proposition \ref{prop:finiteflipping}.  We define its \textit{tournament expansion} by a sequence of rankings $(\mathtt{rk}_{t})_{t\ge t_{0}}$ defined as follows: set the rank of the origin at time $t$ by
				\begin{equation}
					\mathtt{rk}_{t}(0)=\mathcal{E}_{t}(0)=\sum_{s=0}^{t-1} \one(\text{0 is excited at time $s$})
				\end{equation}
				for all $t\ge 0$, and then extend to all sites by a path integral of $-d\xi_{t}$:
				\begin{equation}\label{eq:def_tournament_expansion}
					\mathtt{rk}_{t}(x) -\mathtt{rk}_{t}(0) = -\int_{0}^{x} d\xi_{t}
				\end{equation}  
				The minus sign in front of the path integral reflects the fact that the direction of $d\xi_{t}$ agrees with the direction of excitation, and ranking should increase whenever getting excited. 
				
				
				Our first observation about tournament expansion is that all sites indeed increment their rank by 1 if and only if they get excited:

				\begin{prop}\label{excitation}
					Let $(\xi_{t})_{t\ge 0}$, $0\le t_{0}\le 5\kappa$, and $(\mathtt{rk}_{t})_{t\ge 0}$ be as before. 
					For any $x,y\in \mathbb{Z}$ and $t\ge t_{0}$, we have 
					\begin{equation}\label{flux}
						\int_{x}^{y} \,d\xi_{t+1} - \int_{x}^{y} \,d\xi_{t} = \one(\text{$x$ is excited at time $t$}) - \one(\text{$y$ is excited at time $t$}).
					\end{equation}	
					Furthermore, we have 
					\begin{equation*}
						\mathtt{rk}_{t+1}(x)-\mathtt{rk}_{t}(x) = \one(\text{$x$ is excited at time $t$}).
					\end{equation*}
				\end{prop}

				\begin{proof} 
					The second part of the assertion follows immediately from the first part, definition of rank of the origin, and the following identity 
					\begin{eqnarray*}
						\mathtt{rk}_{t+1}(x)-\mathtt{rk}_{t}(x) = [\mathtt{rk}_{t+1}(0)-\mathtt{rk}_{t}(0)]- \left[ \int_{0}^{x} d\xi_{t+1} - \int_{0}^{x} d\xi_{t} \right].
					\end{eqnarray*}
					Now we show the first part. We may assume $x<y$ since both  sides in (\ref{flux}) are antisymmetric under the exchange of $x$ and $y$. Observe that it suffices to show the assertion for $y=x+1$, since then
					\begin{align*}
						\sum_{z=x}^{y-1}(d\xi_{t+1}-d\xi_{t})(z) &= \sum_{i=0}^{k-1} \mathbf{1}(\text{$z$ is excited at time $t$})-\mathbf{1}(\text{$z+1$ is excited at time $t$})\\
						&= \mathbf{1}(\text{$x$ is excited at time $t$})-\mathbf{1}(\text{$y$ is excited at time $t$}),
					\end{align*} 
					where the left hand side equals that of (\ref{flux}). It remains to verify that for any $x\in \mathbb{Z}$, 
					\begin{equation}\label{eq:pf_flux}
						(d\xi_{t+1}-d\xi_{t})(x)  = \mathbf{1}(\text{$x$ is excited at time $t$})- \mathbf{1}(\text{$x+1$ is excited at time $t$})	 
					\end{equation}
					for any $x\in \mathbb{Z}$. The argument is based on the arrow process; namely, the left hand side of (\ref{eq:pf_flux}) equals the ``flux'' of particles, which we define it to be the net change in the number of right arrows minus left arrows on edge $(x,x+1)$ from time $t$ to $t+1$.  
					
					First assume that $x$ excites $x+1$ at time $t$, so $x$ blinks at time $t$ and is not excited at the same time, so the right hand side of (\ref{eq:pf_flux}) equals $-1$. The assumption yields that edge $(x,x+1)$ has only right arrows, and its bottom right arrow is released at time $t$; since $x$ is not excited at time $t$, there is no right arrow released from $(x-1,x)$, so $(x,x+1)$ loses one right arrow during the transition $\xi_{t}\mapsto \xi_{t+1}$. Then using the arrow dynamics, the left hand side of (\ref{eq:pf_flux}) also equals $-1$, as desired. Second, suppose that $x$ does not excite $x+1$, but $x+1$ is excited by $x+2$ at time $t$. If $x$ is also excited by $x-1$ at time $t$, then no particle leaves the edge $(x,x+1)$, but there is both an incoming left arrow from the right and an right arrow from the left, which annihilate each other. Thus the flux is zero as asserted. Lastly, suppose that $x+1$ is excited by $x+2$ but $x$ is not excited by $x-1$. Then $x$ is not excited at time $t$, so the right hand side of (\ref{eq:pf_flux}) equals $-1$. Indeed, no right arrow is released from $(x-1,x)$ and a left arrow is released from $(x+1,x+2)$ at time $t$, so either this incoming left arrow annihilates the bottom right arrow on $(x,x+1)$ if any, or occupies the empty queue at $(x,x+1)$; in both cases, the flux on $(x,x+1)$ is $-1$. This shows the assertion. 
				\end{proof}

				The previous observation immediately yields that the tournament expansion for the critical excitable media does follow a tournament-like time evolution where local maxima subsume neighboring sites. 
				
				\begin{prop}\label{prop:tourmanent}
					Let $(\xi_{t})_{t\ge 0}$, $0\le t_{0}\le 5\kappa$, and $(\mathtt{rk}_{t})_{t\ge 0}$ be as before. Fix $t\ge t_{0}$ and $x\in \mathbb{Z}$. For CCA, we have
					\begin{equation}
						\mathtt{rk}_{t+1}(x) = \mathtt{rk}_{t}(x)+\mathbf{1}\big( \text{$\mathtt{rk}_{t}(x+1)>\mathtt{rk}_{t}(x)$ or $\mathtt{rk}_{t}(x-1)>\mathtt{rk}_{t}(x)$}\big).
					\end{equation}	
					For FCA, we have 
					\begin{equation}
						\mathtt{rk}_{t+1}(x) = 
						\begin{cases}
							\mathtt{rk}_{t}(x)+1 & \text{if $\exists y\in \{x\pm 1\}$ s.t. $\xi_{t}(y)=0$ and $\mathtt{rk}_{t}(y)>\mathtt{rk}_{t}(x)$}		\\
							\mathtt{rk}_{t}(x) & \text{otherwise}	
						\end{cases}
					\end{equation}	
				\end{prop}  
				
				\begin{proof}
					By construction \eqref{eq:def_tournament_expansion}, 
					\begin{align}
						\mathtt{rk}_{t}(x)=\mathtt{rk}_{t}(x-1)+d\xi_{t}(x-1) = \mathtt{rk}_{t}(x+1)-d\xi_{t}(x).  
					\end{align}
					In the CCA dynamics, recall that $x$ is excited at time $t$ if and only if either $d\xi_{t}(x-1)>0$ or $d\xi_{t}(x)<0$.
					Hence the assertion for CCA follows from the above. Also, for the FCA dynamics, note that  $x$ is excited at time $t$ if and only if there exists a blinking (of color 0) neighbor $y$ of $x$ with $d\xi_{t}(y,x)>0$, where $d\xi_{t}(y,x)=d\xi_{t}(x-1)$ if $y=x-1$ and $d\xi_{t}(y,x)=-d\xi_{t}(x)$ if $y=x+1$. Hence the assertion follows similarly. 
				\end{proof}

				Now we are ready to prove the key lemma in this section. Recall that the quantity $\mathtt{rk}_{t}(0)=\mathcal{E}_{t}(0)$ is a ``temporal'' quantity in the sense that it depends on the history of trajectory $(\xi_{s})_{s\ge 0}$ up to time $t-1$. Define a ``spatial'' quantity $M_{t}(m)$ by the maximum rank at time $t$ in the $m$-ball centered at the origin, i.e., 
				\begin{equation}\label{M_t}
					M_{t_{0}}(m) = \max\{\mathtt{rk}_{t_{0}}(x)\,:\, |x|\le m\} .
				\end{equation} 
				The next observation relates these two quantities.

				\begin{lemma}\label{lem:numberofpulls}
					Let $(\xi_{t})_{t\ge 0}$, $0\le t_{0}\le 5\kappa$, and $(\mathtt{rk}_{t})_{t\ge 0}$ be as before. 
					Then for each $t\ge 0$, we have that for CCA,  
					\begin{equation}\label{eq:rateofrankincrease_general_CCA}
						M_{t_{0}}( \lfloor t / \lfloor \kappa/2 \rfloor \rfloor)\le \mathcal{E}_{t+t_{0}}(0) \le M_{t_{0}}(t),
					\end{equation}
					and for FCA, \begin{equation}\label{eq:rateofrankincrease_general_FCA}
						M_{t_{0}}\left(\left\lfloor \frac{t}{\lfloor \kappa/2 \rfloor(\kappa+1)}\right\rfloor \right) \le \mathcal{E}_{ t+t_{0}}(0) \le M_{t_{0}}(\lfloor t/\kappa \rfloor).
					\end{equation}
					Furthermore, for FCA with $\kappa=3$, we have 
					\begin{equation}\label{eq:rateofrankincrease_general_3FCA}
						\mathcal{E}_{3t+t_{0}}(0) =  M_{t_{0}}(t).
					\end{equation}
				\end{lemma}

				\begin{proof}
					
					Recall that $\mathcal{E}_{t}(0)=\mathtt{rk}_{t}(0)$ for all $t\ge 0$ by construction. Let $v^{+}=1$ for CCA and $v^{+}=1/\kappa$ for FCA as in Lemma \ref{lem:RWcondition}. By Proposition \ref{prop:tourmanent}, each site increments its rank only when in contact with neighbors with higher rank. This happens if and only if an arrow passes through each site. By Lemma \ref{lem:RWcondition}, the speed of particles are at most $v^{+}$. Hence by time $t_{0}+t$, the origin can adopt the highest rank within radius at most $t v^{+}$, which is $M_{t_{0}}(t v^{+})$. This gives the upper bounds on $\mathcal{E}_{t+t_{0}}(0)$ in the inequalities \eqref{eq:rateofrankincrease_general_CCA}   and \eqref{eq:rateofrankincrease_general_FCA}. 
					
					To show the lower bounds,  we first claim that for any $t\ge 0$, $x\in \mathbb{Z}$, we have 
					\begin{equation}
						\mathtt{rk}_{t+\lfloor \kappa/2\rfloor (\kappa+1)}(x)\ge \max \{ \mathtt{rk}_{t}(y)\,|\,  y\in \{x-1, x,x+1\}  \}.
					\end{equation}		
					In words, the rank of a site $x$ exceeds the maximum rank in its 1-ball $\lfloor \kappa/2\rfloor (\kappa+1)$ iterations ago. This gives a lower bound on the growth rate of ranks at all sites, and repeating this inequality for the origin and making a change of variable in time, the lower bounds in \eqref{eq:rateofrankincrease_general_CCA} and \eqref{eq:rateofrankincrease_general_FCA} follow. 
					
					It remains to verify the claim. We may assume without loss of generality that $y=x-1$. If $\mathtt{rk}_{t}(x)\ge \mathtt{rk}_{t}(x-1)$, then by the monotonicity of ranks in time we have 
					\begin{equation}
						\mathtt{rk}_{t+\lfloor \kappa/2\rfloor (\kappa+1)}(x)\ge \mathtt{rk}_{t}(x-1).
					\end{equation}	
					So we may assume $\mathtt{rk}_{t}(x)< \mathtt{rk}_{t}(x-1)$. On the one hand, by construction of the ranking $\mathtt{rk}_{t}$, the rank difference between adjacent sites is upper bounded by $\lVert d\xi_{t} \rVert_{\infty}\le \lfloor \kappa/2\rfloor $, so we have 
					\begin{equation}
						0 \le \mathtt{rk}_{t}(x-1) - \mathtt{rk}_{t}(x) \le \lfloor\kappa/2\rfloor.
					\end{equation}
					On the other hand, by the monotonicty of ranks, $\mathtt{rk}_{s}(x-1)\ge \mathtt{rk}_{t}(x-1)$ for all times $s\ge t$. Hence by Proposition \ref{prop:tourmanent}, for CCA, the rank of $x$ will increment by one at least until it catches up with $\mathtt{rk}_{t}(x-1)$, which would take at most $\lfloor \kappa/2 \rfloor$ iterations since the rank difference is at most $\lfloor \kappa/2 \rfloor$. Thus the rank of $x$ will become at least $\mathtt{rk}_{t}(x-1)$ by time $t+\lfloor \kappa/2 \rfloor$. Similarly, for FCA, $x$ will increment by one whenever $x-1$ has color 0 at least until it catches up with $\mathtt{rk}_{t}(x-1)$. Note that until this happens, $x-1$ is not excited by $x$, so it returns to 0 at least once in every $\kappa+1$ iterations. Thus the rank of $x$ will become at least $\mathtt{rk}_{t}(x-1)$ by time $t+\lfloor \kappa/2 \rfloor(\kappa+1)$. This shows the assertion. 
					
					Lastly, we show \eqref{eq:rateofrankincrease_general_3FCA}. Consider FCA with $\kappa=3$ and $r=1$. For this, we only need to improve the lower bound in \eqref{eq:rateofrankincrease_general_FCA}. Fix $t\ge t_{0}$ and suppose $\mathtt{rk}_{t}(x)< \mathtt{rk}_{t}(x-1)$. Then since $\kappa=3$, we have $\mathtt{rk}_{t}(x)= \mathtt{rk}_{t}(x-1)-1$. For this, we need to have $(\xi_{t}(x-1),\xi_{t}(x))\in \{(1,2), (2,0), (0,1)\}$. If the color pair is either $(2,0)$ or $(0,1)$, then $x$ excites within two iterations and hence $\mathtt{rk}_{t+2}(x)\ge \mathtt{rk}_{t}(x-1)$. If the color pair is $(1,2)$, then $x$ excites within three iterations unless $\xi_{t}(x-2)=0$ so that $x-1$ excites at time $t$ and hence $\xi_{t+1}(x-1)=1$. However, this corresponds to the stack flip event, which we have shown to never occur for all times $t\ge t_{0}$ in Proposition \ref{prop:finiteflipping}. Thus, in all cases, $\mathtt{rk}_{t+3}(x)\ge \mathtt{rk}_{t}(x-1)$. This shows 
					\begin{align}
						\mathtt{rk}_{t+3}(x)\ge \max\{ \mathtt{rk}_{t}(x-1), \, \mathtt{rk}_{t}(x), \,  \mathtt{rk}_{t}(x+1) \}. 
					\end{align}
					Iterating the above inequality shows $\mathcal{E}_{3t+t_{0}}(0) \ge  M_{t_{0}}(t)$. This completes the proof. 
				\end{proof}

				Now we give a proof to Theorem \ref{thm:main2}. 
				
				\begin{proof}[\textbf{Proof of Theorem \ref{thm:main2}.}] Fix $\kappa\ge 3$, and let $(\xi_{t})_{t\ge 0}$, $t_{0}\in [0,5\kappa]$, and $(\mathtt{rk}_{t})_{t\ge t_{0}}$ be as before. 
					Recall that $M_{t_{0}}(m)$ denotes the maximum rank at time $t_{0}$ within distance $m$ from the origin (see \eqref{M_t}). Fix an integer $s\in [t_{0}+1, t]$ and write 
					\begin{align}\label{eq:max_decomposition}
						M_{t_{0}}(t) = \mathtt{rk}_{t_{0}}(0)+ \max \left\{ \max_{-t\le x < -s} \int_{-s}^{x} d\xi_{t_{0}}, \, \max_{-s\le x \le s} \int_{-s}^{x} d\xi_{t_{0}}, \,  \max_{s< x \le t} \int_{s}^{x} d\xi_{t_{0}} \right\}.
					\end{align}
					Note that $\int_{s}^{x}d\xi_{t_{0}} = S_{x}-S_{s}$. Bounding the maximum from above by the sum and from below by one argument and it follows that 
					\begin{align}\label{eq:max_ranking_bds}
						\E\left[  \max_{s< x \le t} S_{x}-S_{s} \right]  \le  \E[M_{t_{0}}(t)] \le \E\left[ \max_{-s\le x \le s} \int_{0}^{x} d\xi_{t_{0}} \right] + 2 \, \E\left[  \max_{s< x \le t} S_{x}-S_{s} \right], 
					\end{align}
					where we have used that the first and the last integrals in \eqref{eq:max_decomposition} have the same distribution, which follows from the distribution of $d\xi_{t_{0}}$ being symmetric about zero. Now using \cite[Theorems 1 and 2]{lyu2017persistence}, we have 
					\begin{align}
						\E\left[ \max_{s< x \le t} S_{x} \right] \sim \gamma_{g} \sqrt{2 t/\pi}, 
					\end{align}
					where $\gamma_{g}^{2}\in (0,\infty)$   denotes the limiting variance of the edge increments introduced in \eqref{eq:limiting_variance}. Combining this with  \eqref{eq:max_ranking_bds}, it follows that 
					\begin{align}
						\gamma_{g} \sqrt{2 /\pi} \le   \liminf_{t\rightarrow\infty} \frac{\E[M_{t_{0}}(t)]}{\sqrt{t}}  \le   \limsup_{t\rightarrow\infty} \frac{\E[M_{t_{0}}(t)]}{\sqrt{t}} \le 2 \gamma_{g} \sqrt{2 /\pi}. 
					\end{align}
					From this and Lemma \ref{lem:numberofpulls}, the $\Theta(\sqrt{t})$ asymptotics of the expected excitation count $\mathcal{E}_{t}(0)$ in 
					\eqref{eq:excitation_rate} follows immediately. This shows \textbf{(i)}. 
					
					Next, in order to show \textbf{(ii)}, we claim that as $t\rightarrow \infty$,
					\begin{equation}\label{eq:max_ranking_claim_lim}
						\frac{1}{\gamma_{g}\sqrt{t}} M_{t_{0}}(t) \overset{d}{\longrightarrow} \max\left\{ \max_{0\le u \le 1} B_{u} ,\, \max_{0\le v \le 1} B_{v}'  \right\}=:\mathtt{M},
					\end{equation}
					where $(B_{u})_{0\le u \le 1}$ and $(B_{v})_{0\le v \le 1}$ are independent standard Brownian motions. Note that by the reflection principle, for any $s\ge 0$, we have 
					\begin{equation}
						\mathbb{P}(\mathtt{M}\ge s) = 1-(1-2\mathbb{P}(Z\ge s))^{2} =  4\mathbb{P}(Z\ge s) \mathbb{P}(Z\le s)
					\end{equation}
					where $Z\sim N(0,1)$ is a standard normal random variable.

					Note that our underlying Markov chain $(\vec{Y}(x))_{x\in \mathbb{Z}}$ (see \eqref{eq:def_color_tuple}) is irreducible on a finite state space. 
					Since $\gamma_{g}^{2}\in (0,\infty)$, the functional central limit theorem holds for the Markov additive process $S_{n}$ in \eqref{eq:def_associated_RW} (e.g., see \cite[Ch. 17]{meyn2012markov}). 
					Hence it follows that 
					\begin{align*}
						\frac{1}{\gamma_{g} \sqrt{t}} \max_{s\le x \le t} \int_{s}^{x} d\xi_{t_{0}} \overset{d}{\longrightarrow} \max_{0\le u \le 1} \, B_{u}, 
					\end{align*}
					Since the underlying chain $(\vec{Y}_{x})_{x\in \Z}$ is reversible, this functional convergence holds for both directions $x\rightarrow \pm \infty$. In particular, we have 
					\begin{equation}\label{eq:max_convergence}
						\max\left\{\frac{1}{\gamma_{g} \sqrt{t}}\max_{s\le x \le t} \int_{s}^{x}d\xi_{t_{0}}, \frac{1}{\gamma_{g} \sqrt{t}}\max_{-t\le x \le -s} \int_{-s}^{x}d\xi_{t_{0}} \right\}\overset{d}{\longrightarrow} \max\left\{\max_{0\le u \le 1}B_{u}, \max_{0\le v \le 1}B_{v}'\right\}=\mathtt{M}
					\end{equation}
					as $t\rightarrow \infty$. Note that the two integrals above are independent by the finite correlation length and $s\ge t_{0}+1$. Hence the limiting Brownian motions are independent. Finally, (\ref{eq:max_decomposition}) yields that the difference between $M_{t_{0}}(t)$ and $\gamma_{g}\sqrt{t}$ times the left hand side of (\ref{eq:max_convergence}) is bounded. Hence Slutzky's theorem yields the desired diffusive limit of $M_{t_{0}}(t)$ as asserted in the claim.
					
					Now the stochastic domination in \eqref{eq:excitation_rate_gen} follows immediately from the claim \eqref{eq:max_ranking_claim_lim} and Lemma \ref{lem:numberofpulls}. Next, suppose $\kappa=3$. Then Lemma \ref{lem:numberofpulls} yields 
					\begin{align*}
						\mathcal{E}_{t+t_{0}}(0) &= M_{t_{0}}(t) \qquad \textup{for CCA}, \\
						\mathcal{E}_{3t+t_{0}}(0) &= M_{t_{0}}(t) \qquad \textup{for FCA}. 
					\end{align*}
					Hence by \eqref{eq:max_ranking_claim_lim}, we deduce 
					\begin{align}
						\frac{1}{\gamma_{g}\sqrt{t}} \mathcal{E}_{t}(0) \overset{d}{\longrightarrow} \begin{cases} \mathtt{M}  & \textup{for CCA} \\
							\frac{1}{\sqrt{3}} \mathtt{M}  & \textup{for FCA},
						\end{cases}
					\end{align}
					where for FCA we also used the monotonicity of $\mathcal{E}_{t}(0)$ in $t$. To conclude, note that $\gamma_{g}^{2}=2/3$ for CCA and $8/27$ for FCA (see \cite[Sec. 4]{lyu2017persistence}). This shows 
					\eqref{eq:3FCA_exact_tournament_eq2}, as desired. This shows \textbf{(ii)}. 
				\end{proof}

				\vspace{3mm}
				
				\section{Overcoupled CCA and FCA} 
				\label{sec:overcoupled_regime}
				
				In this section, we prove the first part of Theorem \ref{thm:main1} \textbf{(iii)} that every site gets excited at least at a linear rate almost surely when $r>\lfloor \kappa/2 \rfloor$ for CCA and $\lfloor \kappa/2 \rfloor< r < \kappa-1$  for FCA.

				\subsection{Stable Periodic Objects and linear excitation rate}
				\label{sec:overcoupled_regime_gen}
				We first consider the overcoupled CCA. For a concrete example, suppose $\kappa=6$ and $r=3$ (see Fig. \ref{fig:CCA_FCA_sim2}). Suppose an edge $(x,x+1)$ has colors $(0,3)$. Then both sites excite each other simultaneously, so after one iteration, they become colors $(1,4)$. The same holds indefinitely, and since the involved sites $x$ and $x+1$  are already excited at a maximal rate, there is no way that their local dynamics can be affected by external configurations. Thus, any edge with color difference $3$ acts as a stable object, advancing its colors by one every iteration. This is an example of \textit{Stable Periodic Orbits} (SPOs), which was first used by Fisch, Gravner, and Griffeath  to analyze two-dimensional CCA \cite{fisch1991cyclic}. 
				Our key observation for overcoupled CCA is that, despite being in one dimension, over-coupling allows SPOs, which will be the source of linear excitation of all sites. This is enough to show Theorem \ref{thm:main1} \textbf{(iii)(a)}.

				\begin{proof}[\textbf{Proof of Theorem \ref{thm:main1} \textbf{(iii)(a)}}] 
					
					Consider the overcoupled CCA. First assume that $\kappa$ is even. In the initial configuration, there is a positive density of edges with color difference $\kappa/2$. By the ergodic theorem,  
					there are infinitely many such edges almost surely. In the overcoupled regime, edges with a color difference $\kappa/2$ maintains the same color difference, where sites in those edges increment their colors by one at every iteration, regardless of configurations outside. 
					
					
					Next, we show that every site eventually increments its color by one every iteration. Let $A$ be the set of sites that are eventually increasing their colors by one at every iteration. Suppose for a contradiction that $A^c$ is  nonempty. Consider a site $v$ at the boundary of $A^c$ and has a neighbor $w$ in $A$. Then there are infinitely many times that $v$ does not increase its color by one but $w$ increases its color by one every iteration. It follows that there exist a time such that $\xi_t(w)=\xi_t(v)+1$. This in turn implies that $\xi_{t+1}(v)=\xi_t(v)+1$ and $\xi_{t+1}(w)=\xi_t(w)+1$ and, by induction, they remain locked in phase, which shows $v\in A$, a contradiction.

					The case of odd $\kappa$ with $r\geq \lceil \frac{\kappa}{2}\rceil$ is similarly dealt with by considering a SPO comprised by a site in state $0$ followed by $\frac{\kappa-1}{2}$. The gap between these two sites remains fixed below $r$ so they both update at a maximal rate (every time step) and hence remain locked in phase.
					
					Lastly, note that since every site eventually excites every iteration, we have 
					\begin{align}
						\lim_{t\rightarrow\infty} \frac{\mathcal{E}_{t}(0)}{t} = 1. 
					\end{align}
					This shows the assertion for the overcoupled CCA.
				\end{proof}



				
				Next, we consider the overcoupled FCA. Unlike in CCA, in FCA even with overcoupling, there is no finite configuration that have maximal excitation internally and cannot be perturbed by external configuration. To see this, note that in FCA each site can excite their neighbors (when color 0) at most once in $\kappa$ iterations, so on $\Z$ which has degree 2, every site can be excited at most twice in every $\kappa$ iterations. Hence for any configuration on a finite subset of sites $\Omega_{0}$, there will be a time that a boundary site in $\Omega_{0}$ is not excited at some color in $[1,r]$, and one can always place a 0 next to it from outside of $\Omega_{0}$ so that we exert additional excitation that is not supported internally by $\Omega_{0}$. 
				
				However, fortunately, the overcoupled FCA on $\Z$ admits a SPO-like object, where a $(\kappa+1)$-periodic `stable core' is protected by `quasi-stable' buffers. More specifically, length-7 color string $rrr0rrr$ is such an object. We establish its properties in the next lemma. 
				
				\begin{lemma}[Quasi-SPO for over-coupled FCA]\label{lem:FCA_SPO}
					Consider FCA with $r>\lceil \kappa/2 \rceil$. Suppose the color configuration on an interval $[x-3,x+3]$ at some time $t$ is $rrr0rrr$. Then the configuration $rr0rr$ on $[x-2,x+2]$ is $(\kappa+1)$-periodic while the colors on sites $x\pm 3$ can switch between $r$ and $r+1$ in period $\kappa+1$. 
				\end{lemma}
				
				\begin{proof} 
					
					For readability we denote color $r$ in bold in this proof. Without loss of generality we test color $\mathbf{r}$ on one edge and $\mathbf{r+1}$ on the other. It will be manifest in the proof that the edges $\{x-3,x+3\}$ are buffer sites and don't interact with the rest of the pattern.

					First we observe the evolution of the middle 5 sites. The main observation is they can never get pulls from the edge sites $\{x-3,x+3\}$ since the 0's there occur when the middle 5 are not "pullable". Furthermore, $\{x-1\}$ get pulled exactly once at time $t$ and otherwise all 5 advance at maximal rate. 
					
					Next, we turn our attention to the buffer sites at the edges. For and edge with color $\mathbf{r+1}$, say we have $X_t(x-3)=r+1$, it advances at maximal rate up until time  $t+(\kappa-r)$, when it is pulled by $\{x-2\}$ so necessarily we will have $X_{t+(\kappa-r+1)}(x-3)=1$. In the $r$ time steps to $t+(\kappa+1)$ it might get a single pull or none whatsoever, since it cannot be coming from the inside neighbor, arriving at colors $\mathbf{r}$ and $\mathbf{r+1}$ respectively. 
					
					As for an edge with color $\mathbf{r}$, say $X_t(x+3)=r$, there are two possibilities. In the first, showcased below, $X_t(x+3)=r$ does not receive a pull at time t, so it would advance at maximal rate to reach $X_{t+(\kappa-r+1)}(x-3)=1$. Then, as argued in the previous paragraph, it would reach colors $\mathbf{r/r+1} $ at time $t+(\kappa+1)$.

					\hspace*{-9mm}\vbox{{\small
							\begin{align*}
								\begin{matrix}
									\textup{time}: & t & & & & &t+(\kappa-r+1) & &  & & &t+(\kappa+1)  \\ 
									\hline
									\vspace{1mm}\\
									x+3 :   &  \mathbf{r} & r+1 & r+2 & \hdots& 0& \mathbf{1}& * &  *&*  & *&   \mathbf{r/r+1}   \\
									x+2 :  & \mathbf{r} & r+1 & r+2 & \hdots&  0 & \mathbf{1} & 1 &2&\hdots  &     r-1 &\mathbf{r}      \\
									x+1 :&\mathbf{r}&r&r+1 &\hdots &\kappa-1&\mathbf{0}&1  &2&\hdots  &  r-1&\mathbf{r} \\
									x :   &  \mathbf{0} & 1 & 2 & \hdots & \kappa-r& \kappa-r+1& \kappa-r+1 &\kappa-r+2  & \hdots &\kappa-1 & \mathbf{0}  \\
									x-1 :  & \mathbf{r} & r & r+1 & \hdots &  \kappa-1 & \mathbf{0} & 1 &2  &\hdots  & r-1  &\mathbf{r}    \\
									x-2 :&\mathbf{r}&r+1&r+2 & \hdots &0&\mathbf{1}&1  &2&\hdots  &  r-1&\mathbf{r}\\
									x-3 :   &  \mathbf{r+1} & r+2 & r+3 & \hdots& 1& \mathbf{1}& *  & * & * &  *&\mathbf{r/r+1}   \\
								\end{matrix}
							\end{align*}
					}}

					Otherwise, $X_t(x+3)=r$ does indeed get pulled at time $t$ and advances at maximal rate thereafter to reach $X_{t+(\kappa-r+1)}(x+3)=0$. In the remaining $r$ steps to $t+(\kappa+1)$ it must advance at maximal rate since the next blinking state $0$ at $\{x+4\}$ must come at $t+(\kappa+1)$ at the earliest. This is because $\{x+4\}$ does receive at least one pull in that time interval, as it does so  at time $t+(\kappa-r+1)$ since $(\kappa-r+1)\leq r$.

					\hspace*{-9mm}\vbox{{\small
							\begin{align*}
								\begin{matrix}
									\textup{time}: & t & & & & &t+(\kappa-r+1) & &  & & &t+(\kappa+1)  \\
									\hline
									\vspace{1mm}\\
									x+3 :   &  \mathbf{r} & r & r+1 & \hdots& \kappa-1& \mathbf{0}& 1 &  2&\cdots  & r-1 &   \mathbf{r}  \\
									x+2 :  & \mathbf{r} & r+1 & r+2 & \hdots&  0 & \mathbf{1} & 1 &2&\hdots  &     r-1 &\mathbf{r}  \\
									x+1 :&\mathbf{r}&r&r+1 &\hdots &\kappa-1&\mathbf{0}&1  &2&\hdots  &  r-1&\mathbf{r} \\
									x :   &  \mathbf{0} & 1 & 2 & \hdots & \kappa-r& \mathbf{\kappa-r+1}& \kappa-r+1 &\kappa-r+2  & \hdots &\kappa-1 & \mathbf{0}   \\
									x-1 :  & \mathbf{r} & r & r+1 & \hdots &  \kappa-1 & \mathbf{0} & 1 &2  &\hdots  & r-1  &\mathbf{r}   \\
									x-2 :&\mathbf{r}&r+1&r+2 & \hdots &0&\mathbf{1}&1  &2&\hdots  &  r-1&\mathbf{r}\\
									x-3 :   &  \mathbf{r+1} & r+2 & r+3 & \hdots& 1& \mathbf{1}& *  & * & * &  *&\mathbf{r/r+1}  \\
								\end{matrix}
							\end{align*}
					}}
				\end{proof}

				Now we show the first part of Theorem \ref{thm:main1} \textbf{(iii)(b)}. 
				
				\begin{proof}[\textbf{Proof of the first part of Theorem \ref{thm:main1} \textbf{(iii)(b)}}] By the Borel-Cantelli lemma the $rrr0rrr$ pattern occurs infinitely often in the random initial configuration. Hence, by the Lemma \ref{lem:FCA_SPO} we have a positive density of $(\kappa+1)$ periodic sites. 
					
					Next, we would like to show that every site excites at a linear rate: 
					\begin{align}\label{eq:linear_excitation_pf}
						\liminf_{t\rightarrow\infty} \frac{1}{t}\mathcal{E}_{t}(0)>0 
					\end{align}
					and that every site is not eventually $\kappa$- or $(\kappa+2)$-periodic. 
					For this, we first argue that the origin will eventually have  periodic dynamics. If the origin belongs to an $rrr0rrr$ string, then we know it will be $(\kappa+1)$ periodic. Otherwise, we can find the nearest such patterns in the initial configuration, so there is an almost surely finite random interval $[-N,M]$ containing the origin, where at the both ends of this interval we have $rrr0rrr$ patterns initially. In Lemma  \ref{lem:FCA_SPO}, we have shown that the dynamics in the interior of such patterns are determined regardless of configurations outside. Thus, the dynamics on the interval $(-N,M)$ is independent of the sites outside this interval. In turn, FCA restricted on $(-N,M)$ is a finite-state deterministic dynamical system, which must be periodic. This shows that the origin must eventually have a periodic dynamics.
					
					Now that we know the origin is eventually periodic, either it is eventually $\kappa$-periodic and excites only finitely many times or it has larger period than $\kappa$ and it gets excited at a linear rate. Thus  we can conclude \eqref{eq:linear_excitation_pf} by showing that the origin cannot be $\kappa$-periodic.  Suppose not. If  its neighbor, say site 1, is not already $\kappa$ periodic eventually, then its phase will drift so that together with the origin they are both in the blinking state, i.e.  $\{\xi_{t}(0),\xi_{t}(1)\}=\{0,0\}$. We have two possibilities. If the origin's neighbor blinks again after $\kappa+1$ steps they will arrive at configuration $\{\xi_{t+\kappa+1}(0),\xi_{t+\kappa+1}(1)\}=\{1,0\}$.  The origin would get a pull at time $\kappa+1$, contradicting $\kappa$ periodicity. Likewise if it takes $\kappa+2$ for the origin neighbor's to blink. Since the origin updates at maximal rate, they will arrive at configuration $\{\xi_{t+\kappa+1}(0),\xi_{t+\kappa+1}(1)\}=\{2,0\}$ contradicting the origin's $\kappa$-periodicity. Indeed, the origin would also get pulled in this case since, in this regime,  $r>\lceil\frac{\kappa}{2}\rceil\geq 2$ for all $\kappa\geq 3$. 
					
					It remains to show that the origin cannot be eventually $(\kappa+2)$-periodic. This argument proceeds similarly by showing that a $(\kappa+2)$-periodic site must have both neighbors also $(\kappa+2)$-periodic. Suppose the origin is $(\kappa+2)$-periodic. If the site to its right is not eventually $(\kappa+2)$-periodic, then it's phase must drift backwards (relative to the origin) by 1 or 2 steps in each interblinking cycle.  Therefore, there must be a time at which  $\{X_{t}(0),X_{t}(1)\}=\{0,0\}$. Looking at the next interblinking period, we have two possibilities:

					\begin{enumerate}
						\item{\textbf{Case 1: Right neighbor blinks in $\kappa$ steps}}: For the origin to have an interblinking period of $\kappa+2$ it must receive an excitation from the right neighbor. By the assumption, this can only happen at time $t+\kappa$. But by then the origin has been excited at most once from the left so either $\xi_{t+\kappa}(0)=0$ or $\xi_{t+\kappa}(0)=\kappa-1$. The former case obviously contradicts $(\kappa+2)$ periodicity. In the latter observe that, since $r<\kappa-1$ the origin updates at time $t+\kappa$ and blinks at time $t+(\kappa_1)$ also contradicting $\kappa+2$ periodicity.



								

						\begin{align*}
							\begin{matrix}
								\textup{time}: & t +(\kappa+2) & & t+\kappa &   &  & &  &t &\\ 
								\hline
								1 : &   * & \mathbf{1}& \mathbf{0}&
								\cdots &&&1 & \mathbf{0} & \\
								0 :  & 1 & \mathbf{0} & \mathbf{\kappa-1} &\cdots &  & & 1&  \mathbf{0} & 
							\end{matrix}
						\end{align*}

						\item{\textbf{Case 2: Right neighbor blinks in $\kappa+1$ steps}}: This gives a similar contradiction. Since $X_t(1)$ doesn't blink until $X_t(0)=\kappa-1$ and $r<\kappa-1$, the origin receives no pull from the right. So it can be pulled at most once from the left, contradicting $\kappa+2$ periodicity.

								

						
						\begin{align*}
							\begin{matrix}
								\textup{time}: & t+(\kappa+2) &t+(\kappa+1) & &   &  & &  &t &\\ 
								\hline
								1 : &  \mathbf{1} & \mathbf{0} & \kappa-1&
								\cdots &&&1 & \mathbf{0}& \\
								0 :  & \mathbf{0} & \mathbf{\kappa-1} & \kappa-2  &\cdots &  & & 1&  \mathbf{0} & 
							\end{matrix}
						\end{align*}

					\end{enumerate}
					Therefore, any neighbor of the origin would eventually have be $\kappa$ or  $\kappa+2$-periodic respectively. By induction, all vertices are eventually $\kappa$-periodic (resp. $\kappa+2$-periodic), but this contradicts existence of a positive density of $(\kappa+1)$-periodic sites we established before. 
				\end{proof}

				\subsection{$(\kappa+1)$-periodicity of the weakly overcoupled FCA} 
				\label{sec:FCA_weakly_overcoupled}
				
				In this section, we consider the overcoupled FCA in the particular case $r=\lceil \frac{\kappa}{2} \rceil+1$ and prove the second part of Theorem \ref{thm:main1} \textbf{(iii)(b)}. We show that all lattice sites are eventually $\kappa+1$ periodic. The cases of even and odd $\kappa$ differ slightly. As usual,  we let $\xi_{t}:\mathbb{Z}\rightarrow \Z/\kappa \Z$ denote the color configuration at time $t$. We will require the following definition.

				\begin{define}[Defects]\label{def:defects}
					We say there is a \textit{defect} at time $t$ on the edge $(x,x+1)$ if:
					\begin{itemize}[itemsep=0.1cm]
						\item For odd $\kappa$:      $\{\xi_{t}(x), \xi_{t}(x+1)\}= \{0,r-1\}$ or $\{0,r\}$.
						
						\item For even $\kappa$:     $\{\xi_{t}(x), \xi_{t}(x+1)\}= \{0,r-2\}$, $\{0,r-1\}$ or $\{0,r\}$. 
					\end{itemize}
				\end{define}

				A crucial feature of defects in the `weakly overcoupled FCA' is that they cannot be created spontaneously. Defects present at the current time must have been preceded by defects at some prior times on the same edges. We establish this fact and characterize all possible ways that defects can arise. Note that the classification depends on the parity of $\kappa$.

				\begin{prop}\label{prop:defect_backtrack}
					Suppose there is a defect at time $t$ on edge $(x,x+1)$ for $t\ge r+1$. Then there must be a defect on edge $(x,x+1)$ at time either $t-(r-1)$, $t-r$ or $t-(r+1)$. More precisely, without loss of generality, suppose $\xi_{t}(x)=0$. Then the dynamics on $(x,x+1)$ during $[t-(r+1), t]$ must be of the following. Bold columns signify defects.

					\begin{description}
						
						\item[(i)] For the odd $\kappa$ case we have the following 4:  
						
						\begin{align}
							&\begin{matrix}
								\textup{time}: & t & & &   & &  &t-r & \\ 
								\hline 
								x : &  \mathbf{0}&\kappa-1 &\kappa-2&
								\cdots &  &r & \mathbf{r} & *\\
								x-1 :  & \mathbf{r-1} & * & *&\cdots &  & 1&  \mathbf{0} & *
							\end{matrix}
							\label{eq:defect_type1}
						\end{align}
						\begin{align}
							&  \begin{matrix}
								\textup{time}: & t & & &   &  & &  &t-r &\\ 
								\hline
								x : &  \mathbf{0} & \kappa-1 & \kappa-2&
								\cdots &&&r & \mathbf{r} & *\\
								x-1 :  & \mathbf{r} & r-1 & r-2 &\cdots &  & & 1&  \mathbf{0} & *
							\end{matrix}
							\label{eq:defect_type2}
						\end{align}
						
						\begin{align}
							&   \begin{matrix}
								\textup{time}: & t & & &     & & &  & &t-r -1\\ 
								\hline
								x : &  \mathbf{0} & \kappa-1 & \kappa-2&
								\cdots & &&r & r & \mathbf{r}\\
								x-1 :  & \mathbf{r} & * &*& \cdots &  & & &  1 & \mathbf{0}
							\end{matrix}
							\label{eq:defect_type3}
						\end{align}
						\begin{align}
							&        \begin{matrix}
								\textup{time}: & t & & &  &    & &  & &t-r -1\\
								\hline
								x : &  \mathbf{0} & \kappa-1 & \kappa-2&
								\cdots & &&r & r-1 & \mathbf{r-1} \\
								x-1 :  & \mathbf{r} & * & * &\cdots &  & & &  1& \mathbf{0}
							\end{matrix}
							\label{eq:defect_type4}
						\end{align}

						\item[(ii)] For the even $\kappa$ case we have the following 9:

						\begin{align}
							&\begin{matrix}
								\textup{time}: & t & &   &  &t-r+1 & \\ 
								\hline
								x : &  \mathbf{0} & \kappa-1 &
								\hdots & r & \mathbf{r} & *\\
								x-1 :  & \mathbf{r-2} & * &\hdots &    1&  \mathbf{0} & *
							\end{matrix}
							\label{eq:defect_type5}
						\end{align}
						\begin{align}
							&\begin{matrix}
								\textup{time}: & t & & &    &t-r+1 & \\ 
								\hline
								x : &  \mathbf{0} & \kappa-1 & 
								\hdots  &r & \mathbf{r} & *\\
								x-1 :  & \mathbf{r-1} & r-2 &\hdots &    1&  \mathbf{0} & *
							\end{matrix}
							\label{eq:defect_type6}
						\end{align}
						\begin{align}
							&\begin{matrix}
								\textup{time}: & t & & &  & &   &t-r & \\ 
								\hline
								x : &  \mathbf{0} & \kappa-1 & &
								\hdots & &r-1 & \mathbf{r-1} & *\\
								x-1 :  & \mathbf{r-1} & * &&\hdots &   & 1&  \mathbf{0} & *
							\end{matrix}
							\label{eq:defect_type7}
						\end{align}
						\begin{align}
							&\begin{matrix}
								\textup{time}: & t & & &  & & &  &t-r & \\ 
								\hline
								x : &  \mathbf{0} & \kappa-1 & &
								\hdots & &r&r & \mathbf{r} & *\\
								x-1 :  & \mathbf{r-1} & * &&\hdots & &*  & 1&  \mathbf{0} & *
							\end{matrix}
							\label{eq:defect_type8}
						\end{align}
						\begin{align}
							&\begin{matrix}
								\textup{time}: & t & & &  & &   &t-r & \\ 
								\hline
								x : &  \mathbf{0} & \kappa-1 & \kappa-2 &
								\hdots & &r-1 & \mathbf{r-1} & *\\
								x-1 :  & \mathbf{r} & r-1 & r-2&\hdots &   & 1&  \mathbf{0} & *
							\end{matrix}
							\label{eq:defect_type9}
						\end{align}
						\begin{align}
							&\begin{matrix}
								\textup{time}: & t & & &  & &&   &t-r & \\ 
								\hline
								x : &  \mathbf{0} & \kappa-1 & \kappa-2 &
								\hdots & &r&r & \mathbf{r} & *\\
								x-1 :  & \mathbf{r} & r-1 & r-2&\hdots &   &2& 1&  \mathbf{0} & *
							\end{matrix}
							\label{eq:defect_type10}
						\end{align}
						
						\begin{align}
							&\begin{matrix}
								\textup{time}: & t & &   & &&   &t-r-1 & \\ 
								\hline
								x : &  \mathbf{0} & \kappa-1 & \kappa-2 &
								\hdots & r-1&r-2 & \mathbf{r-2} & *\\
								x-1 :  & \mathbf{r} & * & *&\hdots &  * & 1&  \mathbf{0} & *
							\end{matrix}
							\label{eq:defect_type11}
						\end{align}
						\begin{align}
							&\begin{matrix}
								\textup{time}: & t & & &  & &  & &t-r-1 & \\ 
								\hline
								x : &  \mathbf{0} & \kappa-1 & \kappa-2 &
								\hdots & &r-1&r-1 & \mathbf{r-1} & *\\
								x-1 :  & \mathbf{r} & * & *&\hdots &  &* & 1&  \mathbf{0} & *
							\end{matrix}
							\label{eq:defect_type12}
						\end{align}
						\begin{align}
							&\begin{matrix}
								\textup{time}: & t & & &  & & &     &&t-r-1 & \\ 
								\hline
								x : &  \mathbf{0} & \kappa-1 & \kappa-2 &
								\hdots &&r&r &r-1 & \mathbf{r-1} & *\\
								x-1 :  & \mathbf{r} & * & *&\hdots &&  *& *& 1&  \mathbf{0} & *
							\end{matrix}
							\label{eq:defect_type13}
						\end{align}
					\end{description}
					Furthermore, there is no spontaneous birth of defects. Specifically, for each edge $e=(x,x+1)$, either there is no defect on $e$ for all times, or there are defects on $e$ up to some finite time and then there is no defect thereafter (`transient defect'), or there are defects on $e$ for infinitely many times ('recurrent defect'). 
				\end{prop}

				\begin{proof}  It suffices to show the second part of the statement. We begin with the odd case. First suppose that  $(\xi_{t}(x),\xi_{t}(x-1))=(0,r-1)$. For the first sequence in \eqref{eq:defect_type1} note that $x$ cannot get pulled by $x-1$ so, backtracking, we should find the blinking state at time $t-r$ or $t-r+1$. Also since $\kappa-(r-1)=r$, we know this must be the color of site $x-1$ at time $t-(r-1)$. If $x$ pulls $x-1$ at that time we will have $\xi_{t-(r-2)}(x-1)=r$ which yields a backtracking contradiction. Thus $x$ blinks at time $t-r$ yielding \eqref{eq:defect_type1}.

					Next suppose that $(\xi_{t}(x-1),\xi_{t}(x))=(0,r)$. Similarly, $x$ does not get pulled by $x-1$ during times $[t,t-r]$ so we should find the blinking state at time $t-r$ or $t-r-1$. In the former case, again since $\kappa-(r-1)=r$, we must have $\xi_{t-(r-1)}(x-1)=r$  resulting in \eqref{eq:defect_type2}, or else we have a backtracking contradiction as before. In the latter case, we have two options. We can arrive at $\xi_{t-(r-1)}(x)=r$ by having $\xi_{t-(r+1)}(x)=r$ as in \eqref{eq:defect_type3} or $\xi_{t-(r+1)}(x)=r-1$ which yields \eqref{eq:defect_type4} .

					As for the even case, we begin by considering $(\xi_{t}(x-1),\xi_{t}(x))=(0,r-2)$. Since $x-1$ does not receive a pull from $x$ during times $[t-r,t]$ we  have it must blink at times $t-r+2$ or $t-r+1$. The former, argued as before, would result in a backtracking contradiction leaving us with \eqref{eq:defect_type5}. 
					
					Next suppose, $(\xi_{t}(x-1),\xi_{t}(x))=(0,r-1)$. With no pull from $x$,  $x-1$ must blink at times $t-r+1$ or $t-r$. In the former case, we must have  $x-1$ updates at a maximal rate and $\xi_{t-r+1}(x)=r$ yielding \eqref{eq:defect_type6}. In the latter case we have two options. We can have $\xi_{t-r}(x)=r-1$  and $x$ receives one pull from $x-1$ at time $t-r$ as in \eqref{eq:defect_type7} or we have $\xi_{t-r}(x)=r$ and $x$ receives consecutive pulls from $x-1$ and $x+1$, respectively, at times $t-r$ and $t-r+1$ yielding \eqref{eq:defect_type8}.

					Finally, suppose $(\xi_{t}(x-1),\xi_{t}(x))=(0,r)$, we have 5 possibilities. Once again, since $x$ doesn't pull $x-1$, the former must blink at times $t-r$ or $t-r-1$. If it blinks at time $t-r$ we have two possibilities. As before, we can have again that  $\xi_{t-r}(x)=r-1$  and $x$ receives one pull from $x-1$ at time $t-r$ yielding \eqref{eq:defect_type9} or that $\xi_{t-r}(x)=r$ and $x$ receives consecutive pulls from $x-1$ and $x+1$ respectively at times $t-r$ and $t-r+1$, which results in \eqref{eq:defect_type10}.
					
					Otherwise, it blinks at time  $t-r-1$ and we have three possibilities. If $\xi_{t-r-1}(x)=r-2$, then $x$ gets a pull at that time step and must update maximally from then on. This gives us sequence \eqref{eq:defect_type11}. If  $\xi_{t-r-1}(x)=r-1$ we have two other possibilities since $x$ must receive at least one pull from $x+1$. This may happen at time $t-r$ yielding \eqref{eq:defect_type12}. Otherwise, it happens at time $t-r+1$ when $\xi_{t-r+1}(x)=r$  resulting in  \eqref{eq:defect_type13}. 
				\end{proof}

				From Proposition \ref{prop:defect_backtrack}, it seems that the combinatorics of defect dynamics is too complicated to handle, since the different local dynamics can be concatenated in all possible ways to generate long-term dynamics. Fortunately, if a recurrent defect is right next to an edge without a recurrent defect, then there is only one type of dynamics possible for the recurrent defect. This is shown in the proposition below. We relegate its proof to the following subsection.

				\begin{prop}\label{prop:defect_recurrence2}
					Suppose an edge $(x-1,x)$ has a recurrent defect and the next edge $(x,x+1)$ does not. Then the dynamics on $(x-1,x)$ must be as follows:

					\begin{enumerate}
						\item For odd $\kappa$ the composition of the sequence \eqref{eq:defect_type2} only, that is, 
						\begin{align}
							\begin{matrix}
								x : & \cdots  &  \mathbf{0} & \kappa-1 & \kappa-2 &\cdots & r & \mathbf{r} & r-1 & r-2 &\cdots & 1& \mathbf{0}  & \cdots  \\
								x-1 :& \cdots  & \mathbf{r} & r-1 & r-2 &\hdots & 1 &  \mathbf{0} & \kappa-1 & \kappa-2&\hdots&r & \mathbf{r} & \cdots
							\end{matrix}
						\end{align}

						\item For even $\kappa$ the composition of sequences \eqref{eq:defect_type6} and \eqref{eq:defect_type9} that is,
						\begin{align}
							\begin{matrix}
								x : & \cdots  &  \mathbf{0} & \kappa-1 & \kappa-2 &\cdots & r-1 & \mathbf{r-1} & r-1 & r-2 &\cdots & 1& \mathbf{0}  & \cdots  \\
								x-1 :& \cdots  & \mathbf{r} & r-1 & r-2 &\hdots & 1 &  \mathbf{0} & \kappa-1 & \kappa-2&\hdots&r & \mathbf{r} & \cdots
							\end{matrix}
						\end{align}

					\end{enumerate}
				\end{prop}

				We now deduce the second part of Theorem \ref{thm:main1}\textbf{(iii)(b)}. 
				

				\begin{proof}[\textbf{Proof of the second part of Theorem \ref{thm:main1} \textbf{(iii)(b)}}]
					Our goal now is to show that  for the `weakly ocercoupled' FCA with $r=\lceil \kappa/2 \rceil+1$, almost surely, every site becomes $(\kappa+1)$-periodic. 
					
					First, we claim that there is a positive density of edges that have no defect for all times or transient defects. To see this, there is a positive density of intervals in the initial configuration of constant color. If this interval is long enough, then the middle edge will not have a defect from time 0 through time 
					$r+1$. Then by Proposition \ref{prop:defect_backtrack}, this middle edge will never have a defect for all times. By translation invariance and the ergodic theorem, we can conclude the claim. 
					
					Second, note that there is a positive density of edges that have recurrent defects and both nodes are eventually $(\kappa+1)$-periodic. Indeed recall that $(\kappa+1)$-periodicity of the interior of the $rrr0rrr$ configuration shown in Lemma \ref{lem:FCA_SPO}. 
					
					Third, we show that every site in a recurrent defect must be eventually $(\kappa+1)$-periodic. Suppose for contradiction that there exists a positive density of edges with recurrent defects and at least one of the end points are not eventually $\kappa+1$ periodic. Then there is an edge $(x-1, x)$ which has a recurrent defect and either $x-1$ or $x$ is not eventually $\kappa+1$ periodic, and the edge $(x,x+1)$ either does not have a recurrent defect or it has recurrent defect such that both $x$ and $x+1$ are eventually $\kappa+1$ periodic. We will get a contradiction from both cases. First, suppose $(x,x+1)$ is not a recurrent defect. Then $(x-1,x)$ is eventually $\kappa+1$ periodic by Proposition \ref{prop:defect_recurrence2}, a contradiction. So we may assume $(x,x+1)$ is a recurrent defect and both $x$ and $x+1$ are eventually $(\kappa+1)$-periodic. Therefore, $x-1$ must be eventually \textit{not} $(\kappa+1)$-periodic. But by Proposition \ref{prop:defect_backtrack}, the dynamics on $(x-1,x)$ must be given by composition of \eqref{eq:defect_type2}. It follows that  $x-1$ is also eventually $\kappa+1$ periodic, a contradiction. 
					
					To summarize, we have shown that there is a positive density of $(\kappa+1)$-periodic recurrent defects and all other edges are not recurrent defects. Consider a contiguous interval $I$ of edges that are not recurrent defects, which is surrounded by $(\kappa+1)$-periodic recurrent defects. We claim that $(\kappa+1)$-periodicity will spread into the interval $I$. To see this, observe that the lack of defects on $I$ allows us to use arrow dynamics as we introduced in Section \ref{sec:arrow_dynamics}. Namely, defects correspond to stacks of arrows that flip back and forth automatically, releasing one arrow every time it flips. Since there is no defect on $I$ after a finite amount of time, all stacks of arrows in $I$ have size $<\lfloor \frac{\kappa-1}{2} \rfloor$ for all sufficiently large times, so arrows in $I$ move ballistically without flipping and annihilate upon collisions. Hence, all arrows initially in $I$ will have annihilated after a finite time and arrows injected into $I$ from the recurrent defects at both ends will keep colliding somewhere in the middle of $I$ (see Figure \ref{fig:CCA_FCA_sim1} for $r=3$ and \ref{fig:CCA_FCA_sim2} for $r=4$). Thus every site in $I$ must be either $(\kappa+1)$- or $(\kappa+2)$-periodic, depending on the timing that it is hit by the arrows from both sides. However, we have shown that it cannot be $(\kappa+2)$-periodic (first part of Thm. \ref{thm:main1}\textbf{(iii)(b)}) since $(\kappa+2)$-periodicity spreads to neighboring sites.  Thus every site in $I$ must be eventually $(\kappa+1)$-periodic. This is enough to conclude. 
				\end{proof}

				\subsection{Proof of Proposition \ref{prop:defect_recurrence2}}
				
				In this section we prove Proposition \ref{prop:defect_recurrence2}. In order to simplify the presentation, we will prove the statement for two concrete cases: $\kappa=5$ and $\kappa=6$. Generalization to other cases will be straightforward. 
				
				\begin{proof}[\textbf{Proof of Proposition \ref{prop:defect_recurrence2} for $\kappa=5$}]
					
					Suppose an edge $(x-1,x)$ has a recurrent defect and the next edge $(x,x+1)$ does not. Suppose $\kappa=5$ and $r=3$. Then there are four defects according to Proposition \ref{prop:defect_backtrack}:
					\begin{align}
						\begin{matrix}
							\textup{time}: & t & & & t-3 \\ 
							\hline
							x : &  \mathbf{0} & 4 & 3 & \mathbf{3} & *\\
							x+1 :  & \mathbf{2} & * & 1 &  \mathbf{0} & *
						\end{matrix}\label{eq:defect_type14} \\
						\begin{matrix}
							\textup{time}: & t & & & t-3 \\ 
							\hline
							x : &  \mathbf{0} & 4 & 3 & \mathbf{3} & *\\
							x+1 :  & \mathbf{3} & 2 & 1 &  \mathbf{0} & *
						\end{matrix} \label{eq:defect_type15}
						\\
						\begin{matrix}
							\textup{time}: & t & & & &t-4 \\ 
							\hline
							x : &  \mathbf{0} & 4 & 3 & 3 & \mathbf{3}\\
							x+1 :  & \mathbf{3} & * & * &  1 & \mathbf{0}
						\end{matrix}\label{eq:defect_type16}
						\\
						\begin{matrix}
							\textup{time}: & t & & & &t-4 \\ 
							\hline
							x : &  \mathbf{0} & 4 & 3 & 2 & \mathbf{2}\\
							x+1 :  & \mathbf{3} & * & * &  1 & \mathbf{0}
						\end{matrix}
						\label{eq:defect_type17}
					\end{align}
					Our goal is to show that the local dynamics on the edge $(x-1,x)$ with a recurrent defect must be given by the concatenation of \eqref{eq:defect_type15} with itself:
					\begin{align}
						\begin{matrix}
							x : & \cdots  &  \mathbf{0} & 4 & 3 & \mathbf{3} & 2 & 1 & \mathbf{0} &  4& 3 & \mathbf{3}  & \cdots  \\
							x-1 :& \cdots  & \mathbf{3} & 2 & 1 &  \mathbf{0} & 4 & 3 & \mathbf{3} & 2 &1 & \mathbf{0} &\cdots
						\end{matrix}
					\end{align}

					Suppose that there exists $M>0$ such that there is no defect on the edge $(x,x+1)$ for all times $t\ge M$. We first rule out \eqref{eq:defect_type17}, which will also rule out \eqref{eq:defect_type14}, since \eqref{eq:defect_type14} must be followed by \eqref{eq:defect_type17} according to Proposition \ref{prop:defect_backtrack}. 
					
					Suppose $(x-1,x)$ has a defect at some time $t\ge M+12$ and its dynamics backtracks as \eqref{eq:defect_type17} from time $t$. There are two cases to consider: $(\xi_{t}(x-1),\xi_{t}(x))=(3,0)$ or $(0,3)$. Suppose the former. Then by Proposition \ref{prop:defect_backtrack}, the local dynamics on $(x-1,x,x+1)$ during $[t-7, t+4]$ must be one of the following two cases:
					\begin{align}
						\begin{matrix}
							\textup{time}: & t & & & &t-4 & & & t-7 & & & t-10 \\ 
							\hline
							x+1 :&&&&&&a& b &  &  &  &  & c& d \\
							x :   &  \mathbf{0} & 4 & 3 & 2 & \mathbf{2}& * & 1 & \mathbf{0} & 4 & 3 & \mathbf{3} & e & f \\
							x-1 :  & \mathbf{3} & * & * &  1 & \mathbf{0} & 4 & 3 & \mathbf{3} & 2 & 1 & \mathbf{0} & 4 & 3  
						\end{matrix}
						\label{eq:defect_local_dyn_pf01}\\[5pt]
						\begin{matrix}
							\textup{time}: & t & & & &t-4 & & & t-7 & & & & t-11 \\ 
							\hline
							x+1 :&&&&&&a& b &  &  &  & c & &  \\
							x :   &  \mathbf{0} & 4 & 3 & 2 & \mathbf{2}& * & 1 & \mathbf{0} & 4 & 3 & 3 & \mathbf{3} & *   \\
							x-1 :  & \mathbf{3} & * & * &  1 & \mathbf{0} & 4 & 3 & \mathbf{3} & *&* & 1 & \mathbf{0} & 4   
						\end{matrix}
						\label{eq:defect_local_dyn_pf02}\\[5pt]
						\begin{matrix}
							\textup{time}: & t & & & &t-4 & & & t-7 & & & & t-11 \\ 
							\hline
							x+1 :&&&&&&a& b &  &  &  &  & c& d \\
							x :   &  \mathbf{0} & 4 & 3 & 2 & \mathbf{2}& * & 1 & \mathbf{0} & 4 & 3 & 2 & \mathbf{2} & e   \\
							x-1 :  & \mathbf{3} & * & * &  1 & \mathbf{0} & 4 & 3 & \mathbf{3} & *&* & 1 & \mathbf{0} & 4   
						\end{matrix}
						\label{eq:defect_local_dyn_pf03}
					\end{align}
					In all cases, $a\ne 0$ since backtracking two steps yields a contradiction. Since $x$ needs one pull during $[t-6,t-5]$ but $x-1$ does not pull $x$ during this period, it follows that $b=0$. In case~\eqref{eq:defect_local_dyn_pf01}, observe that $x$ does not pull $x+1$ during $[t-12,t-7]$, so $x+1$ is pulled at most once during this interval by $x+2$. If $x+1$ is not pulled by $x+2$ during this interval, then $c=0$; otherwise $d=0$, so  $0\in \{c,d\}$. If $c=0$, then $e=3$, which is a contradiction with there being no defect on $(x,x+1)$ at time $t-11\ge M$. Likewise, if $d=0$ we also arrive at a contradiction with there being no defect on $(x,x+1)$ at time $t-12\ge M$. In case~\eqref{eq:defect_local_dyn_pf02}, $x$ must be pulled by $x+1$ at time $t-10$, so $c=0$, which creates a defect on $(x,x+1)$ and leads to a contradiction. In case~\eqref{eq:defect_local_dyn_pf03}, as in case~\eqref{eq:defect_local_dyn_pf01}, $x+1$ is pulled at most once during $[t-12,t-7]$, so either $c=0$ or $d=0$ and $e=2$; either scenario leads to a defect on $(x,x+1)$ and a contradiction. This rules out the case $(\xi_{t}(x-1),\xi_{t}(x))=(3,0)$.
					
					Next, suppose $(\xi_{t}(x-1),\xi_{t}(x))=(0,3)$. Then by Proposition \ref{prop:defect_backtrack}, the local dynamics on $(x-1,x,x+1)$ during $[t-11, t]$ must be one of the following two cases:
					\begin{align}
						\begin{matrix}
							\textup{time}:  & t & & & &t-4 & & & t-7\\ 
							\hline
							x+1 : &&a&b&c& &&& &d&e&f& \\
							x :   &  \mathbf{3} & * & * & 1 & \mathbf{0} & 4 & 3 & \mathbf{3} & 2& 1 & \mathbf{0}  \\
							x-1 : & \mathbf{0} & 4 & 3 &  2 & \mathbf{2} & * & 1 &  \mathbf{0} & 4 &3  & *  
						\end{matrix}
						\label{eq:defect_local_dyn_pf1}\\[5pt]
						\begin{matrix}
							\textup{time}:  & t & & & &t-4 & & & t-7 \\ 
							\hline
							x+1 : &&a&b&c& &&& &d&e&f& \\
							x :   &  \mathbf{3} & * & * & 1 & \mathbf{0} & 4 & 3 & \mathbf{3} & *& * &1 & \mathbf{0}  \\
							x-1 : & \mathbf{0} & 4 & 3 &  2 & \mathbf{2} & * & 1 &  \mathbf{0} & 4 &3 & * & * 
						\end{matrix}
						\label{eq:defect_local_dyn_pf2}
					\end{align}
					In both cases, $0\in \{a,b,c\}$. Note that $b\ne 0$ by two-step backtracking. Also, $a\ne 0$ since otherwise there is a defect on the edge $(x, x+1)$ at time $t-1$, contrary to the hypothesis. Thus we must have $c=0$. For case \eqref{eq:defect_local_dyn_pf1}, $0\notin \{d,e\}$, so we must have $f=0$. This requires $x+1$ to be pulled twice during $[t-10,t-3]$, but it is not pulled by $x$ during this period, a contradiction. For case \eqref{eq:defect_local_dyn_pf2}, $d\ne 0$ since otherwise there is a defect on $(x,x+1)$ at time $t-8$, and $e\ne 0$ by two-step backtracking. So we must have $f=0$, but then this gives the same contradiction as in the previous case. This shows that \eqref{eq:defect_type17} (and hence \eqref{eq:defect_type14}) cannot appear on the edge $(x-1,x)$ after time $M+8$. 
					
					Next, we rule out \eqref{eq:defect_type16}. Suppose for a contradiction that the dynamics on $(x-1,x)$ during $[t-4,t]$ is given by \eqref{eq:defect_type16}. We consider two cases: $(\xi_{t}(x-1),\xi_{t}(x))=(3,0)$ or $(0,3)$. If it is the former, the local dynamics on $(x-1,x,x+1)$ during $[t-4, t]$ is given as 
					\begin{align}\label{eq:defect_local_dyn_pf22}
						\begin{matrix}
							\textup{time}:   & t & & & &t-4 \\ 
							\hline
							x+1 :&  & & & 0 & & \\
							x :   & \mathbf{0} & 4 & 3 & 3 & \mathbf{3}  \\
							x-1 : & \mathbf{3} & * & * &  1 & \mathbf{0} 
						\end{matrix}
						,
					\end{align}
					which contradicts the hypothesis that there is no defect on the edge $(x,x+1)$ after time $M$ and $t\ge M+7$. Next, assume $(\xi_{t}(x-1),\xi_{t}(x))=(0,3)$. By Proposition \ref{prop:defect_backtrack}, the sequence \eqref{eq:defect_type3} can be preceded by one of \eqref{eq:defect_type15}, \eqref{eq:defect_type16}, and \eqref{eq:defect_type17}. Since  we have ruled out \eqref{eq:defect_type17} and \eqref{eq:defect_local_dyn_pf22}, the local dynamics on $(x-1,x,x+1)$ during $[t-11, t]$ is given by one of the two cases below:
					\begin{align}
						\begin{matrix}
							\textup{time}:   & t & & & &t-4 &&&&&& t-10\\ 
							\hline
							x+1 :&  & a& b& c & & &&&  d &e \\
							x :   & \mathbf{3} & * & * &  1 & \mathbf{0} & 4 & 3 & \mathbf{3} & 2 & 1 & \mathbf{0} \\
							x-1 : & \mathbf{0} & 4 & 3 & 3 & \mathbf{3} & 2 &  1 & \mathbf{0} & 4 & 3 & \mathbf{3}
						\end{matrix}
						\label{eq:defect_local_dyn_pf3} \\[5pt]
						\begin{matrix}
							\textup{time}:   & t & & & &t-4  &&&&&&& t-11\\ 
							\hline
							x+1 :&  & a& b& c & & &&& d & e & f\\
							x :   & \mathbf{3} & * & * &  1 & \mathbf{0} & 4 & 3 & \mathbf{3} & * & * & 1 & \mathbf{0} \\
							x-1 : & \mathbf{0} & 4 & 3 & 3 & \mathbf{3} & 2 &  1 & \mathbf{0} & 4 & 3 & 3 & \mathbf{3}
						\end{matrix}
						\label{eq:defect_local_dyn_pf4}.
					\end{align}
					For both cases, we must have $0\in \{a,b,c\}$ and $0\notin \{a,b\}$ since that would create a defect on the edge $(x,x+1)$ at times $\ge t-3\ge M$. So $c=0$. Then for \eqref{eq:defect_local_dyn_pf3}, $0\in \{d,e\}$ since $x+1$ is not pulled by $x$ during $[t-9,t]$. But neither  $d$ nor $e$ can be $0$ for the assumed dynamics at $x$ during $[t-10,t-7]$. Hence \eqref{eq:defect_local_dyn_pf3} is impossible. For \eqref{eq:defect_local_dyn_pf4}, $0\in \{d,e,f\}$ but $0\notin \{d,e\}$ since that will create a defect on the edge $(x,x+1)$ at times $\ge t-9\ge M$. So $f=0$, but then $x+1$ must be pulled twice during $[t-9, t-3]$, which is impossible since it is not pulled by $x$ during this period. This rules out \eqref{eq:defect_type16}.
					
					At this point, we have shown that the sequences \eqref{eq:defect_type14}, \eqref{eq:defect_type16} and \eqref{eq:defect_type17} and their upside-down versions, are transient on the edge $(x-1,x)$. Hence the only recurrent sequence among the four types in Proposition \ref{prop:defect_backtrack} is \eqref{eq:defect_type15}. This shows the assertion.  
				\end{proof}

				\begin{proof}[\textbf{Proof of Proposition \ref{prop:defect_recurrence2} for $\kappa=6$}]
					
					Suppose an edge $(x-1,x)$ has a recurrent defect and the next edge $(x,x+1)$ does not. Suppose $\kappa=6$ and $r=4$. Then there are nine defects according to Proposition \ref{prop:defect_backtrack}:
					\begin{align}
						\begin{matrix}     \textup{time}:  & t &  & &t-3 \\ 
							\hline
							x : &  \mathbf{0} & 5 & 4 & \mathbf{4} \\
							x-1 : & \mathbf{2} & * & 1 &  \mathbf{0}
						\end{matrix}
						\label{eq:defect_type18}
					\end{align}

					\begin{align}
						\begin{matrix}
							\textup{time}:  & t &  & &t-3 \\ 
							\hline
							x : &  \mathbf{0} & 5 & 4 & \mathbf{4}  \\
							x-1 : & \mathbf{3} & 2 & 1 &  \mathbf{0}  
						\end{matrix}
						\label{eq:defect_type19}
					\end{align}

					\begin{align}
						\begin{matrix}
							\textup{time}:  & t &  && &t-4 \\ 
							\hline
							x : &  \mathbf{0} & 5 & 4 & 3 &\mathbf{3} \\
							x-1 : & \mathbf{3} & * & * &  1 &\mathbf{0} 
						\end{matrix}
						\label{eq:defect_type20}
					\end{align}

					\begin{align}
						\begin{matrix} \textup{time}:  & t &  && &t-4 \\ 
							\hline
							x : &  \mathbf{0} & 5 & 4 & 4 &\mathbf{4}  \\
							x-1 : & \mathbf{3} & * & * &  1 &\mathbf{0} 
						\end{matrix}
						\label{eq:defect_type21}
					\end{align}

					\begin{align}
						\begin{matrix} \textup{time}:  & t &  && &t-4 \\ 
							\hline
							x : &  \mathbf{0} & 5 & 4 & 3 &\mathbf{3}  \\
							x-1 : & \mathbf{4} & 3 & 2 &  1 &\mathbf{0} 
						\end{matrix}
						\label{eq:defect_type22}
					\end{align}

					\begin{align}
						\begin{matrix} \textup{time}:  & t &  && &t-4 \\ 
							\hline
							x : &  \mathbf{0} & 5 & 4 & 4 &\mathbf{4}  \\
							x-1 : & \mathbf{4} & 3 & 2 &  1 &\mathbf{0} 
						\end{matrix}
						\label{eq:defect_type23}
					\end{align}
					
					\begin{align}
						\begin{matrix} \textup{time}:  & t &  && &&t-5 \\ 
							\hline
							x : &  \mathbf{0} & 5 & 4 & 3 &3 &\mathbf{3} \\
							x-1 : & \mathbf{4} & * & * &  * &1 &\mathbf{0}
						\end{matrix}
						\label{eq:defect_type24}
					\end{align}
					
					\begin{align}
						\begin{matrix} \textup{time}:  & t &  && &&t-5 \\ 
							\hline
							x : &  \mathbf{0} & 5 & 4 & 3 &2 &\mathbf{2} \\
							x-1 : & \mathbf{4} & * & * &  * &1 &\mathbf{0}
						\end{matrix}
						\label{eq:defect_type25}
					\end{align}
					
					\begin{align}
						\begin{matrix} \textup{time}:  & t &  && &&t-5 \\ 
							\hline
							x : &  \mathbf{0} & 5 & 4 & 4 &3 &\mathbf{3} \\
							x-1 : & \mathbf{4} & * & * &  * &1 &\mathbf{0}
						\end{matrix}
						\label{eq:defect_type26}
					\end{align}
					Our goal is to show that the local dynamics on the edge $(x-1,x)$ with a recurrent defect should be  given by the concatenation of \eqref{eq:defect_type19} and \eqref{eq:defect_type22}:
					\begin{align}
						\begin{matrix}
							x : & \cdots  &  \mathbf{0} & 5 & 4 & \mathbf{4}&3 & 2 & 1 & \mathbf{0} &  5& 4 & \mathbf{4}  & \cdots  \\
							x-1 :& \cdots  & \mathbf{3} & 2 & 1 &  \mathbf{0} & 5& 4 & 3 & \mathbf{3} & 2 &1 & \mathbf{0} &\cdots
						\end{matrix}
					\end{align}

					We first exclude   \eqref{eq:defect_type25}  which also excludes \eqref{eq:defect_type18} since it would have to be preceded by it. 
					
					\begin{align}
						&   \begin{matrix}
							\textup{time}: & t & & & & &t-5 & &  &t-8 & &&& t-12 \\
							\hline
							x+1 :&&& &&&&a&b  &&  & & && c& d \\
							x :   &  \mathbf{0} & 5 & 4 & 3&2 & \mathbf{2}& * & 1 & \mathbf{0} & 5 & 4 & 3& \mathbf{3} & e & f \\
							x-1 :  & \mathbf{4} & * & *& * &  1 & \mathbf{0} & 5 & 4 & \mathbf{4} & 3 &2& 1 & \mathbf{0} & 5 & 4  
						\end{matrix}
						\label{eq:defect_local_dyn_pf5}\\[5pt]
						&  \begin{matrix}
							\textup{time}: & t & & & & &t-5 & &  &t-8 & &&& t-12 \\ 
							\hline
							x+1 :&&& &&&&a&b  &&  & & && &  \\
							x :   &  \mathbf{0} & 5 & 4 & 3&2 & \mathbf{2}& * & 1 & \mathbf{0} & 5 & 4 & 4& \mathbf{4} &  &  \\
							x-1 :  & \mathbf{4} & * & *& * &  1 & \mathbf{0} & 5 & 4 & \mathbf{4} & 3 &2& 1 & \mathbf{0} &  &   
						\end{matrix}
						\label{eq:defect_local_dyn_pf6}\\[5pt]
						&   \begin{matrix}
							\textup{time}: & t & & & & &t-5 & &  &t-8 & &&&& t-13 \\
							\hline
							x+1 :&&& &&&&a&b  &&  & & && &  & \\
							x :   &  \mathbf{0} & 5 & 4 & 3&2 & \mathbf{2}& * & 1 & \mathbf{0} & 5 & 4 & 3& 3& \mathbf{3} &  &  \\
							x-1 :  & \mathbf{4} & * & *& * &  1 & \mathbf{0} & 5 & 4 & \mathbf{4} & * &*& *&1& \mathbf{0} &  &   
						\end{matrix}
						\label{eq:defect_local_dyn_pf7}\\[5pt]
						&    \begin{matrix}
							\textup{time}: & t & & & & &t-5 & &  &t-8 & &&&& t-13 \\
							\hline
							x+1 :&&& &&&&a&b  &&  & & & & c &d \\
							x :   &  \mathbf{0} & 5 & 4 & 3&2 & \mathbf{2}& * & 1 & \mathbf{0} & 5 & 4 & 3& 2& \mathbf{2} & e & f \\
							x-1 :  & \mathbf{4} & * & *& * &  1 & \mathbf{0} & 5 & 4 & \mathbf{4} & * &*& *&1& \mathbf{0} & 5 & 4  
						\end{matrix}
						\label{eq:defect_local_dyn_pf8}\\[5pt]
						&  \begin{matrix}
							\textup{time}: & t & & & & &t-5 & &  &t-8 & &&&& t-13 \\
							\hline
							x+1 :&&& &&&&a&b  &&  & & && &  & \\
							x :   &  \mathbf{0} & 5 & 4 & 3&2 & \mathbf{2}& * & 1 & \mathbf{0} & 5 & 4 & 4& 3& \mathbf{3} &  &  \\
							x-1 :  & \mathbf{4} & * & *& * &  1 & \mathbf{0} & 5 & 4 & \mathbf{4} & * &*& *&1& \mathbf{0} &  & 
						\end{matrix}
						\label{eq:defect_local_dyn_pf9}
					\end{align}
					
					Note that in these 5 cases we must have $b=0$ or else we get a backtracking contradiction at time $t-8$. Note that  \eqref{eq:defect_local_dyn_pf6}, \eqref{eq:defect_local_dyn_pf7}  and  \eqref{eq:defect_local_dyn_pf9}  have apparent defects in edge $(x,x+1)$. For  \eqref{eq:defect_local_dyn_pf5}  since $b=0$ we have $0\in\{c,d\}$ and we get a defect in either case since $\{e,f\}$ will have color 2 or 3. For  \eqref{eq:defect_local_dyn_pf8}  similary $0\in\{c,d\}$ and we have a defect of color gap 2 either way.  Now we consider the orientation reversed (or flipped) case of \eqref{eq:defect_type25}, that is when $\{X_t(x-1),X_t(x)\}=(0,4)$. This will take care of  \eqref{eq:defect_type25}  and  \eqref{eq:defect_type18}.
					

					\begin{align}
						\begin{matrix}
							\textup{time}:  & t & & & &&t-4 & & & t-7\\ 
							\hline
							x+1 : &&a&b&c&d &&& & &e&f&g&h \\
							x :   &  \mathbf{4} & * & *& * & 1 & \mathbf{0} & 5 & 4 & \mathbf{4} &3 & 2& 1 & \mathbf{0}  \\
							x-1 : & \mathbf{0} &5 & 4 & 3 &  2 & \mathbf{2} & * & 1 &  \mathbf{0} & 5 &4  & *  & * 
						\end{matrix}
						\label{eq:defect_local_dyn_pf10}\\[5pt]
						\begin{matrix}
							\textup{time}:  & t & & & &&t-4 & & & t-7 \\ 
							\hline
							x+1 : &&a&b&c&d &&& &&e&f&g&h \\
							x :   &  \mathbf{4} & * & *& * & 1 & \mathbf{0} & 5 & 4 & \mathbf{4} & *& * & *&1 & \mathbf{0}  \\
							x-1 : & \mathbf{0} &5 & 4 & 3 &  2 & \mathbf{2} & * & 1 &  \mathbf{0}  & 5 &4 & * & * & *
						\end{matrix}
						\label{eq:defect_local_dyn_pf11}
					\end{align}

					For \eqref{eq:defect_local_dyn_pf10} note that $0\in \{a,b,c,d\}$. $a=0$, $b=0$ both yield defects on $(x,x+1)$, $c=0$ gives a backtracking contradiction so we have $d=0$. But then we must have $0\in \{f,g\}$, $f=0$ gives a defect and $g=0$ gives a backtracking contradiction. Similarly with \eqref{eq:defect_local_dyn_pf11} we have $d=0$ and thus $0\in \{f,g\}$ yielding a defect either way.

					\eqref{eq:defect_type21}, \eqref{eq:defect_type23}, \eqref{eq:defect_type24}  and \eqref{eq:defect_type26} have apparent defects, so we only need to consider their orientation reversals. We are going to  leave the reversed version of  \eqref{eq:defect_type23} until the very end since we will use the fact that all remaining sequences (and their reversals) except \eqref{eq:defect_type19} and \eqref{eq:defect_type22} have been excluded.

					First we consider \eqref{eq:defect_type21} reversed. Note that it takes us to $(x,x+1)=(4,0)$ and hence it must be followed by \eqref{eq:defect_type22} (which leads us to $(x,x+1)=(0,3)$ and hence we have two options.

					\begin{align}
						\begin{matrix}
							\textup{time}:  & t & & & &t-4 & & && t-8\\ 
							\hline
							x+1 : &&a&b&c &&& & &e&f&g&h \\
							x :   &  \mathbf{3} & * & *&  1 & \mathbf{0} & 5 & 4 & 3& \mathbf{3} &2 & 1 & \mathbf{0}  \\
							x-1 : & \mathbf{0} &5 & 4 & 4 &   \mathbf{4} & 3 &2 & 1 &  \mathbf{0} & 5 &4  & *  & * 
						\end{matrix}
						\label{eq:defect_local_dyn_pf12}\\[5pt]
						\begin{matrix}
							\textup{time}:  & t & & & &t-4 & & && t-8\\ 
							\hline
							x+1 : &&a&b&c&d &&& & &e&f&g&h \\
							x :   &  \mathbf{3} & * & *&  1 & \mathbf{0} & 5 & 4 & 3& \mathbf{3} &* & *& 1 & \mathbf{0}  \\
							x-1 : & \mathbf{0} &5 & 4 & 4 &   \mathbf{4} & 3 &2 & 1 &  \mathbf{0} & 5 &4  & *  & * 
						\end{matrix}
						\label{eq:defect_local_dyn_pf13}
					\end{align}

					Note for both we must have that $c=0$ or else $0\in \{a,b\}$ which either leads to a defect or a backtracking contradiction.  Then since $x$ doesn't pull $x+1$  between times $t-10$ and  $t-4$ we must have $0\in \{f,g\}$ which in the case of \eqref{eq:defect_local_dyn_pf12} leads to a backtracking contradiction and in the case of \eqref{eq:defect_local_dyn_pf13} leads to a defect.

					Next we deal with flipped \eqref{eq:defect_type24}. It has to be preceded by \eqref{eq:defect_type19} then we have 3 options \eqref{eq:defect_type22}, \eqref{eq:defect_type23} and \eqref{eq:defect_type24}. We have:
					\begin{align}
						& \begin{matrix}
							\textup{time}: & t & & & & &t-5 & &  &t-8 & &&&& \\ 
							\hline
							x+1 :&&a& b&c&0&&  && &  &d &e& &\\
							x :   &  \mathbf{4} & *& * &  *& 1& \mathbf{0}& 5 & 4 & \mathbf{4} & 3 & 2 &  1&\mathbf{0} &  &  \\
							x-1 :  & \mathbf{0} & 5 & 4 & 3&  3 & \mathbf{3} & 2 & 1 & \mathbf{0} & 5&4 &3& \mathbf{3} &  &   
						\end{matrix}
						\label{eq:defect_local_dyn_pf18}
					\end{align}
					
					\begin{align}
						&  \begin{matrix}
							\textup{time}: & t & & & & &t-5 & &  &t-8 & &&&&  \\ 
							\hline
							x+1 :&&a& b&c&0&&  &&  & &d &e& &\\
							x :   &  \mathbf{4} & * & * & *& 1& \mathbf{0}& 5 & 4 & \mathbf{4} & 3 &  2& 1&\mathbf{0} &  &  \\
							x-1 :  & \mathbf{0} & 5 & 4 & 3&  3 & \mathbf{3} & 2 & 1 & \mathbf{0} & 5&4& 4 & \mathbf{4} &  &   
						\end{matrix}
						\label{eq:defect_local_dyn_pf19}\\[5pt]
						&  \begin{matrix}
							\textup{time}: & t & & & & &t-5 & &  &t-8 & &&&&  \\ 
							\hline
							x :   &  \mathbf{4} & * & * & *& 1& \mathbf{0}& 5 & 4 & \mathbf{4} & *&  *& *& 1&\mathbf{0} &  &  \\
							x-1 :  & \mathbf{0} & 5 & 4 & 3&  3 & \mathbf{3} & 2 & 1 & \mathbf{0} & 5&4& 3 & 2&\mathbf{2} &  &   \\
							x-2 :&&& &&0&5&4  &&  & &&d &e& &
						\end{matrix}
						\label{eq:defect_local_dyn_pf20}
					\end{align}

					Note in in \eqref{eq:defect_local_dyn_pf18} \eqref{eq:defect_local_dyn_pf19} $0\notin\{a,b,c\}$ or else we get the usual backtracking contradictions/defects. So then $0\in\{d,e\}$ resulting in a backtracking contradiction. Similarly in \eqref{eq:defect_local_dyn_pf20} we must have $\xi_{t-4}(x-2)=0$ which implies $0\in\{d,e\}$ resulting in a backtracking contradiction.

					We now deal with the flipped \eqref{eq:defect_type26}. Note that it has to be preceded by \eqref{eq:defect_type19} as we excluded the other recurrent defect sequences with with (bold) 3's. Then it we have $\xi_{t-3}(x-1)=0$
					and therefore $0\in\{c,d\}$, since there is one pull from $x-1$. $c=0$ will give a backtracking contradiction so $d=4$. This leaves two possibilities. Either \eqref{eq:defect_type19} is preceded by \eqref{eq:defect_type23} or \eqref{eq:defect_type26}. Note that $\xi_{t-4}(x+1)=0$ or else we get defects or backtracking contradictions as usual. Hence $0\in\{a,b\}$. In the case of \eqref{eq:defect_type26} this yields a defect. In the case of \eqref{eq:defect_type23} a backtracking contradiction.
					\begin{align}
						& \begin{matrix}
							\textup{time}: & t & & & & &t-5 & &  &t-8 & &&&& t-13 \\ 
							\hline
							x+1 :&&& &&0&&  &&  & &a &b& &\\
							x :   &  \mathbf{4} & * & * & *& 1& \mathbf{0}& 5 & 4 & \mathbf{4} & * & * & *& 1&\mathbf{0} &  &  \\
							x-1 :  & \mathbf{0} & 5 & 4 & 4&  3 & \mathbf{3} & 2 & 1 & \mathbf{0} & 5&4& 4 &3& \mathbf{3} &  &   \\
							x-2 :&&& &0&&&  &&  & &c &d& &             
						\end{matrix}
						\label{eq:defect_local_dyn_pf16}
					\end{align}
					
					\begin{align}
						&  \begin{matrix}
							\textup{time}: & t & & & & &t-5 & &  &t-8 & &&& t-12 \\ 
							\hline
							x+1 :&&& &&0&&  &&  & &a &b& &\\
							x :   &  \mathbf{4} & * & * & *& 1& \mathbf{0}& 5 & 4 & \mathbf{4} & 3 &  2& 1&\mathbf{0} &  &  \\
							x-1 :  & \mathbf{0} & 5 & 4 & 4&  3 & \mathbf{3} & 2 & 1 & \mathbf{0} & 5&4& 4 & \mathbf{4} &  &   \\
							x-2 :&&& &0&&&  &&  & &c &d& &             
						\end{matrix}
						\label{eq:defect_local_dyn_pf17}
					\end{align}

					Now we are left with  tackling \eqref{eq:defect_type20}. Note that it can be preceded by \eqref{eq:defect_type19}, \eqref{eq:defect_type20} and \eqref{eq:defect_type21}. This last one, and its reversed version, have been shown to be transient so we must consider these two possibilities.
					
					\begin{align}
						& \begin{matrix}
							\textup{time}: & t & &  & &t-4 & &  &t-8 & &&& t-12 \\ 
							\hline
							x :   &  \mathbf{0} & 5 & 4 & 3& \mathbf{3}& 2 & 1 & \mathbf{0} & 5 & 4 & 3& \mathbf{3} &  &  \\
							x-1 :  & \mathbf{3} & * & * &  1 & \mathbf{0} &5  & 4 & \mathbf{4} & 3 &2& 1 & \mathbf{0} & 5 & 4 &3&\mathbf{3} \\
							x-2 :&&& &0&&&  &&  & &a &b& &             
						\end{matrix}
						\label{eq:defect_local_dyn_pf14}\\[5pt]
						&  \begin{matrix}
							\textup{time}: & t & & & &t-4& & &  &t-7 & &&& t-11 \\ 
							\hline
							x+1 :&&& &&&a&b  &&  &&  && c& d \\
							x :   &  \mathbf{0} & 5 & 4 & 3& \mathbf{3}& * & * &1 & \mathbf{0} & 5 & 4 & 3& \mathbf{3} & e & f \\
							x-1 :  & \mathbf{3} & * & * &  1 & \mathbf{0} & 5 & 4 & 3 &\mathbf{3} & * &*& 1 & \mathbf{0} & 5 & 4  
						\end{matrix}
						\label{eq:defect_local_dyn_pf15}
					\end{align}

					In \eqref{eq:defect_local_dyn_pf14} we have that \eqref{eq:defect_type19} has to be preceded  by  \eqref{eq:defect_type22} since the others with 4 have been excluded. Note here $b=0$ would give a backtracking contradiction. $c=0$ would imply $0\in \{e,f\}$ which also yields a backtracking contradiction at $x-1$ either way. $a=0$ would imply $0\in \{d,e\}$ (note the pull from $x-1$). This would also yield backtracking contradictions at $x-1$.

					Note that in \eqref{eq:defect_local_dyn_pf15} a two step backtracking gives that $b\neq 0$ or else we have a contradiction. So $a=0$ and we have a defect at time $t-8$. By reflection symmetry, this takes care of the flipped version too.

					Finally, we can deal with the reversed \eqref{eq:defect_type23}. Note that, since everything else has been excluded (including the non-reversed version of \eqref{eq:defect_type23}),  backtracking from $\{\xi_{t-4}(x-1),\xi_{t-4}(x)\}=\{4,0\}$  only has the option the recurrent defect obtained by the concatenation of \eqref{eq:defect_type19} and \eqref{eq:defect_type22} (which will be indeed the desired configuration by elimination to be the one). \eqref{eq:defect_type23} in fact acts as a "\textit{orientation switch}":
					\begin{align}
						\begin{matrix}
							\textup{time}: &&&&&&&& t & & & &t-4 & & && t-8 \\ 
							\hline
							x :   &  \mathbf{4} & 3 & 2 & 1 & \mathbf{0} & 5 & 4 & \mathbf{4} & 3& 2 &1 & \mathbf{0} &5&4&3&\mathbf{3}&2&1&\mathbf{0}&5&4&3&\mathbf{3} \\
							x-1 : & \mathbf{0} &5 & 4 & 3 & \mathbf{3} & 2 & 1 &  \mathbf{0}  & 5 &4 & 4 & \mathbf{4} & 3&2&1&\mathbf{0}&5&4&\mathbf{4}&3&2&1&\mathbf{0}\\
						\end{matrix}
						\label{eq:defect_local_dyn_pf21}
					\end{align}

					Its appearance forces the reversal of the concatenated pattern of  \eqref{eq:defect_type19} and \eqref{eq:defect_type22} when evolving forwards. See this illustrated above in \eqref{eq:defect_local_dyn_pf21}. For it to be recurrent the color gap 4 defects must alternate between $\{4,0\}$ and $\{0,4\}$ but this can only happen if the non-reversed \eqref{eq:defect_type23} appears as a switch in between two reversed \eqref{eq:defect_type23}. But it was already shown not to be recurrent. 
				\end{proof}

				\section*{Acknowledgements} 
				HL is partially supported by DMS-2232241.

				\bibliographystyle{amsalpha}  
				\bibliography{mybib}

				\newpage
				
				\appendix

				\section{Fisch's Argument for Undercoupled Regimes}
				
				\hypertarget{appendix}{}

				In this appendix, for the reader's convenience, we provide the proof of Lemma \ref{lem:fixation_suff} (restated in Lemma \ref{lem:FCA_fixation_suff_appendix} below). Our argument follows \cite{bramson1989flux} but differs in that we use the arrow dynamics introduced in Section \ref{sec:arrow_dynamics} in order to keep track of influence between sites. 
				
				\begin{lemma}\label{lem:FCA_fixation_suff_appendix} 
					The following two statements combined give us sufficient conditions for fixation.
					
					\begin{enumerate}
						
						\item $\xi_t$ fixates if
						\begin{align*}
							\lim_{n\rightarrow\infty}\P(T(z)<\infty\quad \text{for some}\quad z<-n)=0
						\end{align*}
						
						\item  For an arrow coming from  $z\in \Z^{-}$  to reach the origin we must have an interval $I$ containing $[z,0]$ such that $\sum_{x\in I}\varphi(x)<0$. More precisely, we have the following event inclusion
						\begin{align*}\{T(z)<\infty \quad \text{for some } z<-n\}\subset \left\{\sum_{x\in I}\varphi(x)< 0\quad \text{for some}\quad I\supset [-n,0]\right\}\end{align*}

					\end{enumerate}
					
				\end{lemma}

				\begin{proof} 
					We begin by showing (1). Define:
					\begin{description}[itemsep=0.1cm]
						\item{(1)} $\tau(j):=$time of the $j$-th excitation at the origin
						\item{(2)} $\alpha(j)=$ arrow label responsible for the $j$-th excitation at the origin. We set $\alpha(j):=\infty$ if $\tau(j):=\infty$. 
						\item{(3)} $A:=\{\tau(j)<\infty \hspace{1mm} \forall j\ge 1\}$
						\item{(4)} $G_n:=\{|\alpha(j)|<n \hspace{1mm} \forall j\ge 1\}$.
					\end{description}

					See proof Proposition \ref{prop:clustering_blockade} to argue that infinitely many excitations of the origin involve infinitely many arrows.  Thus the event $A\cap G_n$ has probability zero.  Since the assumption of (1) and reflection symmetry imply that $G_n$ occurs for some $n$ with probability 1 we have that
					
					$$\P(A)=\P(A\cap (\cup_n G_n))=\P(\cup_n(A\cap G_n))=0$$

					\noindent And since $A^c=\{X_t \hspace{2mm} \text{fixates}\}$, we have fixation.
					
					To prove (2) let us denote $H_n:=\{T(z)<\infty\quad \text{for some}\quad z<-n\}$. On $H_n$, let $\rho$ be the first time an active path from $(-\infty, -n)\times \{0\}$ reaches the origin, and let $m_{-}<-n$ be its initial position. If such a path exists, we define $m_{+}\geq 0$ to be the rightmost source of an active path which reaches the origin before time $\rho$. Otherwise we set $m_{+}=0$. Now focus on the arrows at time 0 for locations  $x\in I=(m_{-},m_{+}]$. On $H_n$, each blockade in I, must at some time before $\rho$ be replaced by a live edge, and these live edges originate in $I$ since arrows do not cross. 
					
					By use of the edge elimination function we conclude
					\begin{align}
						H_n&\subset\left\{\#\text{ live edges  $\geq$  "blocakde mass"  at time 0 on some interval $I\supset[-n, 0]$}\right\}\nonumber\\
						&\subset\left\{\sum_{x\in I}\varphi(x)\leq 0 \hspace{1mm}\text{for some interval $I\supset[-n, 0]$}\right\}\nonumber\\
						&\subset\left\{\sum_{x=-l}^m\varphi(x)\leq 0 \quad \text{for some $l\geq n$ and $m\geq 0$}\right\}
					\end{align}

					\vspace{2mm}
					To see these inclusions, suppose (WLOG) we have a right arrow on the leftmost edge $(z,z+1)$. For it to reach the origin, all the blockade mass on I must be cleared. The first inclusion is trivial since each blockade needs at least one arrow to be destroyed.  As defined in (\ref{eq1}), the edge elimination function counts free arrows as negative (regardless of the direction they are pointing, a conservative estimate). More precisely, the $|d\eta|-2r$ term accounts for the fact that each blockade has "capacity" $|d\eta|-r$ and releases $r$ arrows when destroyed. This explains the second inclusion. The last one is trivial. 
				\end{proof}

			\end{document}